\documentclass[12pt,a4paper]{article}
\usepackage{amsmath}
\usepackage{amssymb}\usepackage{bbm}
\usepackage{booktabs} 
\usepackage{threeparttable}
\usepackage{multirow}
\usepackage{graphicx}
\usepackage[figuresright]{rotating}
\usepackage{subfigure}
\usepackage{latexsym,amsmath,amsfonts,amscd,amssymb}
\usepackage[colorlinks,linkcolor=blue,citecolor=blue]{hyperref}%
\usepackage[noend]{algpseudocode}
\usepackage{algorithmicx,algorithm}

\usepackage{epsfig}
\usepackage{epstopdf}
\usepackage{changebar}
\usepackage{indentfirst}
\usepackage{verbatim}
\usepackage{color}
\usepackage{colordvi}
\usepackage{appendix}
\usepackage{cite}
\usepackage[numbers,sort]{natbib}

\newcommand{\dx}{\Delta x}
\newcommand{\dy}{\Delta y}
\newcommand{\dt}{\Delta t}

\newcommand{\bfG}{\mathbf{G}}

\usepackage{indentfirst}
\usepackage{verbatim}
\usepackage{tikz}
\usepackage{array}

\usepackage{amsmath}
\usepackage{theorem}
\theorembodyfont{\normalfont}
\newtheorem{theorem}{Theorem}[section]

\newtheorem{rem}{Remark}[section]

\newtheorem{example}{Example}[section]

\usepackage{pgflibraryarrows}
\usepackage{pgflibrarysnakes}
\usepackage{booktabs,makecell,multirow}

\numberwithin{equation}{section}
\numberwithin{figure}{section}
\numberwithin{table}{section}

\usepackage{soul}
\usepackage{caption}
\begin{document}

\captionsetup[figure]{labelfont={bf},name={Fig.},labelsep=period}

\baselineskip=1.5pc
\begin{center}
	{\Large {\bf A moment-based Hermite WENO scheme with unified stencils for hyperbolic conservation laws}}\footnote{{The research was partially supported by  National Key R$\&$D Program of China [Grant Number 2022YFA1004500]}.
	}
\end{center}

\begin{center}
	Chuan Fan\footnote{Department of  Mathematics, Southern University of Science and Technology, Shenzhen, Guangdong 518055, P.R. China. E-mail: fanc@sustech.edu.cn.},
	Jianxian Qiu\footnote{School of Mathematical Sciences and Fujian Provincial Key Laboratory of Mathematical Modeling and High-Performance Scientific Computing, Xiamen University, Xiamen, Fujian 361005, P.R. China. E-mail: jxqiu@xmu.edu.cn.}, and
	Zhuang Zhao\footnote{Corresponding author. School of Mathematical Sciences,  Xiamen University, Xiamen, Fujian 361005, P.R. China. E-mail: zzhao@xmu.edu.cn.}
\end{center}

\vspace{.05in}
\centerline
{\bf Abstract\ }
\bigskip

In this paper, a fifth-order moment-based Hermite weighted essentially non-oscillatory scheme with unified stencils (termed as HWENO-U) is proposed for hyperbolic conservation laws. The main idea of the HWENO-U scheme is to modify the first-order moment by a HWENO limiter only in the time discretizations  using the same information of spatial reconstructions, in which the limiter not only overcomes spurious oscillations well, but also ensures the stability of the fully-discrete scheme. For the HWENO reconstructions, a new scale-invariant nonlinear weight is designed by incorporating only the integral average values of the solution, which keeps all properties of the original one while is more robust for simulating challenging problems with sharp scale variations. Compared with previous HWENO schemes, the advantages of the HWENO-U scheme are: (1) a simpler implemented process involving  only a single HWENO reconstruction applied throughout the entire procedures  without any modifications for the governing equations; (2) increased efficiency by utilizing the same candidate stencils, reconstructed polynomials, and linear and nonlinear weights in both the HWENO limiter and spatial reconstructions; (3) reduced problem-specific dependencies and improved rationality,  as the nonlinear weights are identical for the function $u$ and its non-zero multiple $\zeta u$. Besides, the proposed scheme retains the advantages of previous HWENO schemes, including compact reconstructed stencils and the utilization of artificial linear weights.  Extensive benchmarks are carried out to validate the accuracy, efficiency, resolution, and robustness of the proposed scheme.

\vspace{.2in}
\vfill
\vfill {\bf Key Words:} hyperbolic conservation laws, Hermite WENO scheme, unified stencils, limiter, scale-invariant nonlinear weight

{\bf AMS(MOS) subject classification}: 65M60, 35L65


\newpage
\baselineskip=2pc
\section{Introduction}

In this paper, we construct a fifth-order  Hermite weighted essentially non-oscillatory scheme with unified {candidate} stencils (termed as HWENO-U) for hyperbolic conservation laws, where both the zeroth- and first-order moments are evolved in time and used in spatial reconstructions. Compared with other moment-based HWENO schemes \cite{ZhaoChenQiu, ZhaoZQiuJX_ArtificalLinearWeight_HWENO2020, LiShuQiu2022fvMRHWENO, TaoLiQiu, ZhaoZQiuJX2023OFHWENO, DBTM}, the HWENO-U scheme adds a  high order modification for the first-order moments in time discretizations by using the same information of spatial reconstructions, which is {simpler} and more efficient for using the same reconstructed polynomials, smooth indicators, linear and nonlinear weights {in the entire procedures}. HWENO schemes are constructed on the basis of weighted essentially non-oscillatory (WENO) schemes, and WENO schemes have been widely applied for hyperbolic conservation laws in the past three decades. The first WENO scheme was proposed by Liu {et al.} \cite{loc} in 1994, where they combined all candidate stencils of essentially non-oscillatory (ENO) schemes \cite{ho,h1,heoc} to achieve a third-order accuracy in the finite volume version. Next, Jiang and Shu developed a fifth-order finite difference WENO scheme  \cite{js} in 1996, in which they gave a general definition for the smoothness indicators and nonlinear weights, and the fifth-order finite volume WENO scheme was presented by Shu \cite{s2} in 1998. After that, WENO schemes have been further developed in \cite{hs,lpr,cd1,ZStetra,ccd,ZQd,bgsAO,ZSWENOMR}, and a recent review can be {found} in \cite{shu}.

The {fundamental} difference between WENO and HWENO schemes is spatial discretizations, where WENO schemes only use the information of solutions, but HWENO schemes can use additional information in each cell, such as the derivatives or first-order moments of solutions. Hence, HWENO schemes can use more compact stencils than WENO schemes on the same order accuracy, resulting in {more minor} numerical errors in smooth cases and {fewer} transition points near discontinuities based on the comparisons in \cite{ZhangZhao}. However, HWENO schemes are less robust than WENO schemes as the derivatives or first-order moments may become quite large near discontinuities. For example, using the same thought of the first one-dimensional HWENO scheme \cite{QSHw1}, the first two-dimensional HWENO scheme \cite{QSHw} gave poor resolutions for the double Mach and forward step problems, though this drawback was solved later in \cite{ZQHW} by using more complicated techniques to reconstruct the derivative terms. {The common point of the {HWENO schemes} } \cite{QSHw1,QSHw,ZQHW} is to use different {stencils or techniques} in the discretization of the governing and derived equations by avoiding discontinuities, {which} also has been used in the subsequent HWENO schemes \cite{LLQ,ZA,TaoLiQiu,CaiZhangQiu2016PPFVHWENOEuler,LiMRHW1}.
However, reducing the use of derivatives or first-order moments alone is not sufficient to control oscillations effectively. For instance, additional techniques such as positive-preserving limiters and {a} smaller time step are required in the first finite difference HWENO scheme \cite{LLQ}. Furthermore, the selection of optimal stencils and approximated methods often heavily relies on numerical experiences.

To enhance the robustness of HWENO schemes, Zhao et al. \cite{ZhaoChenQiu} proposed an alternative approach to control the derivatives or first-order moments as limiters in the discontinuous Galerkin (DG) method \cite{DG2}, which can effectively overcome oscillations even with a normal time step. This moment-based HWENO scheme also can be viewed as a $P_1P_M$ method, as defined by Dumbser et al. \cite{DBTM}. The key feature of the HWENO scheme \cite{ZhaoChenQiu} is the separation of limiters and spatial reconstructions into two distinct parts, while the limiters in DG methods \cite{Zhongshu, LuoBaum, Zhulimiter} and the spatial reconstructions in WENO schemes \cite{cd1,ccd,ZQd,bgsAO} have been extensively studied over the past three decades. Consequently, constructing HWENO schemes {with the proposed framework}  \cite{ZhaoChenQiu} benefits from the wealth of mature references in these two fields, reducing reliance on numerical experiences. Later, by modifying the first-order moments in advance as \cite{ZhaoChenQiu}, the modified HWENO scheme with artificial linear weights \cite{ZhaoZQiuJX_ArtificalLinearWeight_HWENO2020}, the positivity-preserving HWENO scheme \cite{FanZhangQiu2021PPhybridHWENO}, the Hermite TENO scheme \cite{weHTENO}, the multi-resolution HWENO scheme \cite{LiShuQiu2022fvMRHWENO}, and the finite difference HWENO combined with limiter scheme \cite{ZhangZhao} have been developed to solve hyperbolic conservation laws. However, {the proposed framework} \cite{ZhaoChenQiu} still utilizes two sets of stencils as the first HWENO scheme \cite{QSHw1}, and using two sets of stencils means repetitive algorithms, and double or triple computational costs. Recently,  Zhao and Qiu \cite{ZhaoZQiuJX2023OFHWENO} designed a sixth-order  HWENO scheme by introducing damping terms in the first-order moment equations as the oscillation-free DG methods \cite{OFDG1,OFDG2}. This approach allows for the use of unified stencils in spatial reconstructions,  which is easier to implement and have higher efficiency. However, the presence of damping terms in \cite{ZhaoZQiuJX2023OFHWENO} has  significant impacts on the stability, particularly when simulating strong shocks and extreme problems with highly stiff damping terms, which leads to a small time step restriction and requires the (modified) exponential Runge-Kutta (RK)  time discretization \cite{HuangShuexp}.

In this paper, we mainly focus on developing a practical moment-based HWENO scheme with unified candidate stencils. Based upon previous studies of the HWENO schemes \cite{ZhaoZQiuJX_ArtificalLinearWeight_HWENO2020, ZhaoChenQiu, LiShuQiu2022fvMRHWENO, FanZhangQiu2021PPhybridHWENO, ZhaoZQiuJX2023OFHWENO, TaoLiQiu},  it has been observed that the first-order moments tend to become large near discontinuities, which potentially impacts the robustness of HWENO schemes. To address this issue, various approaches have been introduced in the aforementioned HWENO schemes, such as reducing the utilization of the first-order moments optimally, controlling  the first-order moments near discontinuities before spatial reconstructions, or introducing damping terms in the first-order moment equations. These schemes share a common characteristic where all first-order moments are utilized in the spatial reconstructions. The main reason is that omitting the utilization of the first-order moment on the central cell will lead to instabilities for HWENO schemes based on our mathematical analysis. More intuitively, we take the one-dimensional case in the target cell $I_i$ as an example. If we directly use the values of the zeroth-order moments $\{\overline{u}_{i-1},\overline{u}_{i},\overline{u}_{i+1}\}$ and the first-order moments $\{\overline{v}_{i-1},\overline{v}_{i+1} \}$ to discretize the space, the fully-discrete scheme will be unstable by using the forward Euler or third-order SSP RK time discretization \cite{s2}, proved in Theorem \ref{sec2:thm1}. Taking into account  the symmetry of the stencils and the formulation of the  Lax-Friedrichs scheme \cite{LaxFridrichs1954}:
\begin{equation}\label{LaxF} \frac{\overline{u}^{n+1}_i-\left(\frac{\overline{u}^n_{i-1}+\overline{u}^n_{i+1}}{2}\right)}{\dt}+\frac{\overline{u}^{n}_{i+1}-\overline{u}^{n}_{i-1}}{2\dx}=0.
\end{equation}
It is well known that the scheme $\frac{\overline{u}^{n+1}_i-\overline{u}^{n}_i}{\dt}+\frac{\overline{u}^{n}_{i+1}-\overline{u}^{n}_{i-1}}{2\dx}=0$
is  unstable, but the Lax-Friedrichs scheme  \eqref{LaxF} is stable as $\overline{u}^{n}_i$  is modified by $\frac{\overline{u}^n_{i-1}+\overline{u}^n_{i+1}}{2}$ in the time discretization. Inspired by this point, we also introduce a modification for the first-order moment in time discretizations using the information provided by spatial discretizations, and the proposed scheme is proved to be stable through analyses in Subsection \ref{sec2:stabilityHWENO} using the Fourier method. {Besides, it is worth noting that this modification step, which utilizes the same information as these for spatial discretizations, plays a vital role in the adoption of unified stencils throughout the entire procedures.} To overcome spurious oscillations near discontinuities, we use the HWENO method in the modified HWENO (HWENO-M) scheme \cite{ZhaoZQiuJX_ArtificalLinearWeight_HWENO2020} to modify the first-order moments and {perform spatial discretizations, where} the linear weights can be any positive numbers as long as their sum is one. Differently, the modification and spatial discretizations are combined into a single step for they use the same information, resulting in unified candidate stencils in the HWENO-U scheme, which simplifies the implementation process of \cite{ZhaoZQiuJX_ArtificalLinearWeight_HWENO2020} and enhances the computational efficiency. Furthermore, we also improve the nonlinear weights in the HWENO-M scheme to make them scale-invariant. For the reasonability,  the function $u$ and its non-zero multiple $\zeta u$  should have the same nonlinear weights on the same cells. Conversely,  the nonlinear weights in the HWENO-M scheme lose this basic property. Although this scale-dependent nonlinear weight has no obvious {differences} on the simulations of benchmark tests shown in the various subsquen WENO schemes \cite{ZQd,ZhuJunQiuJX2017fvWENOZQ, tezq, trzq},  the results may generate oscillations in simulating large-scale problems \cite{ChenWu,CaiQiuWu}, and the similar phenomenon also occurs in simulating extreme problems based on our numerical experiments. To inherit the advantages of the nonlinear weights in \cite{ZQd,ZhuJunQiuJX2017fvWENOZQ, tezq, trzq, ZhaoZQiuJX_ArtificalLinearWeight_HWENO2020} and make them scale-invariant, we bring the integral average values of solutions into the original definition, which also can be viewed as a normalization procedure for the nonlinear weights. This minor modification has no impacts on the accuracy firstly, and it also can enhance resolutions and is more {robust} for simulating challenging problems with sharp scale variations. In short,  the HWENO-U scheme uses unified {stencils} in the entire procedures, which avoids repetitive algorithms and enhances computational efficiencies. {Furthermore, the HWENO-U scheme has the capability to simulate extreme problems by directly incorporating a  positive-preserving technique from \cite{FanZhangQiu2021PPhybridHWENO, ZhangS}, which is simpler and more practical compared to the other HWENO scheme with unified stencils \cite{ZhaoZQiuJX2023OFHWENO}, since the proposed scheme avoids the introduction of additional parameters and stiff terms.} Besides, the designed scale-invariant nonlinear weight is more reasonable and {robust} in numerical simulations. These advantages will be demonstrated in the next algorithm descriptions and numerical tests.

The paper is organized {as follows: Section 2 presents} the detailed implementation of the HWENO-U scheme in one- and two-dimensional cases, and provides a stability analysis for the linear scheme using the Fourier method. In Section 3, extensive benchmarks  are conducted to illustrate the numerical accuracy, high resolution, and robustness of the proposed scheme. Finally, concluding remarks are given in Section 4.

\section{Description of HWENO-U scheme}\label{sec2:HWENO-U}

{This section contains three subsection. {In the first and second subsections, we present} the detailed procedures of the moment-based fifth-order HWENO-U scheme in the one- and two-dimensional cases, respectively, in which the high order HWENO modification for the first-order moments and the spatial reconstructions use the same information, such as candidate stencils, reconstructed polynomials, smoothness indicators, linear and nonlinear weights. Remarkably, incorporating the modification for the first-order moments in the time discrete stage is essential to ensure the stability of the HWENO-U scheme, therefore, we give a stability analysis in the last subsection to illustrate it.}

\subsection{One-dimensional case}\label{sec2:1dcase}
Consider one-dimensional scalar hyperbolic conservation laws
\begin{equation}\label{sec2:1dHCLS}
	\begin{cases}
		u_t + f(u)_x=0,
		\\
		u_0(x)=u(x,0).
	\end{cases}
\end{equation}
For simplicity, we consider a uniform partition of a given domain $[a,b]$, $a=x_{\frac12}<x_{\frac32}<\cdots<x_{N_x+\frac12}=b$. Let $I_{i}=[x_{i-\frac12},x_{i+\frac12}]$ denote a computational cell with its length $\dx=x_{i+\frac12}-x_{i-\frac12}$ and its center $x_i=\frac12(x_{i-\frac12}+x_{i+\frac12})$. By multiplying equation \eqref{sec2:1dHCLS} with a test function $\phi(x) \in \text{span}\left\{\frac{1}{\dx}, \frac{x - x_i}{(\dx)^2}\right\}$, integrating over interval $I_i$ and using integration by parts, we have
\begin{equation}\label{sec2:WeakForm1d}
	\left\{
	\begin{aligned}
		&\frac{{\rm d} \bar{u}_i(t)}{{\rm d} t}=-\frac 1 {\dx} \left ( f(u(x_{i+\frac12},t))-f(u(x_{i+\frac12},t))\right ),\\
		&\frac{{\rm d} \bar{v}_i(t)}{{\rm d} t}=
		- \frac 1 {2\dx} \left ( f(u(x_{i-\frac12},t)) +f(u(x_{i+\frac12},t))\right ) +\frac 1 {(\dx)^2} \int_{I_i}f(u)\mathrm{d}x,
	\end{aligned}
	\right.
\end{equation}
where $ \bar{u}_i(t) =\frac{1}{\dx}\int_{I_i}u(x,t)dx $ and $\bar{v}_i(t) = \frac{1}{\dx} \int_{I_i} u(x,t) \frac{x-x_i}{\dx}dx$ are the zeroth-order moment (cell-average) and first-order moment in the cell $I_i$, respectively.

Let $\{\hat{x}^G_i\}^4_{G=1}$ denote four Gauss-Lobatto points in a cell $I_i$ with the corresponding weights $\{\hat{\omega}_G\}^4_{G=1}$  on the interval $[-\frac12,\frac12]$. The value of the flux function  $f(u(x_{i+\frac12},t))$ is approximated by a high order Lax-Friedrichs numerical flux $\hat{f}_{i+\frac12}$ and the integral term $\int_{I_{i}}f(u)dx$ is approximated by a 4-point Legendre Gauss-Lobatto quadrature formula. Consequently, a conservative semi-discrete scheme is defined as
\begin{equation}\label{sec2:1dSemi_HWENO}
	\left\{
	\begin{aligned}
		&\frac{{\rm d} \bar{u}_i(t)}{{\rm d} t}=-\frac{1}{\dx}\;(\hat{f}_{i+\frac12}-\hat{f}_{i-\frac12}) \triangleq \mathcal{F}^1_i(	\bar{u},	\bar{v}), 
		\\
		&\frac{{\rm d} \bar{v}_i(t)}{{\rm d} t}=-\frac{1}{2\dx}(\hat{f}_{i-\frac12}+\hat{f}_{i+\frac12})+\frac{1}{\dx} \sum\limits_{G=1}^4 \hat{\omega}_G f(u(\hat{x}^G_i,t)) \triangleq \mathcal{F}^2_i(	\bar{u},	\bar{v} ),
	\end{aligned}
	\right.
\end{equation}
where $\mathcal{F}^1_i(	\bar{u},	\bar{v} )$ and  $\mathcal{F}^2_i(	\bar{u},	\bar{v} )$ are the  right-hand terms, and the numerical flux $\hat{f}_{i+\frac12}$ is defined as
\begin{equation}
	\hat{f}_{i+\frac12}=\frac12\left[f(u^-_{i+\frac12}) + f(u^+_{i+\frac12})  - {\alpha} ( u^+_{i+\frac12} - u^-_{i+\frac12} ) \right],
\end{equation}
with $\alpha=\max\limits_u|f'(u)|$. The superscribes $``-"$ and $``+"$ of $u^\pm_{i+\frac12}$ represent the left-hand and right-hand limits at the interface $x_{i+\frac12}$, respectively. The Gauss-Lobatto points $\{\hat{x}^G_i\}^4_{G=1}$ are
\[
\hat{x}^1_i=x_{i-\frac12},~\hat{x}^2_i=x_{i-\frac{\sqrt{5}}{10}},~\hat{x}^3_i=x_{i+\frac{\sqrt{5}}{10}},~\hat{x}^4_i=x_{i+\frac12},
\] with $x_{i+\ell}=x_i+\ell\dx$, the normalized weights $\hat{\omega}_1=\hat{\omega}_4=\frac{1}{12}$ and $\hat{\omega}_2=\hat{\omega}_3=\frac{5}{12}$.

The equations \eqref{sec2:1dSemi_HWENO} also are the semi-discrete form of the {$\mathbbm{P}^1(I_i)$ DG finite element} method, but for  moment-based HWENO schemes, a Hermite reconstruction is used to approximate the values  $u^\pm_{i\mp\frac12}$ and $u_{i\pm\frac{\sqrt{5}}{10}}$ based on the zeroth- and first-order moments in the cells $\{I_{i-1}, I_{i}, I_{i+1}\}$. In the following, we will outline the detailed steps of 1D HWENO-U scheme,  based on the set of  values $\{\bar{u}_{i-1}, \bar{u}_{i}, \bar{u}_{i+1}, \bar{v}_{i-1}, \bar{v}_{i+1}\}$.

\textbf{Step 1.} Reconstruct a quartic polynomial $p_0(x)$ and two linear polynomials $\{p_m(x)\}_{m=1}^2$.

Firstly, we consider a large stencil $S_0=\{I_{i-1},I_i,I_{i+1}\}$ and two small stencils $S_1=\{I_{i-1},I_i\}$, $S_2=\{I_{i},I_{i+1}\}$. {A quartic polynomial $p_0(x)$ is reconstructed by a Hermite reconstruction on $S_0$, satisfying}
\begin{equation}
		\frac{1}{\dx}\int_{I_k}p_0(x)\mathrm{d}x=\bar{u}_k,~k=i-1,i,i+1,~~
		\frac{1}{\dx}\int_{I_k}p_0(x)\frac{x-x_k}{\dx}\mathrm{d}x=\bar{v}_k,~k=i-1,i+1.
\end{equation}
Two linear polynomials {$\{p_m(x)\}_{m=1}^2$ are obtained} by a linear reconstruction based on $S_1$ and $S_2$, respectively, having
\begin{equation}
		\frac{1}{\dx}\int_{I_k} p_1(x) \mathrm{d}x=\bar{u}_k,~k=i-1,i;~~
		\frac{1}{\dx}\int_{I_k} p_2(x) \mathrm{d}x=\bar{u}_k,~k=i,i+1.
\end{equation}
{Then we rewrite $p_0(x)$ as
\begin{equation}\label{sec21:linearweights}
	{p_0(x)=\gamma_0\left(\frac{1}{\gamma_0}p_0(x) - \frac{\gamma_1}{\gamma_0}p_1(x) - \frac{\gamma_2}{\gamma_0}p_2(x)\right) + \gamma_1p_1(x)  + \gamma_2p_2(x), \gamma_0 \neq 0. }
\end{equation}
To ensure the next WENO procedure stable, $\{\gamma_m\}^2_{m=0}$ are  positive with {$\sum_{m=0}^2 \gamma_m =1$}. }

\textbf{Step 2.} Compute smoothness indicators $\{\beta_m\}_{m=0}^2$ to measure the level of smoothness for the functions $\{p_m(x)\}_{m=0}^2$ in the cell $I_i$, which is  defined as in the classical WENO scheme \cite{js}, satisfying
\begin{equation}\label{sec2:def_beta}
		\beta_m=\sum\limits_{l=1}^{r} \int_{I_i}\dx^{2l-1}(\frac{{\rm d}^l p_m(x)}{{\rm d} x^l})^2 \mathrm{d}x,m=0,1,2,
\end{equation}
where $r$ is the degree of the polynomials $p_m(x)$.  Let $p_m(x)=\sum\limits_{l=0}^r c_{m,l} (\frac{x-x_i}{\dx})^l$, the explicit expressions of the smoothness indicators are
\begin{equation}\label{sec2:beta1d}
	\left \{
	\begin{aligned}
		\beta_0=&(c_{0,1} + \frac{1}{4} c_{0,4})^2 + \frac{13}{3}(c_{0,2} + \frac{63}{130}c_{0,4})^2 + \frac{781}{20} c_{0,3}^2 + \frac{1421461}{2275}c_{0,4}^2,\\
		\beta_m=&c_{m,1}^2, m=1,2,
	\end{aligned}
	\right.
\end{equation}
where the coefficients of the polynomials $\{p_m(x)\}^2_{m=0}$ are listed in Appendix \ref{sec:Appendix1}.

\textbf{Step 3.} Compute nonlinear weights based on linear weights and smoothness indicators.
As in the WENO scheme of Zhu and Qiu (WENO-ZQ) \cite{ZhuJunQiuJX2017fvWENOZQ}, we also introduce a new parameter $\tau$ to measure the absolute difference between $\beta_{0}$, $\beta_{1}$ and $\beta_{2}$ as
\begin{equation}
	\tau=\left(\frac{|\beta_{0}-\beta_{1}|+|\beta_{0}-\beta_{2}|}{2}\right)^2.
\end{equation}
Differently, we add the integral average values  of solutions  into the original definition of the nonlinear weights, that is  $u_{ave}$, having
\begin{equation}\label{sec2:non_weight_1d}
	\omega_m=\frac{\widetilde\omega_m}{\widetilde\omega_{0}+\widetilde\omega_{1}+\widetilde\omega_{2}},
	\ \mbox{with} \ \widetilde\omega_{m}=\gamma_{m}(1+\frac{\tau}{\beta_{m} u_{ave}+\varepsilon}), \ m=0,1,2,
\end{equation}
where $u_{ave} =(\frac{|\overline u_{i-1}|+|\overline u_{i}|+|\overline u_{i+1}|}{3})^2$,
and $\varepsilon$ is set as $10^{-10}$ for avoiding zero denominator. {It is worth mentioning that the nonlinear weights \eqref{sec2:non_weight_1d} are scale-invariant, since the function $u$ and its non-zero multiple $\zeta u$ have the same $\widetilde\omega_{m}$ on the same cells.}

\begin{rem}%
{The original  $\widetilde\omega_{m}$ is   defined as $\gamma_{m}(1+\frac{\tau}{\beta_{m} +\varepsilon})$ for the WENO and HWENO schemes with artificial linear weights \cite{ZhuJunQiuJX2017fvWENOZQ, ZhaoZQiuJX_ArtificalLinearWeight_HWENO2020}, which is also a special case of the formula \eqref{sec2:non_weight_1d} when $u_{ave}=1$. However, the original nonlinear weight in \cite{ZhuJunQiuJX2017fvWENOZQ, ZhaoZQiuJX_ArtificalLinearWeight_HWENO2020} depends on the scale of functions as its $\widetilde\omega_{m}$ is $\gamma_{m}(1+\frac{\zeta^2\tau}{\beta_{m} +\varepsilon})$ for the function $\zeta u$, and this scale-dependent  nonlinear weight has obvious impacts on the simulations of problems with different scales,} as demonstrated in Examples \ref{Example:Lax1D} and \ref{Example:Sedov2D}. Multiplying by $u_{ave}$ in \eqref{sec2:non_weight_1d} is actually a normalization procedure for the nonlinear weights, and it will not destroy the properties of original one, such as the accuracy and resolution, as $u_{ave}$ is of order $\mathcal{O}(1)$. Besides, the formula \eqref{sec2:non_weight_1d} is still suitable when $u_{ave}$ = 0, {that is $\bar{{u}}_{i-1}=\bar{{u}}_{i}=\bar{{u}}_{i+1}=0$, where the solution is smooth in the target cell and the nonlinear weights will degenerate to linear weights simultaneously.}
\end{rem}

Finally, through replacing a part of linear weights  in \eqref{sec21:linearweights} by the nonlinear weights \eqref{sec2:non_weight_1d}, we obtain a nonlinear HWENO reconstructed polynomial $u_i(x)$ for  $u(x)$.  Additionally, a high order modification $\hat{v}_{i}$ for the first-order moment $\bar{v}_{i}$ is obtained using the same $p_m(x)$, $\gamma_m$, and $\omega_m$ simultaneously, which is only used in the next time discretization, having
\begin{equation}\label{sec2:1dHWENORecon}
	\begin{cases}
		u_i(x)=\omega_0\left(\frac{1}{\gamma_0}p_0(x) - \frac{\gamma_1}{\gamma_0}p_1(x) - \frac{\gamma_2}{\gamma_0}p_2(x)\right) + \omega_1p_1(x)  + \omega_2p_2(x),
		\\
		\hat{v}_i=\frac{1}{\Delta x} \int_{I_{i}} u_i(x)\frac{x-x_{i}} {\Delta x}dx=\omega_0\left(\frac{1}{\gamma_0}q_0 - \frac{\gamma_1}{\gamma_0}q_1 - \frac{\gamma_2}{\gamma_0}q_2\right) + \omega_1q_1  + \omega_2q_2,
	\end{cases}
\end{equation}
where $q_m=\frac{1}{\Delta x} \int_{I_{i}} p_m(x)\frac{x-x_{i}} {\Delta x}dx$. Then, the required Gauss-Lobatto point values are evaluated by
\begin{equation*}
	u^+_{i-\frac12} = u_i(x_{i-\frac12}),~u_{i\pm\frac{\sqrt{5}}{10}} = u_i(x_{i\pm\frac{\sqrt{5}}{10}}),~u^-_{i+\frac12} = u_i(x_{i+\frac12}).
\end{equation*}

\textbf{Step 4.} Time discretizations for the semi-discrete scheme \eqref{sec2:1dSemi_HWENO}.

To construct a stable scheme, we modify  the first-order moments in time discretizations as the
Lax-Friedrichs scheme \cite{LaxFridrichs1954} on the basis of  the third-order SSP RK method \cite{s2}, then, the fully-discrete one-dimensional HWENO-U scheme for Eq. \eqref{sec2:1dSemi_HWENO} is written as
\begin{equation}\label{sec2:1dRK3_HWENO}
	\left\{
	\begin{aligned}
		\begin{bmatrix}
			\bar{u}^{(1)}_i \\ \bar{v}^{(1)}_i
		\end{bmatrix} &=~~
		\begin{bmatrix}
			\bar{u}^{n}_i \\ {\hat{v}^{n}_i}
		\end{bmatrix}
		+\dt \begin{bmatrix}
			\mathcal{F}^1_i(\bar{u}^{n},	\bar{v}^{n}) \\
			\mathcal{F}^2_i(\bar{u}^{n},	\bar{v}^{n})
		\end{bmatrix},
		\\
		\begin{bmatrix}
			\bar{u}^{(2)}_i \\ \bar{v}^{(2)}_i
		\end{bmatrix} &=\frac34
		\begin{bmatrix}
			\bar{u}^{n}_i \\  {\hat{v}^{n}_i}
		\end{bmatrix}
		+\frac14( \begin{bmatrix}
			\bar{u}^{(1)}_i \\ {\hat{v}^{(1)}_i}
		\end{bmatrix}
		+\dt \begin{bmatrix}
			\mathcal{F}^1_i(\bar{u}^{(1)},	\bar{v}^{(1)}) \\
			\mathcal{F}^2_i(\bar{u}^{(1)},	\bar{v}^{(1)})
		\end{bmatrix}  ),
		\\
		\begin{bmatrix}
			\bar{u}^{n+1}_i \\ \bar{v}^{n+1}_i
		\end{bmatrix}&=\frac13
		\begin{bmatrix}
			\bar{u}^{n}_i \\ {\hat{v}^{n}_i}
		\end{bmatrix}
		+\frac23( \begin{bmatrix}
			\bar{u}^{(2)}_i \\ {\hat{v}^{(2)}_i}
		\end{bmatrix}
		+\dt \begin{bmatrix}
			\mathcal{F}^1_i(\bar{u}^{(2)},	\bar{v}^{(2)}) \\
			\mathcal{F}^2_i(\bar{u}^{(2)},	\bar{v}^{(2)})
		\end{bmatrix}  ),
	\end{aligned}
	\right.
\end{equation}
where $\hat{v}^{n}_i$, $\hat{v}^{(1)}_i$ and $\hat{v}^{(2)}_i$ represent the high order modification of  $\bar{v}^{n}_i$, $\bar{v}^{(1)}_i$ and $\bar{v}^{(2)}_i$, respectively,  obtained  by the formula \eqref{sec2:1dHWENORecon}. This modified  time discretization has also been used in the finite difference HWENO schemes \cite{FanZhaoXQ,ZhangZhao}, but these two schemes are unable to utilize unified stencils in {both the spatial discretizations and the modification of $v_i$}.

\subsection{Two-dimensional case}\label{sec2:2dcase}
Consider two-dimensional scalar hyperbolic conservation laws
\begin{equation}\label{sec2:2dHCLS}
	\begin{cases}
		u_t+ f(u)_x+g(u)_y=0, \\
		u(x,y,0)=u_0(x,y). \\
	\end{cases}
\end{equation}
For simplicity, we also consider a uniform partition of a given domain $[a,b]\times[c,d]$ with computational cells $I_{i,j}=[x_{i-\frac12},x_{i+\frac12}]\times [y_{j-\frac12},y_{j+\frac12}]$ for $i=1,\ldots,N_x,j=1,\ldots,N_y$. The mesh sizes are $\dx=x_{i+\frac12}-x_{i-\frac12}$ and  $\dy=y_{j+\frac12}-y_{j-\frac12}$, and $(x_i,y_j)$ is the center of the cell $I_{i,j}$ with $x_i=\frac12(x_{i-\frac12}+x_{i+\frac12})$ and $y_j= \frac12(y_{j-\frac12}+y_{j+\frac12})$. Define $I_i= [x_{i-\frac12},x_{i+\frac12}]$ and $I_j=[y_{j-\frac12},y_{j+\frac12}]$.
After multiplying the equation \eqref{sec2:2dHCLS} by a test function $\phi(x,y)\in span\{\frac{1}{\dx\dy},\frac{x-x_i}{(\dx)^2\dy},\frac{y-y_j}{\dx(\dy)^2}\}$, integrating over the cell $I_{i,j}$, and using the integration by parts, we have
\begin{equation}\label{sec2:WeakForm2d}
	\left\{
	\begin{aligned}
		\frac{{\rm d} \bar u_{i,j}(t)}{{\rm d} t}&=
		-\frac{1} {\dx \dy}  \int_{I_j}\left[ f\left(u(x_{i+\frac12},y,t)\right)- f\left(u(x_{i-\frac12},y,t)\right)\right]\mathrm{d}y
		\\&
		\quad -\frac{1} {\dx\dy}  \int_{I_i}\left[ g\left(u(x,y_{j+\frac12},t)\right)- g\left(u(x,y_{j-\frac12},t)\right)\right]\mathrm{d}x,
		\\		
		\frac{{\rm d} \bar v_{i,j}(t)}{{\rm d} t}&=
		-\frac{1} {2\dx\dy}  \int_{I_j}\left[ f\left(u(x_{i-\frac12},y,t)\right)+ f\left(u(x_{i+\frac12},y,t)\right)\right]\mathrm{d}y
		+\frac{1}{(\dx)^2 \dy}\int_{I_{i,j}}f(u)\mathrm{d}x\mathrm{d}y
		\\&
		-\frac{1}{\dx\dy}\int_{I_i}\left[g\left(u(x,y_{j+\frac12},t)\right)-g\left(u(x,y_{j-\frac12},t)\right)\right]\frac{(x-x_i)}{\dx}\mathrm{d}x,
		\\
		\frac{{\rm d} \bar w_{i,j}(t)}{{\rm d} t}&=
		-\frac{1}{\dx\dy}  \int_{I_j}\left[ f\left(u(x_{i+\frac12},y,t)\right)- f\left(u(x_{i-\frac12},y,t)\right)\right]\frac{(y-y_j)}{\dy}\mathrm{d}y
		\\&
		\quad-\frac{1} {2\dx\dy}  \int_{I_i}\left[ g\left(u(x,y_{j-\frac12},t)\right)+ g\left(u(x,y_{j+\frac12},t)\right)\right]dx+\frac 1{\dx (\dy)^2}\int_{I_{i,j}}g(u)\mathrm{d}x\mathrm{d}y,
	\end{aligned}
	\right.
\end{equation}
where $\bar{u}_{i,j}(t)=\frac{1}{\dx\dy}\int_{I_{i,j}} u(x,y,t)\mathrm{d}x\mathrm{d}y$, $\bar{v}_{i,j}(t)=\frac{1}{\dx\dy}\int_{I_{i,j}}u(x,y,t)\frac{x-x_i}{\dx} \mathrm{d}x\mathrm{d}y$ and $\bar {w}_{i,j}(t)=$ $\frac{1}{\dx\dy}\int_{I_{i,j}}u(x,y,t)\frac{y-y_j}{\dy} \mathrm{d}x\mathrm{d}y$ are the zeroth-order moment (cell-average), the first-order moment in the $x$-direction and the first-order moment in the $y$-direction, respectively.

Let $\{\hat{x}^G_i\}^3_{G=1}$ and $\{\hat{y}^G_j\}^3_{G=1}$ denote three Gauss points in the intervals $I_i$ and $I_j$, respectively, and $\{\hat{\omega}_G\}^3_{G=1}$ are the weights of Gauss quadrature formula on a interval $[-\frac12,\frac12]$, i.e. 
\[
\hat{x}^{1}_i=x_{i-\frac{\sqrt{15}}{10}},~\hat{x}^{2}_i=x_{i},~\hat{x}^{3}_i=x_{i+\frac{\sqrt{15}}{10}},~
\hat{y}^{1}_j=y_{j-\frac{\sqrt{15}}{10}},~\hat{y}^{2}_j=y_{j},~\hat{y}^{3}_j=y_{j+\frac{\sqrt{15}}{10}},
\]
with the normalized weights $\hat{\omega}_{1,3}=\frac{5}{18}$ and $\hat{\omega}_2=\frac{4}{9}$.
We use the Gauss quadrature formula to approximate the integral terms over $I_i$, $I_j$ and $I_{i,j}$, and apply high order Lax-Friedrichs numerical fluxes to reconstruct the values of flux functions $f(u(x_{i+\frac12},y,t))$ and $g(u(x,y_{j+\frac12},t))$ at specified points, then a conservative semi-discrete scheme is defined as
\begin{equation}\label{sec2:2dSemi_HWENO}
	\left\{
	\begin{aligned}
		\frac{{\rm d} \bar{{u}}_{i,j}(t)}{{\rm d} t}=& -\frac{1}{\dx}\sum\limits_{G=1}^3\hat{\omega}_G(\hat{f}_{i+\frac12,G}-\hat{f}_{i-\frac12,G})
		-\frac{1}{\dy}\sum\limits_{G=1}^3\hat{\omega}_G(\hat{g}_{G,j+\frac12}-\hat{g}_{G,j-\frac12})
		\triangleq \mathcal{F}^1_{i,j}( \bar{{u}}, \bar{{v}},\bar{{w}} ), \\
		\frac{{\rm d} \bar{{v}}_{i,j}(t)}{{\rm d} t}=& -\frac{1}{2\dx}\sum\limits_{G=1}^3\hat{\omega}_G(\hat{f}_{i-\frac12,G}+\hat{f}_{i+\frac12,G})
		+\frac{1}{\dx}\sum\limits_{G=1}^3\sum\limits_{H=1}^3\hat{\omega}_G\hat{\omega}_H f(u(\hat{x}_i^{G},\hat{y}_j^H))
		\\&-\frac{1}{\dy}\sum\limits_{G=1}^3\hat{\omega}_G\frac{\hat{x}_i^G-x_i}{\dx}(\hat{g}_{G,j+\frac12}-\hat{g}_{G,j-\frac12})
		\triangleq \mathcal{F}^2_{i,j}( \bar{{u}}, \bar{{v}},\bar{{w}} ), \\
		\frac{{\rm d} \bar{{w}}_{i,j}(t)}{{\rm d}t}=&
		-\frac{1}{\dx}\sum\limits_{G=1}^3\hat{\omega}_G\frac{\hat{y}_j^G-y_j}{\dy}(\hat{f}_{i+\frac12,G}-\hat{f}_{i-\frac12,G})
		-\frac{1}{2\dy}\sum\limits_{G=1}^3\hat{\omega}_G(\hat{g}_{G,j-\frac12}+\hat{g}_{G,j+\frac12})
		\\&+\frac{1}{\dy}\sum\limits_{G=1}^3\sum\limits_{H=1}^3\hat{\omega}_G\hat{\omega}_H g(u(\hat{x}_i^{G},\hat{y}_j^H))\triangleq \mathcal{F}^3_{i,j}( \bar{{u}}, \bar{{v}},\bar{{w}} ),
	\end{aligned}
	\right.
\end{equation}
where $\mathcal{F}^1_{i,j}( \bar{{u}}, \bar{{v}},\bar{{w}} )$, $\mathcal{F}^2_{i,j}( \bar{{u}}, \bar{{v}},\bar{{w}} )$ and $\mathcal{F}^3_{i,j}( \bar{{u}}, \bar{{v}},\bar{{w}} )$ are the right-hand terms. The numerical fluxes $\hat{f}_{i+\frac12,G}$ and $\hat{g}_{G,j+\frac12}$ are used to approximate the values of $f(u(x_{i+\frac12},y,t))$ and $g(u(x,y_{j+\frac12},t))$ at the points $\{\hat{x}^G_i\}^3_{G=1}$ and $\{\hat{y}^G_j\}^3_{G=1}$, respectively, defined as
\begin{equation*}\label{2dflx}
	\hat{f}_{i+\frac12,G}=\frac12 \left[f(u^-_{i+\frac12,G})+f(u^+_{i+\frac12,G}) - \alpha_1 (u^+_{i+\frac12,G}-u^-_{i+\frac12,G}) \right],
\end{equation*}
\begin{equation*}\label{2dfly}
	\hat{g}_{G,j+\frac12}=\frac12\left[g(u^-_{G,j+\frac12})+g(u^+_{G,j+\frac12}) - \alpha_2 (u^+_{G,j+\frac12}-u^-_{G,j+\frac12} \right],
\end{equation*}
with $\alpha_1= \max\limits_u|f'(u)|$ and $\alpha_2=\max\limits_u|g'(u)|$. $\{u^\pm_{i+\frac12,G}\}^3_{G=1}$ and $\{u^\pm_{G,j+\frac12}\}^3_{G=1}$ are the values of $u(x,y,t)$ at the points $\{(x_{i+\frac12},\hat{y}^G_j)\}$ and $\{(\hat{x}^G_i,y_{j+\frac12})\}^3_{G=1}$, respectively.
The superscribes $``-"$ and $``+"$ of $\{u^\pm_{i+\frac12,G}\}^3_{G=1}$ represent the limits from the left and right sides at the interface $x_{i+\frac12}$, respectively. Similarly, the superscribes $``-"$ and $``+"$ of $\{u^\pm_{G,j+\frac12}\}^3_{G=1}$ indicate the limits from the bottom and top sides at the interface $y_{j+\frac12}$, respectively.

{The equations \eqref{sec2:2dSemi_HWENO} also can be expressed as the semi-discrete form of the $\mathbbm{P}^1(I_{i,j})$ DG  finite element method}. However, for moment-based  HWENO schemes, a Hermite reconstruction is used to approximate $\{u^\pm_{i\mp\frac12,G}\}^3_{G=1}$, $\{u^\pm_{G,j\mp\frac12}\}^3_{G=1}$ and $\{u(\hat{x}^G_i,\hat{y}^H_j)\}^3_{G,H=1}$ based on the zeroth- and first-order moments in the cells $\{I_{i-1,j-1},$ $I_{i,j-1},$ $I_{i+1,j-1},$ $I_{i-1,j},$ $I_{i,j},$ $I_{i+1,j},$ $I_{i-1,j+1},$ $I_{i,j+1},$ $I_{i+1,j+1}\}$.
To simplify the representation, we rebel the cell $I_{i,j}$ and its adjacent cells as $I_1,...,I_9$, e.g., $I_{i,j}\triangleq I_5$. Let $\{\bar{{u}}_k,\bar{{v}}_k,\bar{{w}}_k\}$ denote the zeroth- and first-order moments of the cell $I_k$, e.g., $\{\bar{{u}}_{i,j}\triangleq \bar{{u}}_5, \bar{{v}}_{i,j}\triangleq \bar{{v}}_5,  \bar{{w}}_{i,j}\triangleq \bar{{w}}_5 \}$. 	

Similar to the  one-dimensional case for utilizing unified stencils to construct a stable scheme, the first-order moments $\bar{{v}}_5$ and $\bar{{w}}_5$ in the central mesh are no longer used in spatial reconstructions. Instead, their high order modified terms $\hat{{v}}_5$ and $\hat{{w}}_5$ are also obtained by using the same information from the spatial reconstructions and incorporated only into time discretizations.
Next, we will provide the detailed reconstructed procedures for $u(x,y)$ at specific points, and the modified terms $\hat{{v}}_5$ and $\hat{{w}}_5$, based on the values $\{\bar{{u}}_1,\ldots,\bar{{u}}_9, \bar{{v}}_k,\bar{{w}}_k\}_{k=2,4,6,8}$.

\textbf{Step 1.} Reconstruct a quartic polynomial $p_0(x,y)$ and four linear polynomials $\{p_m(x,y)\}_{m=1}^{4}$.

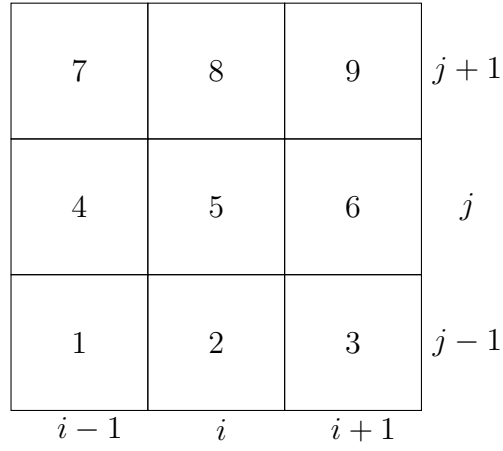
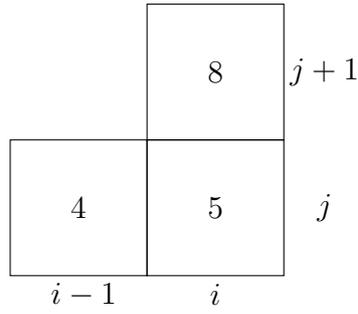
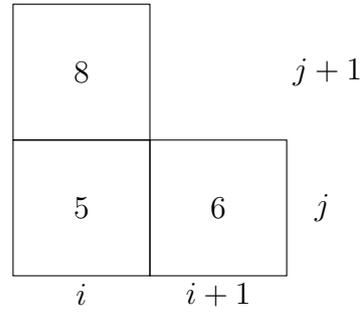
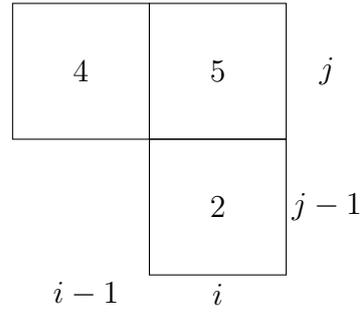
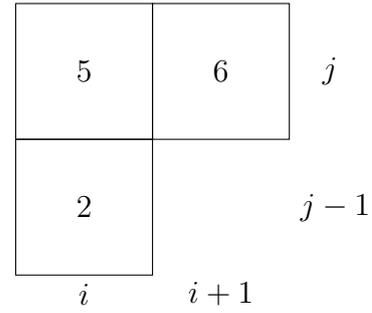
\begin{figure}
	\centering
	\subfigure[~$S_0$]{
		\tikzset{global scale/.style={scale=#1,every node/.append style={scale=#1}}}
		\centering
		\begin{tikzpicture}[global scale = 1]
			\draw(0,0)rectangle+(1.8,1.8);\draw(1.8,0)rectangle+(1.8,1.8);\draw(3.6,0)rectangle+(1.8,1.8);
			\draw(0,1.8)rectangle+(1.8,1.8);\draw(1.8,1.8)rectangle+(1.8,1.8);\draw(3.6,1.8)rectangle+(1.8,1.8);
			\draw(0,3.6)rectangle+(1.8,1.8);\draw(1.8,3.6)rectangle+(1.8,1.8);\draw(3.6,3.6)rectangle+(1.8,1.8);
			\draw(0.9,0.9)node{1};
			\draw(2.7,0.9)node{2};
			\draw(4.5,0.9)node{3};
			\draw(0.9,2.7)node{4};
			\draw(2.7,2.7)node{5};
			\draw(4.5,2.7)node{6};
			\draw(0.9,4.5)node{7};
			\draw(2.7,4.5)node{8};
			\draw(4.5,4.5)node{9};
			\draw(0.9,-0.25)node{$~~ i-1$};
			\draw(2.7,-0.25)node{$~i$};
			\draw(4.5,-0.25)node{$~~i+1$};
			\draw(5.8,0.9)node{$~~~j-1$};
			\draw(5.8,2.7)node{$~~~j$};
			\draw(5.8,4.5)node{$~~~j+1$};
		\end{tikzpicture}
	}
	
	\subfigure[~$S_3$]{
		\tikzset{global scale/.style={scale=#1,every node/.append style={scale=#1}}}
		\centering
		\begin{tikzpicture}[global scale = 1]
			\draw(0+0,0+14.4)rectangle+(1.8,1.8);  \draw(0.9+0,0.9+14.4)node{4};
			\draw(1.8+0,0+14.4)rectangle+(1.8,1.8);\draw(2.7+0,0.9+14.4)node{5};
			\draw(1.8+0,1.8+14.4)rectangle+(1.8,1.8);\draw(2.7+0,2.7+14.4)node{8}; 			
			
			\draw(0.9+0,-0.25+14.4)node{$~i-1$};\draw(2.7+0,-0.25+14.4)node{$i$};
			\draw(4+0,0.9+14.4)node{$~~j$};\draw(4+0,2.7+14.4)node{$~~j+1$};
		\end{tikzpicture}
	}
	\subfigure[~$S_4$]{
		\tikzset{global scale/.style={scale=#1,every node/.append style={scale=#1}}}
		\centering
		\begin{tikzpicture}[global scale = 1]
			\draw(0+5.0,0+14.4)rectangle+(1.8,1.8);  \draw(0.9+5.0,0.9+14.4)node{5};
			\draw(1.8+5.0,0+14.4)rectangle+(1.8,1.8);\draw(2.7+5.0,0.9+14.4)node{6};
			\draw(0+5.0,1.8+14.4)rectangle+(1.8,1.8);\draw(0.9+5.0,2.7+14.4)node{8};
			\draw(0.9+5.0,-0.25+14.4)node{$i$};\draw(2.7+5.0,-0.25+14.4)node{$i+1$};
			\draw(4+5.0,0.9+14.4)node{$~j$};\draw(4+5.0,2.7+14.4)node{$~~j+1$};
		\end{tikzpicture}
	}
	\\
	~\subfigure[~$S_1$]{
		\tikzset{global scale/.style={scale=#1,every node/.append style={scale=#1}}}
		\centering
		\begin{tikzpicture}[global scale = 1]
			\draw(1.8+0,0+9.6)rectangle+(1.8,1.8);\draw(2.7+0,0.9+9.6)node{2};
			\draw(0+0,1.8+9.6)rectangle+(1.8,1.8);\draw(0.9+0,2.7+9.6)node{4};
			\draw(1.8+0,1.8+9.6)rectangle+(1.8,1.8);\draw(2.7+0,2.7+9.6)node{5};			
			\draw(0.9+0,-0.25+9.6)node{$~i-1$};\draw(2.7+0,-0.25+9.6)node{$i$};
			\draw(4.0+0,0.9+9.6)node{$~~j-1$};\draw(4.0+0,2.7+9.6)node{$~~j$};
		\end{tikzpicture}
	}		
	\subfigure[~$S_2$]{
		\tikzset{global scale/.style={scale=#1,every node/.append style={scale=#1}}}
		\centering
		\begin{tikzpicture}[global scale = 1]
			\draw(0+5.0,0+9.6)rectangle+(1.8,1.8);  \draw(0.9+5.0,0.9+9.6)node{2};
			\draw(0+5.0,1.8+0+9.6)rectangle+(1.8,1.8);\draw(0.9+5.0,2.7+9.6)node{5};
			\draw(1.8+5.0,1.8+9.6)rectangle+(1.8,1.8);\draw(2.7+5.0,2.7+9.6)node{6};
			\draw(0.9+5.0,-0.25+9.6)node{$i$};\draw(2.7+5.0,-0.25+9.6)node{$i+1$};
			\draw(4.0+5.0,0.9+9.6)node{$~~~j-1$};\draw(4.0+5.0,2.7+9.6)node{$~~j$};
		\end{tikzpicture}
	}	
	\caption{The big stencil $S_0$, small stencils $\{S_m\}_{m=1}^4$ and their respective labels.}
	\label{sec2:stenciL2d}
\end{figure}

Firstly, we consider a big stencil $S_0$ and four small stencils $\{S_m\}^4_{m=1}$ shown in Fig. \ref{sec2:stenciL2d}. Here, we use the values $\{\bar{{u}}_1,\ldots,\bar{{u}}_9, \bar{{v}}_k,\bar{{w}}_k\}_{k=2,4,6,8}$, $\{\bar{{u}}_k\}_{k=2,4,5}$, $\{\bar{{u}}_k\}_{k=2,5,6}$, $\{\bar{{u}}_k\}_{k=4,5,8}$, and $\{\bar{{u}}_k\}_{k=5,6,8}$ in the stencils $\{S_m\}^4_{m=0}$, respectively, then,  a quartic polynomial $p_0(x,y)\in \mathbbm{P}^4(I_{i,j})$ is obtained by a Hermite reconstruction based on $S_0$, satisfying
\begin{equation}\label{sec:p0(x,y)}
		\begin{aligned}
			&\frac{1}{\dx\dy}\int_{I_k}p_0(x,y) \mathrm{d}x\mathrm{d}y=\bar{{u}}_k,~k=1,...,9,\\
			&\frac{1}{\dx\dy}\int_{I_{k}}p_0(x,y)\frac{x-x_{k}}{\dx}\mathrm{d}x\mathrm{d}y=\bar{{v}}_{k},~k={2,4,6,8},\\
			&\frac{1}{\dx\dy}\int_{I_{k}}p_0(x,y)\frac{y-y_{k}}{\dy}\mathrm{d}x\mathrm{d}y=\bar{{w}}_{k}, ~k={2,4,6,8},
		\end{aligned}
\end{equation}
and four linear polynomials $\{p_m(x,y)\}^4_{m=1} \in \mathbbm{P}^1(I_{i,j})$ are obtained by a linear reconstruction, satisfying
\begin{equation*}
		\frac{1}{\dx\dy}\int_{I_k}p_m(x,y) \mathrm{d}x\mathrm{d}y=\bar{{u}}_k,
\end{equation*}
for
\[
\begin{aligned}
	~m=1,k=2,4,5;~m=2,k=2,5,6;\\~m=3,k=4,5,8;~m=4,k=5,6,8.
\end{aligned}
\]
The quartic polynomial $p_0(x,y)$ can be uniquely determined  by requiring it to exactly match $\{\bar{{u}}_1,\ldots,\bar{{u}}_9, \bar{{v}}_4,\bar{{v}}_6,\bar{{w}}_2, \bar{{w}}_8\}$ with the least square methodology in \cite{hs,ZhaoZQiuJX_ArtificalLinearWeight_HWENO2020}, while the four polynomials $\{p_m(x,y)\}^4_{m=1}$ can be directly obtained by solving $3\times3$ linear systems.

{Still, $p_0(x,y)$  can be written as
\begin{equation}\label{sec22:linearweights}
	{p_0(x,y) = \gamma_0 \left( \frac 1 {\gamma_0}p_0(x,y) - \sum\limits_{m=1}^{4}\frac {\gamma_m} {\gamma_0} p_m(x,y) \right) + \sum\limits_{m=1}^{4}\gamma_n p_m(x,y), \quad \gamma_0\neq0,}
\end{equation}
where $\{\gamma_m\}^4_{m=0}$ also are arbitrary positive linear weights with $\sum^4_{m=0}\gamma_m$=1.}

\textbf{Step 2.} Compute the smoothness indicators of $\{p_m(x,y)\}_{m=0}^4$ by the definition as in \cite{hs,ZhaoZQiuJX_ArtificalLinearWeight_HWENO2020}, given by
\begin{equation}\label{sec2:def_beta2}
	\beta_m= \sum_{|l|=1}^r|I_{i,j}|^{|l|-1} \int_{I_{i,j}}\left( \frac {\partial^{|l|}}{\partial x^{l_1}\partial y^{l_2}}p_m(x,y)\right)^2 \mathrm{d}x\mathrm{d}y, \quad m=0,...,4,
\end{equation}
where $l=(l_1,l_2)$, $|l|=l_1+l_2$ and $r$ is the degree of $p_m(x, y)$. Similar to the one-dimensional case, let $p_m(x,y)=\sum\limits_{n=0}^r c_{m,n} \phi_n(x,y)$, where the basis functions $\phi_n(x,y)$ are defined as
\begin{equation*}
	\begin{aligned}
		&\phi_0=1,~\phi_1=\xi_i,~\phi_2= \eta_j,~
		\phi_3=\xi_i^2,~
		\phi_4=\xi_i\eta_i,
		\phi_5=\eta_j^2,~
		\phi_6=\xi_i^3,~
		\phi_7=\xi_i^2\eta_j,~~
		\\&
		\phi_8=\xi_i\eta_j^2,~
		\phi_9=\eta_j^3,~
		\phi_{10}=\xi_i^4,
		\phi_{11}=\xi_i^3\eta_j,~
		\phi_{12}= \xi_i^2\eta_j^2,~
		\phi_{13}=\xi_i\eta_j^3,~
		\phi_{14}=\eta_j^4,~
		\phi_{15}=\xi_i^5,\cdots,
	\end{aligned}
\end{equation*}	
with $\xi_i=\frac{x-x_i}{\dx}$ and $\eta_j=\frac{y-y_j}{\dy}$. Then the explicit expression of the smoothness indicators are
\begin{equation}\label{sec2:beta2d}
\left\{
\begin{aligned}
	\beta_0&=
	\frac{1}{2}(c_{0,1} + \frac{1}{2}c_{0,6})^2 + \frac{1}{2}(c_{0,1} + \frac{1}{6}c_{0,8})^2  + \frac{1}{2}(c_{0,2} + \frac{1}{6}c_{0,7})^2 + \frac{1}{2}(c_{0,2} + \frac{1}{2}c_{0,9})^2 + \frac{13}{6}(c_{0,3} + \frac{63}{65}c_{0,10})^2 \\&
	 + \frac{13}{6}(c_{0,3} + \frac{1}{6}c_{0,12})^2 + \frac{7}{12}(c_{0,4} + \frac{17}{35}c_{0,11})^2 + \frac{7}{12}(c_{0,4} + \frac{17}{35}c_{0,13})^2 + \frac{13}{6}(c_{0,5} + \frac{1}{6}c_{0,12})^2 \\& +\frac{13}{6}(c_{0,5} + \frac{63}{65}c_{0,14})^2 + \frac{3119}{80}(c_{0,6} + \frac{5}{9357}c_{0,8})^2 +\frac{ 3379}{720}(c_{0,7} + \frac{15}{3379}c_{0,9})^2 + \frac{2634769}{561420} c_{0,8}^2 \\&   + \frac{2634769}{67580}c_{0,9}^2   +\frac{5676583}{9100}(c_{0,10} + \frac{3185}{11353166}c_{0,12})^2 + \frac{709573}{16800}(c_{0,11} + \frac{1155}{709573}c_{0,13})^2 \\& + \frac{230094013357}{12261419280}(c_{0,12} + \frac{2145748374}{230094013357}c_{0,14})^2 + \frac{31468281769}{745051650}c_{0,13}^2 + \frac{25118160529227568}{40266452337475}c_{0,14}^2,\\
	\beta_m&= c^2_{m,1} + c^2_{m,2},m=1,2,3,4,
\end{aligned}
\right.
\end{equation}
where the coefficients of the polynomials $\{p_m(x)\}^4_{m=0}$ are listed in Appendix \ref{sec:Appendix1}.

\textbf{Step 3.} Compute nonlinear weights based on linear weights and smoothness indicators.
Similar to the one-dimensional case, we also use a parameter $\tau$ to measure the overall difference between $\{\beta_{m}\}_{m=0}^4$,
\begin{equation}
	\label{tao4}
	\tau=\left(\frac{|\beta_{0}-\beta_{1}|+|\beta_{0}-\beta_{2}|+|\beta_{0}-\beta_{3}|+|\beta_{0}-\beta_{4}|}{4}\right)^2,
\end{equation}
then we compute the nonlinear weights by
\begin{equation}\label{sec2:non_weight_2d}
	\omega_m=\frac{\widetilde\omega_m}{
		\widetilde\omega_0+\ldots+\widetilde\omega_4},
\ \mbox{with} \ \widetilde\omega_{m}=\gamma_m(1+\frac{\tau}{\beta_{m}u_{ave}+\varepsilon}),~ m=0,\ldots,4,
\end{equation}		
where $u_{ave} =(\frac{|\overline u_2|+|\overline u_4|+|\overline u_5|+|\overline u_6|+|\overline u_8|}{5})^2$, and $\varepsilon=10^{-10}$ is to avoid zero denominator.
 The nonlinear weights also preserve the scale-invariant property for adding $u_{ave}$, while the original formulation employed in the WENO and HWENO schemes with artificial linear weights \cite{ZhuJunQiuJX2017fvWENOZQ, ZhaoZQiuJX_ArtificalLinearWeight_HWENO2020} (the special case of $u_{{ave}}=1$ in Eq. \eqref{sec2:non_weight_2d}) lacks this fundamental property. Consequently, the original definition for the  nonlinear weights in \cite{ZhuJunQiuJX2017fvWENOZQ, ZhaoZQiuJX_ArtificalLinearWeight_HWENO2020} leads to noticeable oscillations when simulating the Sedov blast wave problem, as exhibited in Example \ref{Example:Sedov2D}.

 Finally,  by replacing a part of the linear weights in \eqref{sec22:linearweights} by the nonlinear weighs \eqref{sec2:non_weight_2d}, a nonlinear HWENO reconstructed polynomial $u_{i,j}(x,y)$ is obtained for $u(x,y)$. Also, using the same polynomials, linear and nonlinear weights, the high order modification  $\hat{{v}}_{i,j}$ for the first-order moment $\bar{{v}}_{i,j}$ and the high order modification $\hat{{w}}_{i,j}$ for the first-order moment $\bar{{w}}_{i,j}$ are obtained simultaneously, but the modified values are only used in the following time discretizations, having
\begin{equation}\label{sec2:2dHWENORecon}
	\begin{cases}
		u_{i,j}(x,y) = \omega_0 \left( \frac 1 {\gamma_0}p_0(x,y) - \sum\limits_{m=1}^{4}\frac {\gamma_m} {\gamma_0} p_m(x,y) \right) + \sum\limits_{m=1}^{4}\omega_n p_m(x,y),
		\\
		\hat{{v}}_{i,j}=\frac{1}{\dx\dy}\int_{I_{i,j}} u_{i,j}(x,y)\frac{x-x_{i}}{\dx}\mathrm{d}x\mathrm{d}y=\omega_0\left(\frac{1}{\gamma_0}q^v_0 - \sum\limits_{m=1}^4 \frac{\gamma_m}{\gamma_0}q^v_m\right) +\sum\limits_{m=1}^4\omega_{m}q^v_m,~
		\\
		\hat{{w}}_{i,j}=\frac{1}{\dx\dy}\int_{I_{i,j}} u_{i,j}(x,y)\frac{y-y_{j}}{\dy}\mathrm{d}x\mathrm{d}y\omega_0\left(\frac{1}{\gamma_0}q^w_0 - \sum\limits_{m=1}^4 \frac{\gamma_m}{\gamma_0}q^w_m\right) +\sum\limits_{m=1}^4\omega_{m}q^w_m,
	\end{cases}
\end{equation}
where $q^v_m=\frac{1}{\dx\dy}\int_{I_{i,j}} p_m(x,y)\frac{x-x_{i}}{\dx}\mathrm{d}x\mathrm{d}y$ and  $q^w_m=\frac{1}{\dx\dy}\int_{I_{i,j}} p_m(x,y)\frac{y-y_{j}}{\dy}\mathrm{d}x\mathrm{d}y$. Then, the values at specific points  that we need are computed as below:
\begin{equation*}
	u^\mp_{i\pm\frac12,G} = u_{i,j}(x_{i\pm\frac12},\hat{y}_j^G),~
	u^\mp_{G,j\pm\frac12} = u_{i,j}(\hat{x}_i^G,y_{i\pm\frac12}),~
	u(\hat{x}^G_i,\hat{y}^{H}_j) = u_{i,j}(\hat{x}_i^G,\hat{y}_j^{H}),~
	G,H={1,2,3}.
\end{equation*}

\textbf{Step 4.} Time discretizations for the semi-discrete scheme \eqref{sec2:2dSemi_HWENO}.

As in the one-dimensional case, the modified third-order SSP RK method \eqref{sec2:1dRK3_HWENO} is also used to solve the two-dimensional semi-discrete scheme \eqref{sec2:2dSemi_HWENO}. Differently, the involved variables are $(\overline{u}_{i,j},\overline{v}_{i,j},\overline{w}_{i,j})$. Also,
the modified terms $\hat{{v}}_{i,j}$ and $\hat{{w}}_{i,j}$ are treated solely as time stage values, and they are obtained by the formula \eqref{sec2:2dHWENORecon}.


\begin{rem}
	For one- and two-dimensional compressible Euler equations, the HWENO procedures are used in cooperation with the local characteristic decomposition to avoid spurious oscillations, which is similar to the classical WENO scheme \cite{js}. Besides, the computation of $u_{ave}$ in \eqref{sec2:non_weight_1d} and \eqref{sec2:non_weight_2d} is also implemented in the local characteristic direction for $u_{ave}$ relies on the reconstructed variable.
\end{rem}

\subsection{Stability analysis}\label{sec2:stabilityHWENO}
In this subsection, we present the stability analysis for the proposed HWENO-U scheme by the Fourier analysis method. This potent technique for stability analysis depends heavily on the assumption of uniform meshes and periodic boundary conditions. Additionally, it is only effective for the linear scheme used to solve a scalar linear equation.

For simplicity of analysis, we consider the one-dimensional  linear equation
\begin{equation}\label{sec3:linearEq}
	u_t+au_x=0, x\in[0,2\pi], t>0,
\end{equation}
with constant coefficient $a$. Assume $a=1$, then the semi-discrete finite volume HWENO scheme \eqref{sec2:1dSemi_HWENO} reads
\begin{equation}\label{sec3:1dlinear_HWENO}
	\left\{
	\begin{aligned}
		&\frac{{\rm d} \bar{u}_i(t)}{{\rm d} t}=-\frac{1}{\dx}(u^-_{i+\frac12}-u^-_{i-\frac12}),
		\\
		&\frac{{\rm d} \bar{v}_i(t)}{{\rm d} t}=-\frac{1}{2\dx}(u^-_{i-\frac12}+u^-_{i+\frac12})+\frac{1}{\dx} \bar{u}_i.
	\end{aligned}
	\right.
\end{equation}
Here, we first use the moments $\{\bar{u}_{i-1}, \bar{u}_{i}, \bar{u}_{i+1}, \bar{v}_{i-1}, \bar{v}_{i+1}\}$ to reconstruct $u^-_{i+\frac12}$ linearly, e.g., $u^-_{i+\frac12}=\frac{269}{456}\bar{u}_{i-1} +\frac{7}{12}\bar{u}_{i} -\frac{79}{456}\bar{u}_{i+1}+\frac{177}{76}\bar{v}_{i-1}+\frac{63}{76}\bar{v}_{i+1}$. Substituting $u^-_{i+\frac12}$ into equations \eqref{sec3:1dlinear_HWENO} gives
\begin{equation}\label{sec3:1dlinearHWENO_vector}
	\frac{{\rm d}\emph{\textbf{u}}_i(t)}{{\rm d}t}=\frac{1}{\dx}(\mathbf{A}\emph{\textbf{u}}_{i-2} +\mathbf{B}\emph{\textbf{u}}_{i-1} +\mathbf{C}\emph{\textbf{u}}_{i} + \mathbf{D}\emph{\textbf{u}}_{i+1}),
\end{equation}
where $\emph{\textbf{u}}_i(t)=(\bar{u}_i(t),\bar{v}_i(t))^\mathsf{T}$, $\mathbf{A}$, $\mathbf{B}$, $\mathbf{C}$ and $\mathbf{D}$ are $2\times2$ constant matrices given by
\[
\mathbf{A}=\left[ \begin {array}{cc} -{\frac{79}{456}}&-{\frac{63}{76}}
\\ \noalign{\medskip}{\frac{79}{912}}&{\frac{63}{152}}\end {array}
\right],
\mathbf{B}=\left[ \begin {array}{cc} {\frac{115}{152}}&{\frac{63}{76}}
\\ \noalign{\medskip}-{\frac{187}{912}}&{\frac{63}{152}}\end {array}
\right],
\mathbf{C}=\left[ \begin {array}{cc} {\frac{1}{152}}&-{\frac{177}{76}}
\\ \noalign{\medskip}{\frac{377}{912}}&{\frac{177}{152}}\end {array}
\right],
\mathbf{D}=\left[ \begin {array}{cc} -{\frac{269}{456}}&{\frac{177}{76}}
\\ \noalign{\medskip}-{\frac{269}{912}}&{\frac{177}{152}}\end {array}
\right].
\]
For the stability of scheme \eqref{sec3:1dlinear_HWENO}, we have the following conclusions.

\begin{figure}[!ht]
	\centering
	\subfigure[{The polar plot of $\lambda_{1,2}(\widetilde{\mathbf{G}})$ in \eqref{sec3:lambda12}, $K\in [0,2\pi]$}]
	{\includegraphics[width=8cm,angle=0]{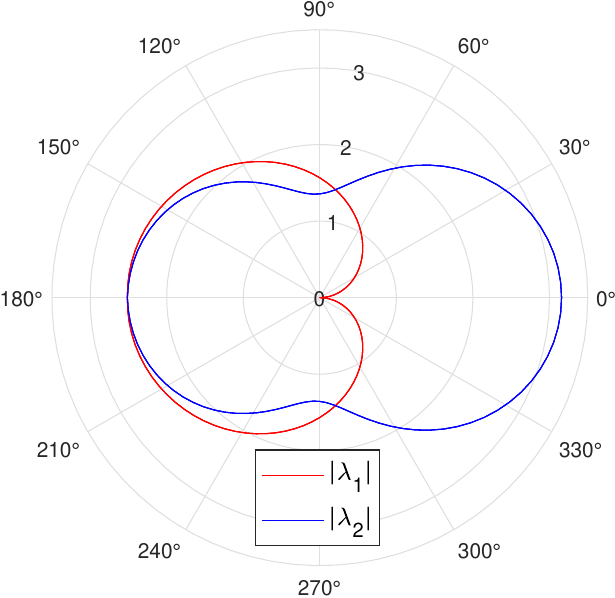}}
	\subfigure[The polar plot of {spectral radius $\rho(\hat{\bfG}_3)$} in \eqref{RK3_G_HWENO_U} with different  {Courant-Friedrichs-Lewy {(CFL)}} number for the HWENO-U scheme.]
	{\includegraphics[width=8cm,angle=0]{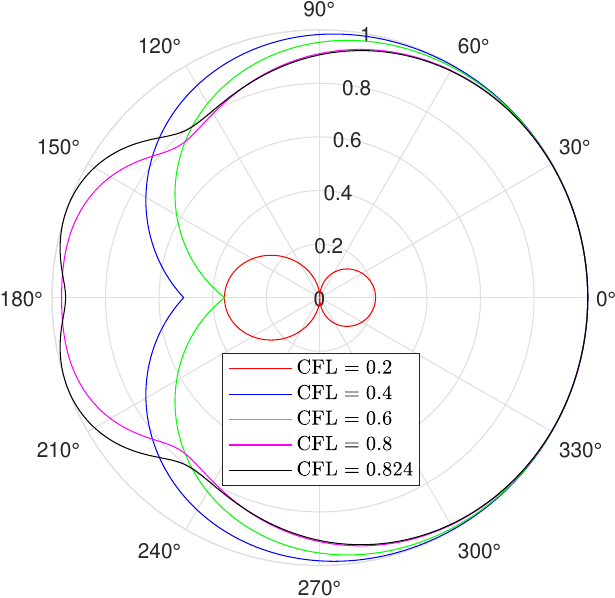}}
	\caption{ The polar plot of the eigenvalues and amplification matrix. }
	\label{CFL_HWENO_U}
\end{figure}

\begin{theorem}\label{sec2:thm1}
	\emph{
		Combining the semi-discrete HWENO scheme \eqref{sec3:1dlinear_HWENO} with either the stand forward  Euler or   third-order SSP RK time discretization \cite{s2}, the resulting schemes are both unstable.
	}
	
	{Proof.}
	The semi-discrete HWENO scheme \eqref{sec3:1dlinear_HWENO} with the forward Euler time method is
	\begin{equation}\label{sec3:1dlinearHWENO_Euler}
		\emph{\textbf{u}}_i^{n+1} =\emph{\textbf{u}}_i^{n}+ \frac{\dt}{\dx}(\mathbf{A}\emph{\textbf{u}}_{i-2}^{n}+\mathbf{B}\emph{\textbf{u}}_{i-1}^{n}+\mathbf{C}\emph{\textbf{u}}_{i}^{n}+\mathbf{D}\emph{\textbf{u}}_{i+1}^{n}),
	\end{equation}	
	To apply the {von-Neumann analysis}, we have an assumption on the solution
	\begin{equation}\label{sec3:1dFourier1}
		\emph{\textbf{u}}_{i}^{n} = \hat{\emph{\textbf{u}}}^{n}e^{\sigma ik\dx},
	\end{equation}
	where $\sigma$ is the imaginary unit satisfying $\sigma^2=-1$, and  $k$ is the wave number.
	We expect that
	\begin{equation}\label{sec3:1dFourier2}
		\emph{\textbf{u}}_{i}^{n+1} = \hat{\emph{\textbf{u}}}^{n+1}e^{\sigma ik\dx},
	\end{equation}
	where $\hat{\emph{\textbf{u}}}^{n+1}={\mathbf{G}_1}\hat{\emph{\textbf{u}}}^{n}$ and ${\mathbf{G}_1}$ is the amplification matrix. Substituting \eqref{sec3:1dFourier1} and \eqref{sec3:1dFourier2} into \eqref{sec3:1dlinearHWENO_vector} gives
	\begin{equation}\label{RK1_G}
		{ {\mathbf{G}_1}}= {\mathbf{I}}+\frac{\dt}{\dx}\widetilde{\mathbf{G}},
	\end{equation}
	where $\widetilde{\mathbf{G}}= \mathbf{A} e^{-2\sigma K}+\mathbf{B}e^{-\sigma K}+\mathbf{C}+\mathbf{D}e^{\sigma K}, K= k\dx\in [0,2\pi]$ is a simplified wave number. The two eigenvalues of $\widetilde{\mathbf{G}}$ are 
	\begin{equation}\label{sec3:lambda12}
		\begin{aligned}
			\lambda_{1,2}&=\frac{1}{456}(267+131e^{\sigma iK}+267e^{-\sigma iK}+(55\pm\sqrt{\theta}) e^{-2\sigma iK}),
		\end{aligned}
	\end{equation}
	where $\theta =17161e^{6\sigma iK}+554226e^{5\sigma iK}-343029e^{4\sigma iK}+329356e^{3\sigma iK}-71709e^{2\sigma iK}+29370e^{\sigma iK}+3025 \neq 0,~ \forall K\in [0,2\pi]$. According to the von-Neumann stability analysis of Section 2.2 in \cite{Strikwerda}, {the necessary condition of stability for the scheme \eqref{sec3:1dlinearHWENO_Euler} is  the spectral radius $\rho(\bfG_1)\le 1$, which is equivalent to $\max\limits_{i=1,2}|\lambda_{i}(\bfG_1)|\le 1$}, $ \forall K\in [0,2\pi]$. However, we can see that $\lambda_{1,2}$ in \eqref{sec3:lambda12} violate
this condition  from the left of Fig. \ref{CFL_HWENO_U}. Therefore, the scheme \eqref{sec3:1dlinearHWENO_Euler} is unstable. Similarly, when combining the third-order SSP RK method, the amplification matrix of the resulting HWENO scheme is
\begin{equation}\label{RK3_G}
	{\mathbf{G}_3}=\frac13 \mathbf{I}+\frac12(\mathbf{I} +\frac{\dt}{\dx}  \widetilde{\mathbf{G}})\mathbf{I} +\frac16(\mathbf{I} +\frac{\dt}{\dx}  \widetilde{\mathbf{G}})^3.
\end{equation}
Obviously, {the spectral radius $\rho(\bfG_3)> 1$} in this case for the third-order RK method that is a convex combination of forward  Euler methods,  indicating that the resulting scheme is unstable. $\square$ 
\end{theorem}

Notice that the first-order moment $\bar{v}_i$ is not used to approximate $u^-_{i+\frac12}$ in the spatial discretizations of the scheme \eqref{sec3:1dlinear_HWENO}, which makes the two fully-discrete schemes above unstable. However, we find that by combining the high order modification of the first-order moments in the time discretizations, the new fully-discrete HWENO scheme becomes stable even though the same saptial discretizations are employed. The provable process is presented below.

\begin{theorem}\label{sec2:thm2}
	\emph{When using the modified third-order SSP RK time time discretization \eqref{sec2:1dRK3_HWENO} to solve the scheme \eqref{sec3:1dlinear_HWENO}, {the necessary condition of stability for the resulting HWENO scheme is  {$0< \frac{\dt}{\dx}$} {$\lesssim0.824$}.}}
	
	{Proof.}
	Firstly, by employing the high order modification of the first-order moment, we utilize the modified forward Euler time-marching method to resolve the scheme \eqref{sec3:1dlinear_HWENO}, namely,
	\begin{equation}\label{sec3:HWENO1d_Euler}
		\emph{\textbf{u}}^{n+1}_i =
		\begin{bmatrix}
			\bar{u}^{n}_i \\ \hat{v}^{n}_i
		\end{bmatrix}
		+ \frac{\dt}{\dx}(\mathbf{A}\emph{\textbf{u}}_{i-2}^{n}+\mathbf{B}\emph{\textbf{u}}_{i-1}^{n}+\mathbf{C}\emph{\textbf{u}}_{i}^{n}+\mathbf{D}\emph{\textbf{u}}_{i+1}^{n}),
	\end{equation}
	where $\hat{v}^{n}_i=\frac{5}{76}(\bar{u}^{n}_{i+1}-\bar{u}^{n}_{i-1})-\frac{11}{38}(\bar{v}^{n}_{i+1}+\bar{v}^{n}_{i-1})$. Through applying the von-Neumann analysis, we obtain the amplification matrix $\hat{\mathbf{G}}_1 = \hat{\mathbf{A}} +\frac{\dt}{\dx}  \widetilde{\mathbf{{G}}}$ with
	\[
	\hat{\mathbf{A}}=\begin{bmatrix}
		1 & 0 \\ \frac{5}{76}(e^{\sigma iK}-e^{-\sigma iK}) & - \frac{11}{38}(e^{\sigma iK}+e^{-\sigma iK})
	\end{bmatrix}.
	\]
	Similarly, when using the third-order SSP RK time discretization \eqref{sec2:1dRK3_HWENO}, the amplification matrix becomes
	\begin{equation}\label{RK3_G_HWENO_U}
		\mathbf{\hat{G}}_3=\frac13\hat{\mathbf{A}}+\frac12(\hat{\mathbf{A}} +\frac{\dt}{\dx}  \widetilde{\mathbf{G}})\hat{\mathbf{A}}+\frac16(\hat{\mathbf{A}} +\frac{\dt}{\dx}  \widetilde{\mathbf{G}})^3,
	\end{equation}
    The polar plot of {spectral radius $\rho(\mathbf{\hat{G}}_3)$ }is presented in the right of Fig. \ref{CFL_HWENO_U} for the HWENO-U scheme with different CFL numbers. Therefore, {the necessary condition of stability for the resulting HWENO scheme is  $\rho(\mathbf{\hat{G}}_3)\le 1$, which is equivalent to $0<\frac{\dt}{\dx}$ $\lesssim 0.824$}. This value can be numerically determined  by sampling 10000 points for $K \in [0,2\pi]$.
	$\square$ 
\end{theorem}

\section{Numerical tests}

In this section, we present the numerical results of the benchmark and extreme examples to verify the fifth-order
accuracy, efficiency, high resolution, and robustness of the proposed HWENO-U scheme. For comparisons, we mainly consider the proposed HWENO-U, HWENO-M \cite{ZhaoZQiuJX_ArtificalLinearWeight_HWENO2020}, and WENO-ZQ \cite{ZhuJunQiuJX2017fvWENOZQ,trzq} schemes since the three schemes have the fifth-order accuracy and use arbitrary positive linear weights in spatial reconstructions. {Particularly, the results of the WENO-ZQ scheme  are computed by the methods of the structured finite volume version  \cite{ZhuJunQiuJX2017fvWENOZQ} and the unstructured finite volume version \cite{trzq} in the one- and two-dimensional cases, respectively.} For the HWENO-U scheme, the linear weights of the low-degree polynomials are set as $1/400$ both in one- and two-dimensional cases, and
the remaining linear weight is assigned to the high-degree polynomial, ensuring that their sum equals one.
For fair comparisons, the linear weights of the HWENO-M and WENO-ZQ schemes are chosen as they suggested from \cite{ZhuJunQiuJX2017fvWENOZQ,ZhaoZQiuJX_ArtificalLinearWeight_HWENO2020,trzq}.  Besides, a positivity-preserving (PP) limiter will be used to improve the robustness of the HWENO-U and HWENO-M schemes in some two-dimensional extreme problems. If not, the two schemes cannot work since  negative densities or pressures will arise, and we refer to \cite{FanZhangQiu2021PPhybridHWENO, CaiZhangQiu2016PPFVHWENOEuler} for the PP researches of finite volume HWENO schemes. The CFL number is set as 0.6. To compare the computational cost, we utilize the programming language Fortran 95 to execute our simulations on the environment of Inter(R) Xeon (R) Gold 6130 CPU @ 2.10 GHz.

\subsection{Accuracy tests}\label{sec3:AccuracyTests}
{In this subsection, we first verifies the fifth-order accuracy of the HWENO-U scheme.} Then, the comparisons of computational costs and errors for the HWENO-U, HWENO-M and WENO-ZQ schemes are presented to demonstrate that the HWENO-U scheme behaves better performances than the other two schemes. {To avoid the machine error of too little computational time, we take the average time of multi-calculations as the final  CPU time in Examples \ref{Example:Burgers1DTestOrder}-\ref{Example:Burgers2DTestOrder}.} {To have a fair comparison, the WENO-ZQ scheme uses a true two-dimensional reconstruction as in \cite{trzq} instead of the dimensional-by-dimensional approach  \cite{ZhuJunQiuJX2017fvWENOZQ} in the two-dimensional case, as the HWENO-U scheme uses a true two-dimensional reconstruction too. Differently, the WENO-ZQ scheme uses a wider stencil to reconstruct a bivariate quartic polynomial.}

\begin{example}\label{Example:Burgers1DTestOrder}
		We solve the one-dimensional nonlinear Burgers' equation
		\begin{equation*}
			u_t+(\frac{u^2}{2})_x=0,~0<x<2,
		\end{equation*}
		with periodic boundary conditions up to the time $T=0.5/\pi$ when the solution is still smooth. The initial condition is $u(x,0)= 0.5+\sin(\pi x)$. The numerical errors and CPU time of the HWENO-U, HWENO-M, and WENO-ZQ schemes are presented in Table \ref{sec3:1dBurgerTest}, which shows the schemes all achieve the fifth-order accuracy. More explicitly, with denser meshes (e.g., $\ge 200$), the CPU time ratio of HWENO-U/WENO-ZQ is about 1.739, whereas the $L^\infty$ error ratio is around 1/12.510, and the CPU time ratio of HWENO-M/WENO-ZQ is around 2.528, but the $L^\infty$ error ratio is almost 1/11.529. Since the product of the CPU time ratio and error ratio is less than 1, it shows that the HWENO schemes are more precise than the WENO-ZQ scheme at the same CPU cost.  More intuitively, we can see it from Fig. \ref{sec3:1dCPU_LinftyL1_1dBurgers} that the HWENO-U scheme is more efficient than the HWENO-M and WENO-ZQ schemes. {Besides, we also present the numerical errors and orders of the first-order moments for the HWENO-U scheme in Table \ref{sec3:1dBurgersTest_v}, which demonstrates that modifying the first-order moments in the time level does not destroy their final accuracy.}
		\begin{table}[t]
			\centering
			\fontsize{10.5}{13}\selectfont
			\begin{threeparttable}
				\caption{Example \ref{Example:Burgers1DTestOrder}. One-dimensional Burgers' equation: $L^\infty$ and $L^1$ errors, orders and CPU time of the HWENO-U, HWENO-M, and WENO-ZQ schemes. }
				\label{sec3:1dBurgerTest}
				{\begin{tabular}{ccccccccccccccc}
						\toprule
						\multirow{1}{*}{Meshes}
						& ${L^\infty}$ error&Order&${L^1}$ error&Order &CPU\cr
						\midrule
						{HWENO-U}&&&&&&\\
						40& 	2.61E-04&   $- $&    3.74E-05&   $- $	&	6.53E-04	\\
						80& 	2.16E-06&   6.96&    2.44E-07&   7.26	&	4.51E-03	\\
						120&	2.99E-07&   4.87&    2.51E-08&   5.61	&	1.58E-02	\\
						160&	7.13E-08&   4.99&    5.43E-09&   5.32	&	4.00E-02	\\
						200&	2.32E-08&   5.03&    1.76E-09&   5.05	&	8.40E-02	\\
						240&	9.28E-09&   5.03&    7.01E-10&   5.06	&	1.56E-01	\\
						\midrule
						{HWENO-M}&&&&&&\\
						40& 	6.17E-05&   $- $&    6.49E-06&   $- $	&	9.23E-04	\\
						80&	    2.31E-06&   4.74&    1.77E-07&   5.20	&	6.50E-03	\\
						120&	3.23E-07&   4.86&    2.42E-08&   4.90	&	2.29E-02	\\
						160&	7.70E-08&   4.98&    5.70E-09&   5.03	&	5.83E-02	\\
						200&	2.51E-08&   5.02&    1.88E-09&   4.97	&	1.22E-01	\\
						240&	1.01E-08&   5.02&    7.57E-10&   5.00	&	2.27E-01	\\
						\midrule
						{WENO-ZQ}&&&&&&\\
						40& 	4.79E-04&   $- $&    4.67E-05&   $- $	&	3.98E-04	\\
						80&	    2.45E-05&   4.29&    1.96E-06&   4.57	&	2.64E-03	\\
						120&	3.56E-06&   4.75&    2.76E-07&   4.84	&	9.15E-03	\\
						160&	8.71E-07&   4.90&    6.65E-08&   4.95	&	2.30E-02	\\
						200&	2.88E-07&   4.96&    2.22E-08&   4.93	&	4.83E-02	\\
						240&	1.17E-07&   4.93&    8.95E-09&   4.98	&	8.97E-02	\\
						\bottomrule
				\end{tabular}}
			\end{threeparttable}
		\end{table}
		\begin{figure}[!ht]
			\centering
			\subfigure{\includegraphics[width=8cm,angle=0]{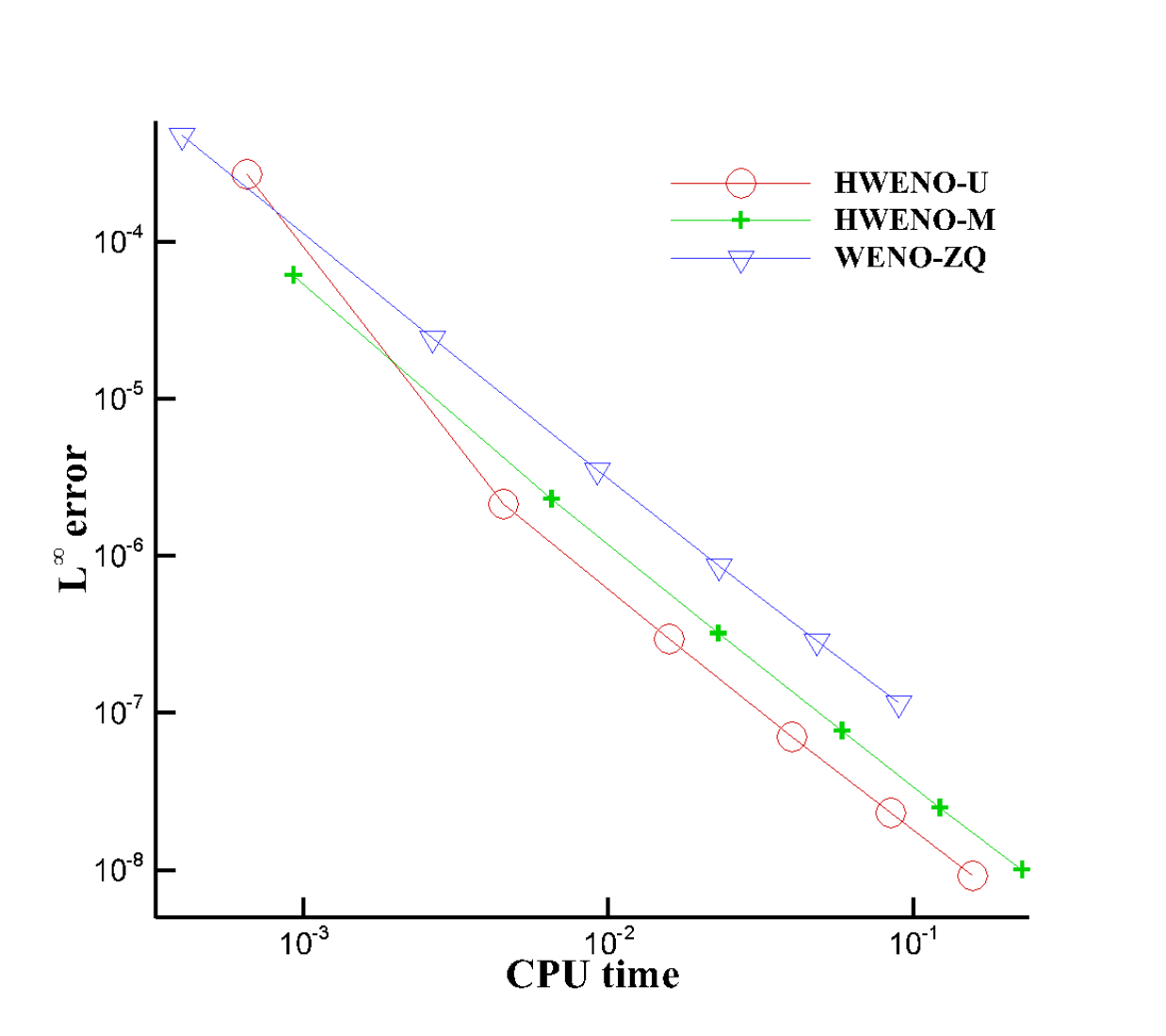}}
			\subfigure{\includegraphics[width=8cm,angle=0]{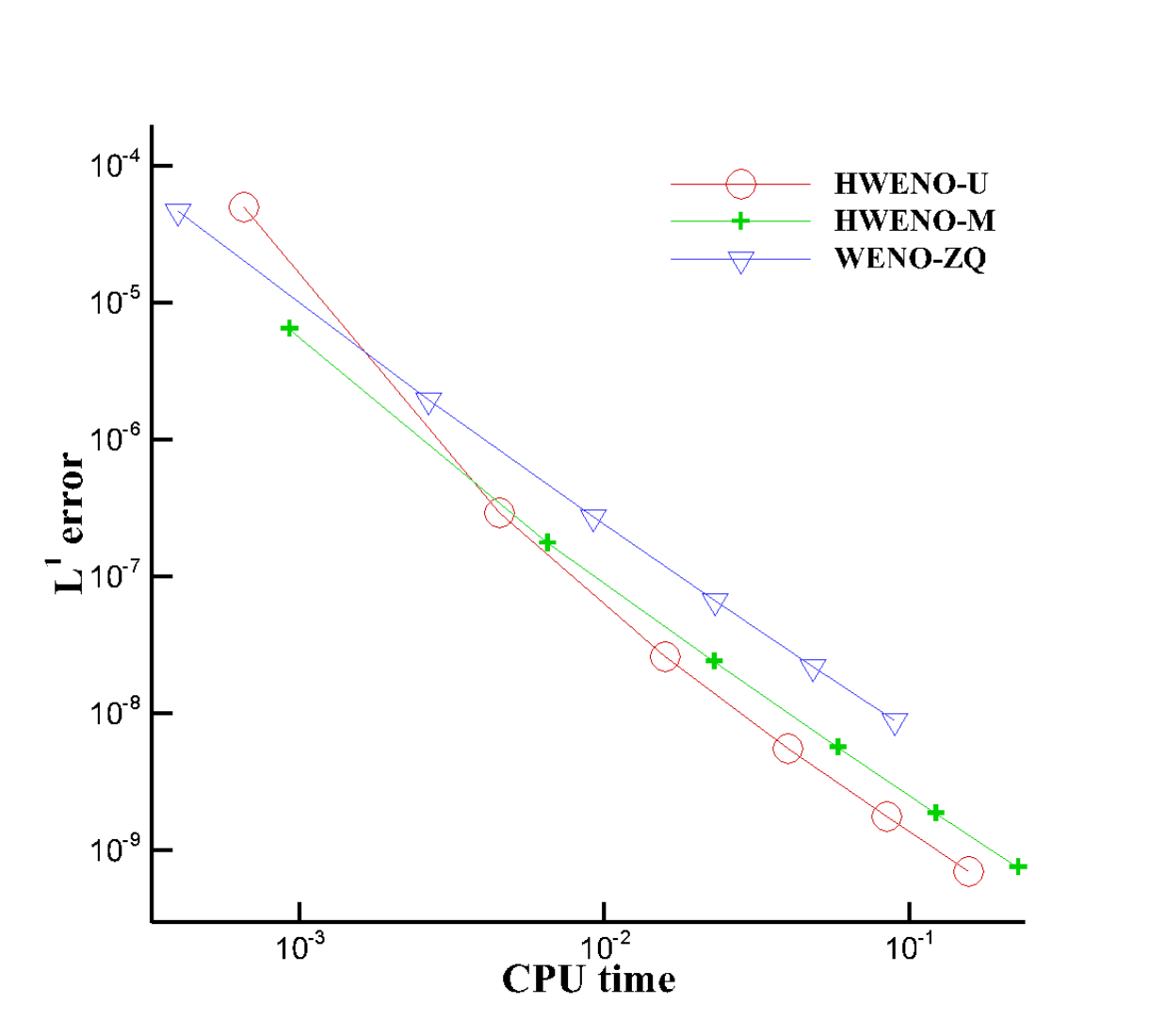}}
			\caption{ Comparison of $L^\infty$, $L^1$ errors and CPU time for Example \ref{Example:Burgers1DTestOrder}. }
			\label{sec3:1dCPU_LinftyL1_1dBurgers}
		\end{figure}
		\begin{table}[!ht]
			\centering
			\fontsize{10.5}{13}\selectfont
			\begin{threeparttable}
				\caption{Example \ref{Example:Burgers1DTestOrder}. $L^\infty$ and $L^1$ errors, and orders of the HWENO-U scheme for the first-order moment. }
				\label{sec3:1dBurgersTest_v}
				{\begin{tabular}{ccccccccccccccc}
						\toprule
						\multirow{1}{*}{Meshes}
						& ${L^\infty}$ error&Order&${L^1}$ error&Order \cr
						\midrule
						40&     3.50E-05&   $- $&    2.04E-06&   $- $\\
						80&	    9.33E-07&   5.23&    5.76E-08&   5.15\\
						120&	1.18E-07&   5.11&    7.37E-09&   5.07\\
						160&	2.72E-08&   5.10&    1.69E-09&   5.12\\
						200&	8.22E-09&   5.36&    5.27E-10&   5.21\\
						240&	3.37E-09&   4.89&    2.13E-10&   4.96\\
						\bottomrule
				\end{tabular}}
			\end{threeparttable}
		\end{table}
\end{example}
\begin{example}\label{Example:Euler1DTestOrder}
		We solve one-dimensional compressible Euler equations		
		\begin{equation*}
			\frac{\partial}{\partial{t}}\begin{bmatrix}	\rho \\ \rho\mu\\E \end{bmatrix} +
			\frac{\partial}{\partial{x}}\begin{bmatrix}	\rho\mu \\ \rho\mu^2+p\\ \mu(E+p) \end{bmatrix} = 0,
		\end{equation*}
		where $\rho$ is the density,  $\mu$ is the velocity, $E$ is the  total energy and $p$ is the pressure. The initial condition is $(\rho,\mu,p,\gamma)=(1+0.2\sin(\pi x),1,1,1.4)$ on the domain $[0,2]$ with periodic boundary conditions. The final time is $T=2$, and the exact solutions are $(\rho,\mu,p)=(1+0.2\sin(\pi (x-T)),1,1)$. The numerical errors and CPU time of the HWENO-U, HWENO-M, and WENO-ZQ schemes are presented in Table \ref{sec3:1dEulerTest}, illustrating the schemes all achieve the fifth-order accuracy.
		More explicitly, on the denser meshes (e.g., $\ge 200$), the CPU time ratio of HWENO-U/WENO-ZQ is about 1.695, whereas the $L^\infty$ error ratio is around 1/10.905, and the CPU time ratio of HWENO-M/WENO-ZQ is around 2.488, but the $L^\infty$ error ratio is almost 1/11.897.  These data demonstrates that at the same CPU cost, the HWENO-U scheme is more accurate than the HWENO-M and WENO-ZQ schemes, which also can be more intuitively observed from Fig. \ref{sec3:1dCPU_LinftyL1_1dEuler}.
		\begin{table}[t]
			\centering
			\fontsize{10.5}{13}\selectfont
			\begin{threeparttable}
				\caption{Example \ref{Example:Euler1DTestOrder}. One-dimensional Euler equations: $L^\infty$ and $L^1$ errors, orders and CPU time of the HWENO-U, HWENO-M, and WENO-ZQ schemes. }
				\label{sec3:1dEulerTest}
				{\begin{tabular}{cccccccccccccc}
						\toprule
						{Meshes} & ${L^\infty}$ error&Order&${L^1}$ error&Order &CPU \cr
						\midrule
						{HWENO-U}\\
						40& 	4.48E-06&   $- $&    7.98E-07&   $- $	&	5.39E-02	\\
						80& 	3.78E-08&   6.89&    8.08E-09&   6.63	&	3.89E-01	\\
						120&	2.92E-09&   6.31&    1.04E-09&   5.05	&	1.36E+00	\\
						160&	5.41E-10&   5.86&    2.46E-10&   5.01	&	3.44E+00	\\
						200&	1.57E-10&   5.56&    8.07E-11&   5.01	&	7.22E+00	\\
						240&	5.87E-11&   5.38&    3.24E-11&   5.00	&	1.34E+01	\\
						\midrule
						{HWENO-M}\\
						40& 	4.45E-07&   $- $&    2.54E-07&   $- $	&	7.74E-02	\\
						80& 	1.26E-08&   5.14&    7.89E-09&   5.01	&	5.58E-01	\\
						120&	1.64E-09&   5.04&    1.04E-09&   5.01	&	1.98E+00	\\
						160&	3.87E-10&   5.02&    2.46E-10&   5.00	&	5.04E+00	\\
						200&	1.26E-10&   5.01&    8.05E-11&   5.00	&	1.06E+01	\\
						240&	5.08E-11&   5.00&    3.23E-11&   5.00	&	1.96E+01	\\
						\midrule
						{WENO-ZQ}\\
						40& 	4.71E-06&   $- $&    2.97E-06&   $- $ 	&	3.28E-02	\\
						80& 	1.47E-07&   5.00&    9.34E-08&   4.99	&	2.29E-01	\\
						120&	1.93E-08&   5.00&    1.23E-08&   5.00	&	8.02E-01	\\
						160&	4.59E-09&   5.00&    2.92E-09&   5.00	&	2.03E+00	\\
						200&	1.50E-09&   5.00&    9.57E-10&   5.00	&	4.26E+00	\\
						240&	6.04E-10&   5.00&    3.85E-10&   5.00	&	7.88E+00	\\
						\bottomrule	
				\end{tabular}}
			\end{threeparttable}
		\end{table}
		\begin{figure}[!ht]
			\centering
			\subfigure{\includegraphics[width=8cm,angle=0]{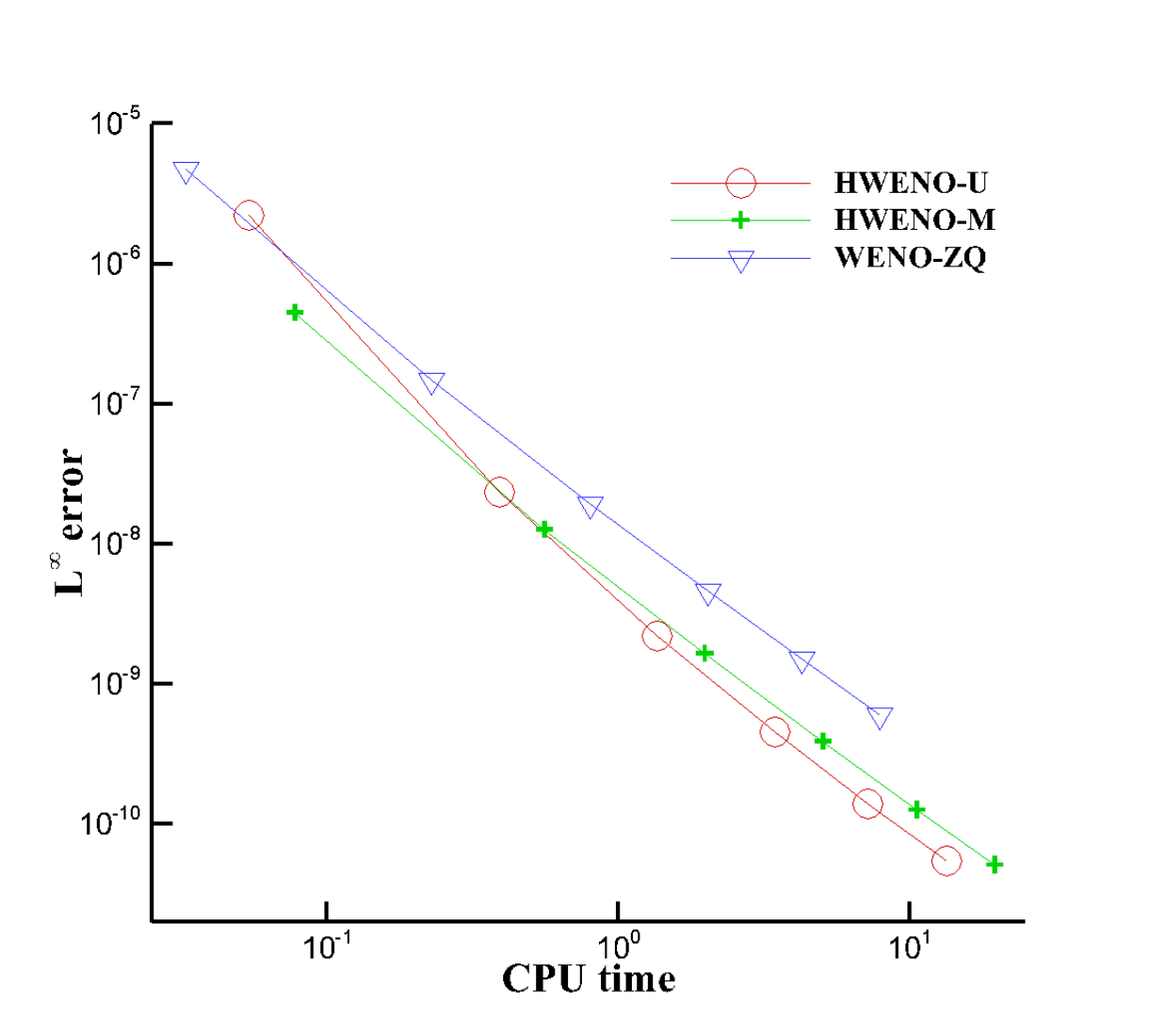}}
			\subfigure{\includegraphics[width=8cm,angle=0]{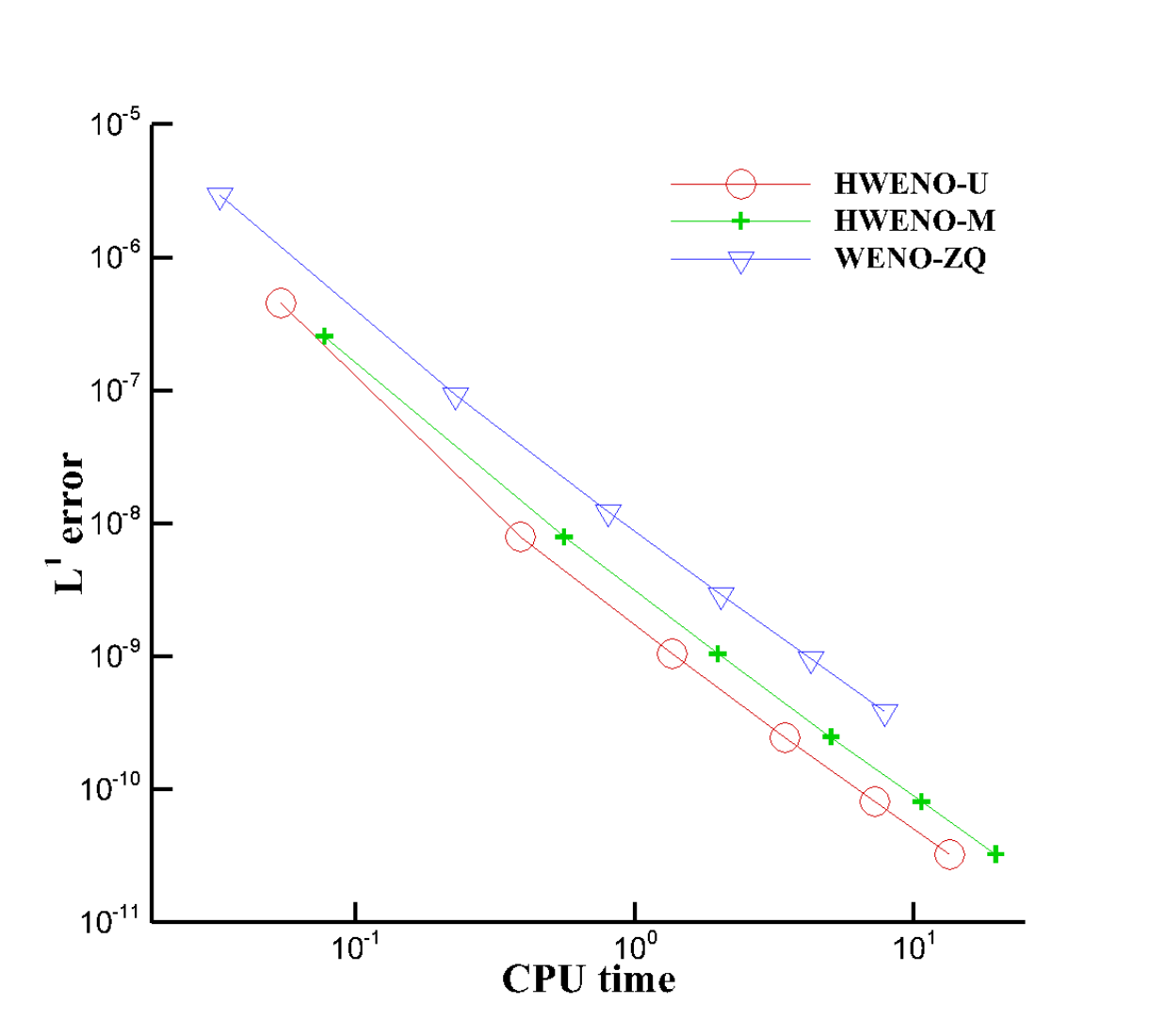}}
			\caption{ Comparison of $L^\infty$, $L^1$ errors and CPU time for Example \ref{Example:Euler1DTestOrder}. }
			\label{sec3:1dCPU_LinftyL1_1dEuler}
		\end{figure}
\end{example}
\begin{example}\label{Example:Burgers2DTestOrder}
		We consider two-dimensional nonlinear Burgers' equation
		\begin{equation*}
		u_t+(\frac{u^2}{2})_x+(\frac{u^2}{2})_y=0,~0<x<4,~0<y<4,
		\end{equation*}
		with the initial condition $u(x,y,0)= 0.5+\sin(\pi (x+y)/2)$ and periodic boundary conditions in $x$ and $y$ directions. Up to the final time $T=0.5/\pi$, the solution is still smooth. The $L^\infty$ and $L^1$ norms of numerical errors and CPU time of the HWENO-U, HWENO-M and WENO-ZQ schemes are shown in Table \ref{sec3:2dBurgersTest},  illustrating  that the three schemes achieve the fifth-order accuracy. {More explicitly, with denser meshes (e.g., $\ge 200$), the CPU time ratio of HWENO-U/WENO-ZQ is about 1.242, whereas the $L^\infty$ error ratio is around 1/64.085, and the CPU time ratio of HWENO-M/WENO-ZQ is around 1.534, but the $L^\infty$ error ratio is almost 1/74.196. This data shows that the HWENO-U scheme is more precise than the HWENO-M and WENO-ZQ schemes at the same CPU cost, which can be intuitively seen from Fig. \ref{sec3:1dCPU_LinftyL1_2dBurgers}.} With the mesh gets denser, we can intuitively observe that the HWENO-U and HWENO-M schemes have similar numerical errors, but the HWENO-U scheme has slightly less computational time. {Besides, the numerical errors and orders of the first-order moments in the $x$ and $y$ directions for the HWENO-U scheme  are presented in Table \ref{sec3:2dBurgersTest_v}, in which  the first-order moments also have the fifth-order accuracy as that in the one-dimensional case. {Note that the errors and orders of the first-order moments in the $x$ and $y$ directions are identical because of the symmetry  solution.}}
		\begin{table}[t]
			\centering
			\fontsize{9.5}{13}\selectfont
			\begin{threeparttable}
				\caption{Example \ref{Example:Burgers2DTestOrder}. Two-dimensional Burgers' equation: $L^\infty$ and $L^1$ errors, orders and CPU time of the HWENO-U, HWENO-M  and WENO-ZQ schemes. }
				\label{sec3:2dBurgersTest}
				{\begin{tabular}{ccccccccccccc}
						\toprule
						{Meshes} & ${L^\infty}$ error&Order&${L^1}$ error&Order &CPU \cr
						\midrule
						{HWENO-U}\\
						$40\times40$&	1.83E-03&   $- $&    1.53E-04&   $- $	&	2.81E-01	\\
						$80\times80$&	2.61E-05&   6.13&    2.47E-06&   5.96	&	3.75E+00	\\
						$120\times120$&	1.43E-06&   7.17&    7.60E-08&   8.59	&	1.89E+01	\\
						$160\times160$&	1.21E-07&   8.57&    9.68E-09&   7.16	&	6.23E+01	\\
						$200\times200$&	3.12E-08&   6.09&    2.34E-09&   6.35	&	1.63E+02	\\
						$240\times240$&	9.52E-09&   6.50&    8.67E-10&   5.46	&	3.69E+02	\\						
						\midrule
						{HWENO-M}\\
						$40\times40$&	5.67E-05&   $- $&    5.80E-06&   $- $	&	3.43E-01	\\
						$80\times80$&	2.27E-06&   4.64&    1.76E-07&   5.04	&	4.61E+00	\\
						$120\times120$&	3.04E-07&   4.96&    2.26E-08&   5.06	&	2.33E+01	\\
						$160\times160$&	7.18E-08&   5.02&    5.46E-09&   4.94	&	7.70E+01	\\
						$200\times200$&	2.38E-08&   4.95&    1.78E-09&   5.03	&	2.01E+02	\\
						$240\times240$&	9.52E-09&   5.02&    7.13E-10&   5.01	&	4.56E+02	\\						
						\midrule
						{WENO-ZQ}\\
						$40\times40$&	1.82E-03&    ---&    2.16E-04&    -- 	&	2.37E-01	\\
						$80\times80$&	1.24E-04&   3.87&    1.07E-05&   4.33	&	3.12E+00	\\
						$120\times120$&	1.96E-05&   4.56&    1.58E-06&   4.72	&	1.56E+01	\\
						$160\times160$&	5.09E-06&   4.68&    3.98E-07&   4.80	&	5.16E+01	\\
						$200\times200$&	1.74E-06&   4.80&    1.33E-07&   4.90	&	1.33E+02	\\
						$240\times240$&	7.17E-07&   4.86&    5.49E-08&   4.86	&	2.93E+02	\\
						\bottomrule
				\end{tabular}}
			\end{threeparttable}
		\end{table}
		\begin{figure}[ht]
			\centering
			\subfigure{\includegraphics[width=8cm,angle=0]{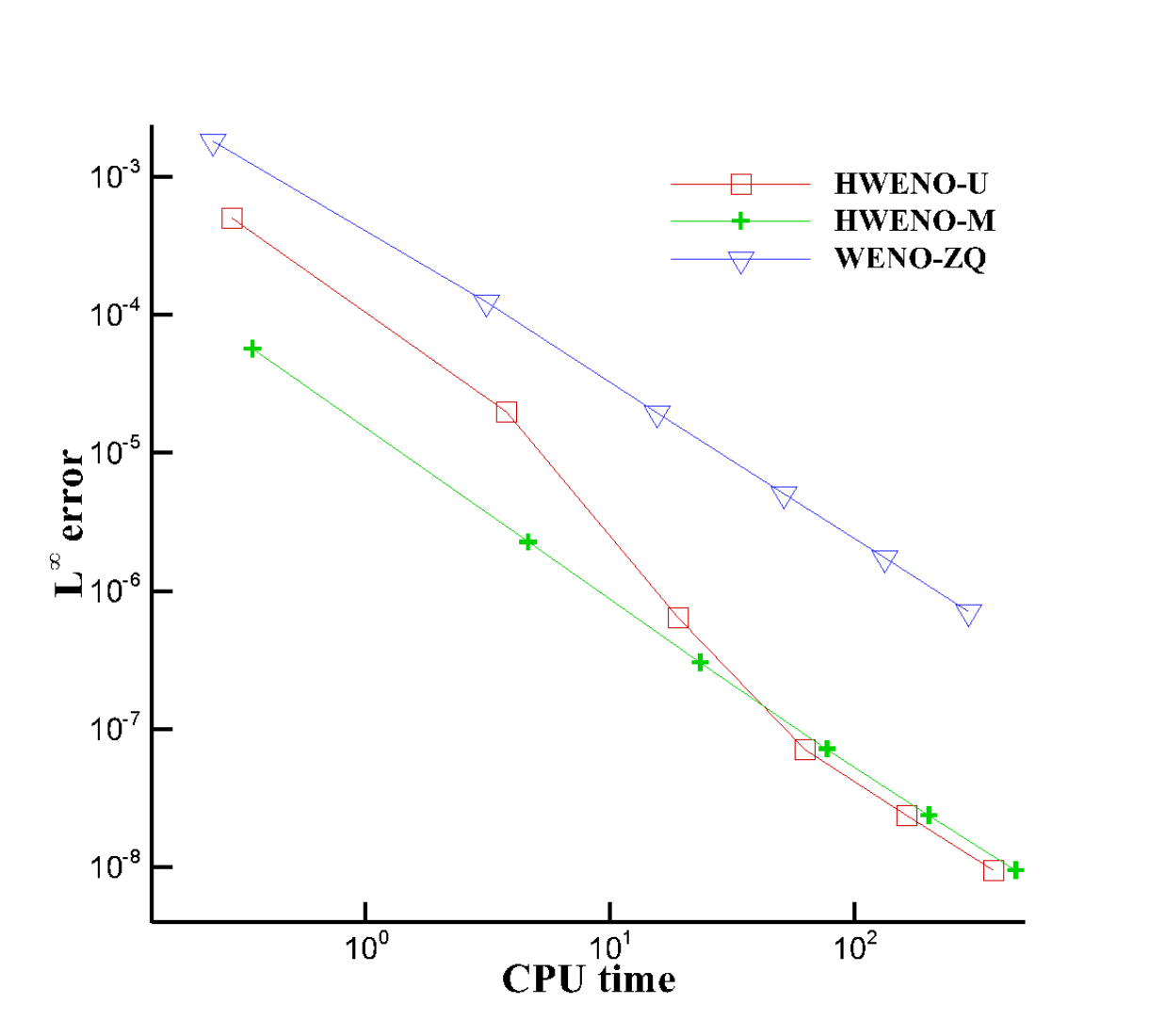}}
			\subfigure{\includegraphics[width=8cm,angle=0]{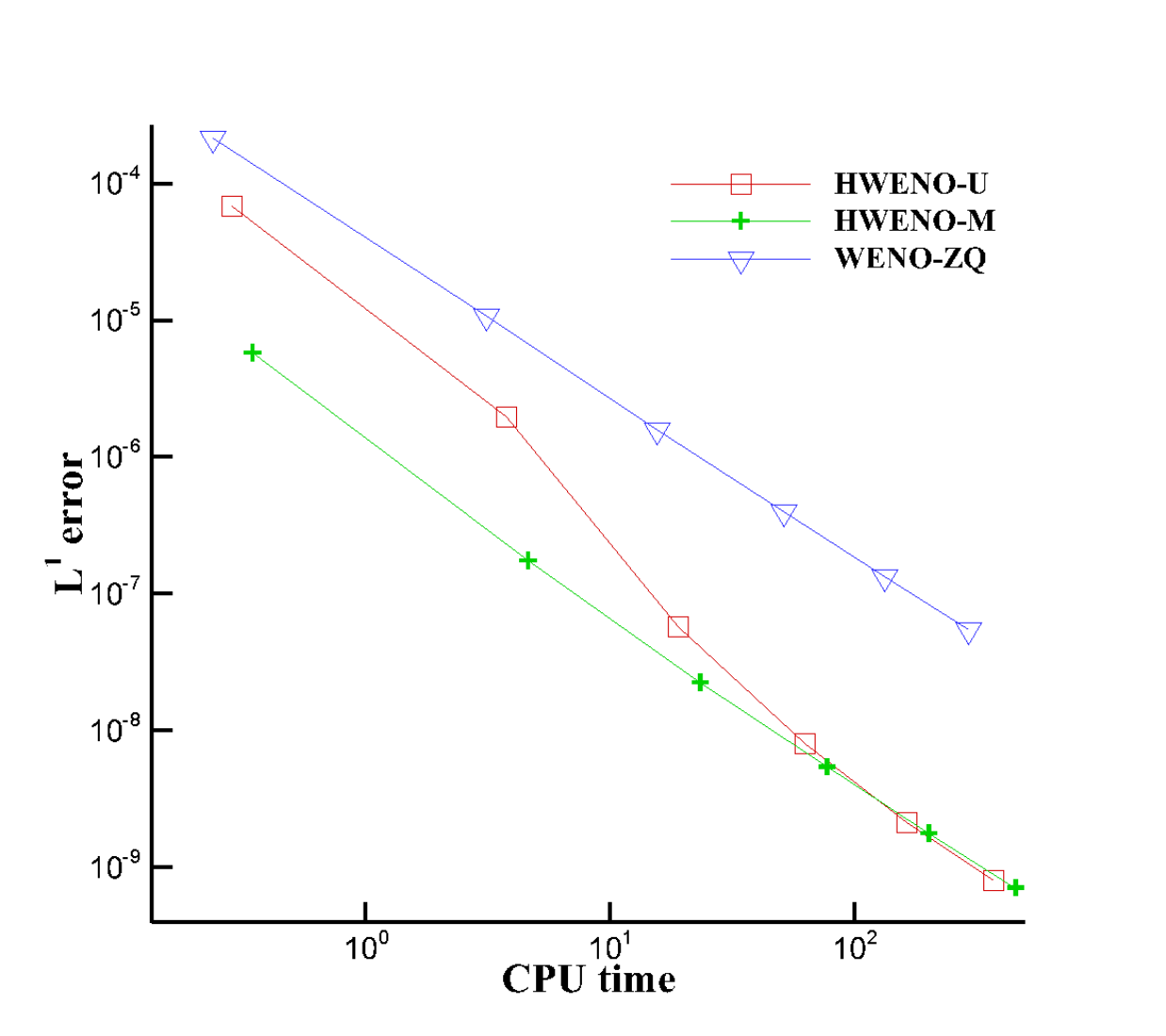}}
			\caption{ Comparison of $L^\infty$, $L^1$ errors and CPU time for Example \ref{Example:Burgers2DTestOrder}. }
			\label{sec3:1dCPU_LinftyL1_2dBurgers}
		\end{figure}
		\begin{table}[ht]
			\centering
			\fontsize{9.5}{13}\selectfont
			{\begin{threeparttable}
				\caption{{Example \ref{Example:Burgers2DTestOrder}. Two-dimensional Burgers' equation: $L^\infty$ and $L^1$ errors, orders of the HWENO-U scheme for the first-order moments in the $x$ and $y$ directions.}}
				\label{sec3:2dBurgersTest_v}
				{\begin{tabular}{ccccccccccccc}
						\toprule
						\multirow{2}{*}{Meshes}&
						\multicolumn{2}{l}{$x$-direction} &&& \multicolumn{2}{l}{$y$-direction}\cr
						\cmidrule(l){2-5} \cmidrule(l){6-9}
						& ${L^\infty}$ error&Order&${L^1}$ error&Order &  ${L^\infty}$ error&Order&${L^1}$ error&Order \cr
						\midrule
						$40\times40$&	9.98E-04&   $- $&    6.29E-05&   $- $		&	9.98E-04&   $- $&    6.29E-05&   $- $	\\
						$80\times80$&	4.40E-06&   7.83&    4.49E-07&   7.13		&	4.40E-06&   7.83&    4.49E-07&   7.13	\\
						$120\times120$&	1.59E-07&   8.19&    1.52E-08&   8.35		&	1.59E-07&   8.19&    1.52E-08&   8.35	\\
						$160\times160$&	2.64E-08&   6.25&    2.23E-09&   6.68		&	2.64E-08&   6.25&    2.23E-09&   6.68	\\
						$200\times200$&	8.38E-09&   5.14&    6.00E-10&   5.88		&	8.38E-09&   5.14&    6.00E-10&   5.88	\\
						$240\times240$&	3.31E-09&   5.09&    2.33E-10&   5.20		&	3.31E-09&   5.09&    2.33E-10&   5.20	\\
						\bottomrule
				\end{tabular}}
			\end{threeparttable}}
		\end{table}
\end{example}

\begin{example}\label{Example:Euler2DTestOrder}
		We solve two-dimensional compressible Euler equations		
		\begin{equation*}
			\frac{\partial}{\partial{t}}\begin{bmatrix}	\rho \\ \rho\mu \\ \rho\nu  \\E \end{bmatrix} +
			\frac{\partial}{\partial{x}}\begin{bmatrix}	\rho\mu \\ \rho\mu^2+p\\ \rho\mu\nu\\ \mu(E+p) \end{bmatrix}+
			\frac{\partial}{\partial{y}}\begin{bmatrix}	\rho\nu \\ \rho\mu\nu \\ \rho\nu^2+p\\ \nu(E+p) \end{bmatrix} = 0,
		\end{equation*}
		where $\rho$ is the density,  $\mu$ and $\nu$  the velocity in $x$ and $y$ directions respectively, $E$ is the total energy and $p$ is the pressure. The initial condition is $(\rho,\mu,\nu,p,\gamma)=(1+0.2\sin(\pi(x+y)),1,1,1,1.4)$ and the computational domain is $[0,4]\times[0,4]$ with periodic boundary conditions in all directions. We compute the solution up to time $T=2$, and the exact solutions are $(\rho,\mu,\nu,p)=(1+0.2\sin(\pi (x+y-2T)),1,1,1)$. The $L^\infty$ and $L^1$ errors are presented in Table \ref{sec3:2dEulerTest}, showing that the three schemes achieve the fifth-order accuracy. More explicitly, with denser meshes (e.g., $\ge 200$), the CPU time ratio of HWENO-U/WENO-ZQ is about 1.252, whereas the $L^\infty$ error ratio is around 1/14.443, and the CPU time ratio of HWENO-M/WENO-ZQ is around 1.378, but the $L^\infty$ error ratio is almost 1/78.888. This data shows that the HWENO-U scheme is more precise than the HWENO-M and WENO-ZQ schemes at the same CPU cost, which can be intuitively seen from Fig. \ref{sec3:1dCPU_LinftyL1_2dBurgers}. Compared to the WENO-ZQ scheme, the HWENO-U and HWENO-M schemes require the computation of two extra first-order moment equations, yet their CPU costs only increase by no more than 40\% due to the repeated utilization of numerical fluxes on the boundary in the zeroth- and first-order moment equations. Overall, both the HWENO-U and HWENO-M schemes demonstrate superior computational efficiency than the WENO-ZQ scheme. It is worth noting that despite using unified stencils throughout the entire procedures, the HWENO-U and HWENO-M schemes still have similar numerical errors and comparable CPU time. This can be attributed to the fact that the HWENO-M scheme modifies the first-order moments in a dimension-by-dimensional manner, resulting in computational cost savings, particularly for high-dimensional systems. However, extending this dimension-by-dimensional approach to unstructured meshes is not straightforward. In contrast, the framework of the HWENO-U scheme is specifically designed to be well-suited for unstructured cases, and the relevant researches  are ongoing.
		\begin{table}[ht]
			\centering
			\fontsize{9.0}{13}\selectfont
			\begin{threeparttable}
				\caption{Example \ref{Example:Euler2DTestOrder}. Two-dimensional Euler equations: $L^\infty$ and $L^1$ errors, orders and CPU time of the HWENO-U, HWENO-M  and WENO-ZQ schemes. }
				\label{sec3:2dEulerTest}
				{\begin{tabular}{ccccccccccccc}
						\toprule
						{Meshes} & ${L^\infty}$ error&Order&${L^1}$ error&Order &CPU \cr
						\midrule
						{HWENO-U}\\
						$40\times40$&	1.33E-04&   $- $&    2.23E-05&   $- $	&	1.24E+02	\\
						$80\times80$&	8.41E-07&   7.31&    8.17E-08&   8.09	&	1.72E+03	\\
						$120\times120$&	4.26E-08&   7.36&    3.82E-09&   7.55	&	1.02E+04	\\
						$160\times160$&	5.28E-09&   7.26&    5.98E-10&   6.45	&	3.65E+04	\\
						$200\times200$&	1.08E-09&   7.11&    1.67E-10&   5.72	&	9.68E+04	\\
						$240\times240$&	3.06E-10&   6.92&    6.45E-11&   5.21	&	2.14E+05	\\
						\midrule
						{HWENO-M}\\
						$40\times40$&	8.85E-07&   $- $&    5.01E-07&   $- $	&	1.34E+02	\\
						$80\times80$&	2.52E-08&   5.14&    1.57E-08&   5.00	&	1.86E+03	\\
						$120\times120$&	3.25E-09&   5.04&    2.06E-09&   5.00	&	1.09E+04	\\
						$160\times160$&	7.68E-10&   5.02&    4.89E-10&   5.00	&	3.96E+04	\\
						$200\times200$&	2.52E-10&   5.00&    1.60E-10&   5.00	&	1.06E+05	\\
						$240\times240$&	1.01E-10&   4.99&    6.43E-11&   5.00	&	2.37E+05	\\
						\midrule
						{WENO-ZQ}\\
						$40\times40$&	1.69E-04&   $- $&    3.74E-05&   $- $	&	9.86E+01	\\
						$80\times80$&	2.56E-06&   6.04&    1.19E-06&   4.97	&	1.30E+03	\\
						$120\times120$&	2.80E-07&   5.46&    1.57E-07&   4.99	&	7.12E+03	\\
						$160\times160$&	6.25E-08&   5.21&    3.74E-08&   5.00	&	2.54E+04	\\
						$200\times200$&	2.00E-08&   5.12&    1.22E-08&   5.00	&	7.48E+04	\\
						$240\times240$&	7.92E-09&   5.07&    4.92E-09&   5.00	&	1.77E+05	\\
						\bottomrule
				\end{tabular}}
			\end{threeparttable}
		\end{table}
		\begin{figure}[!ht]
			\centering
			\subfigure{\includegraphics[width=8cm,angle=0]{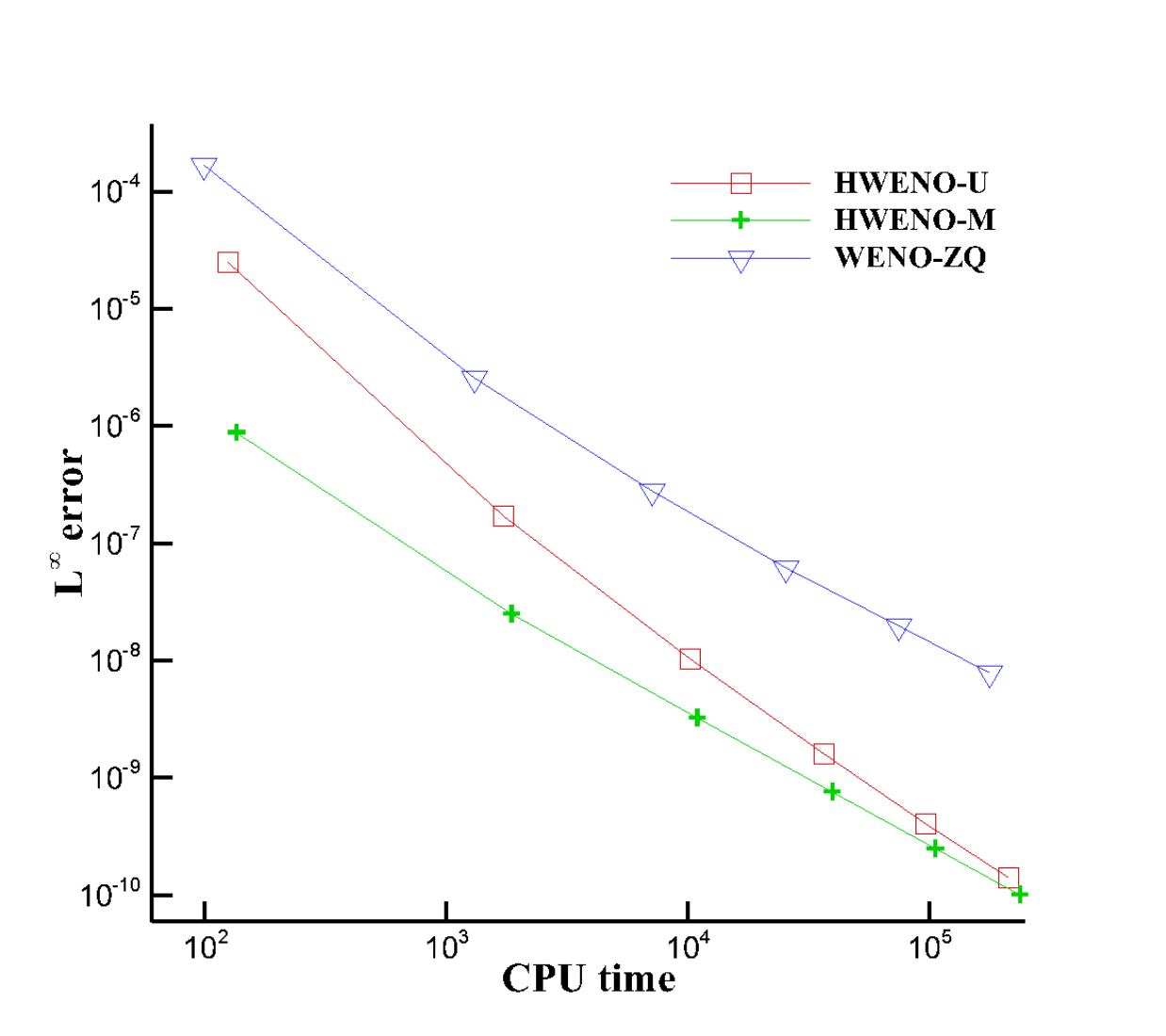}}
			\subfigure{\includegraphics[width=8cm,angle=0]{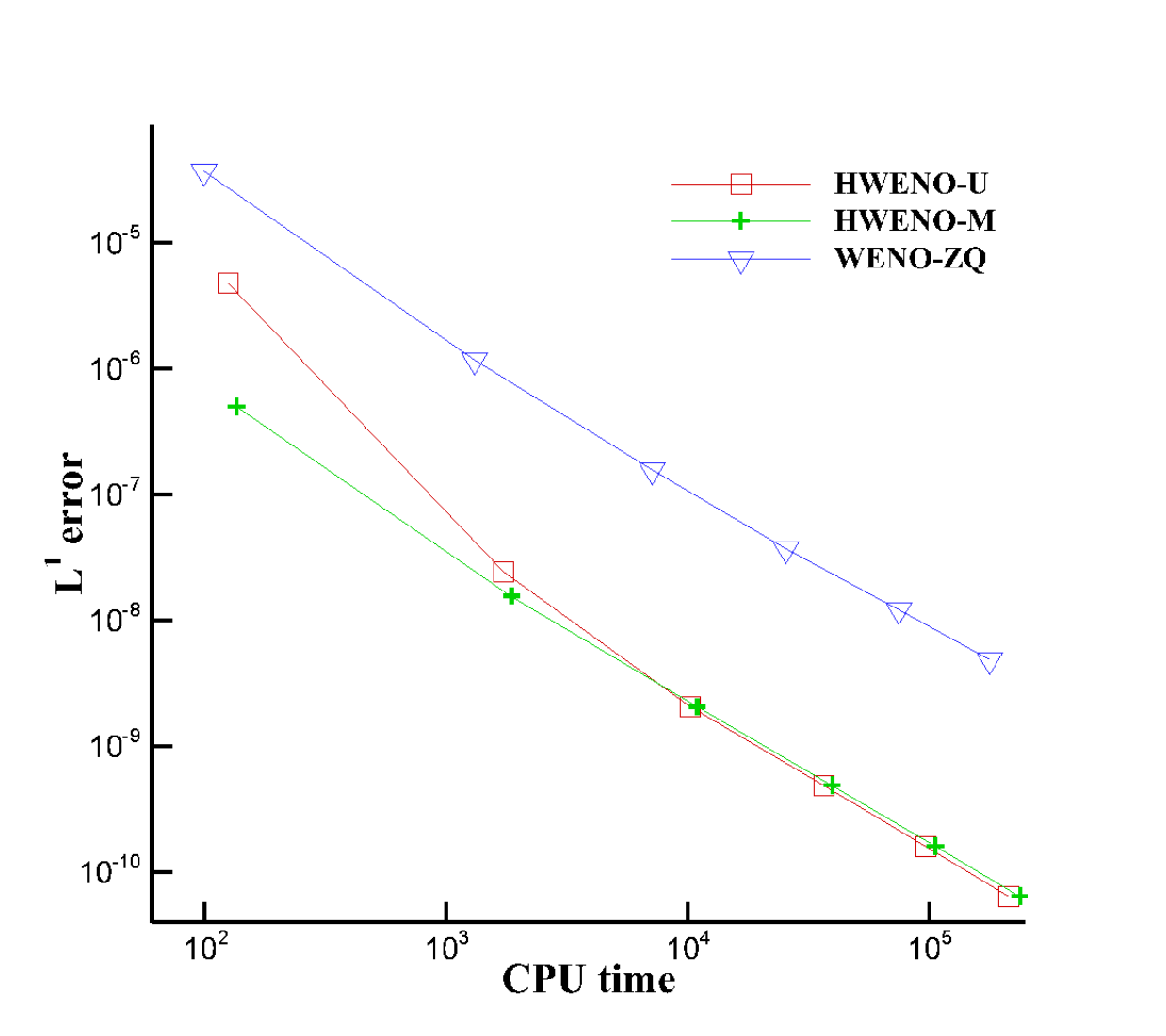}}
			\caption{ Comparison of $L^\infty$, $L^1$ errors and CPU time for Example \ref{Example:Euler2DTestOrder}. }
			\label{sec3:1dCPU_LinftyL1_2dEuler}
		\end{figure}
\end{example}

\subsection{Non-smooth tests}\label{sec3:NonsmoothTests}
In this subsection, we compare the performance of the HWENO-U, HWENO-M, and WENO-ZQ schemes in capturing shocks by simulating some  benchmark and extreme problems.
\begin{example}\label{Example:Lax1D}
\begin{figure}[t]
	\centering
	\subfigure[Density with 100 cells]{\includegraphics[width=7.5cm,angle=0]{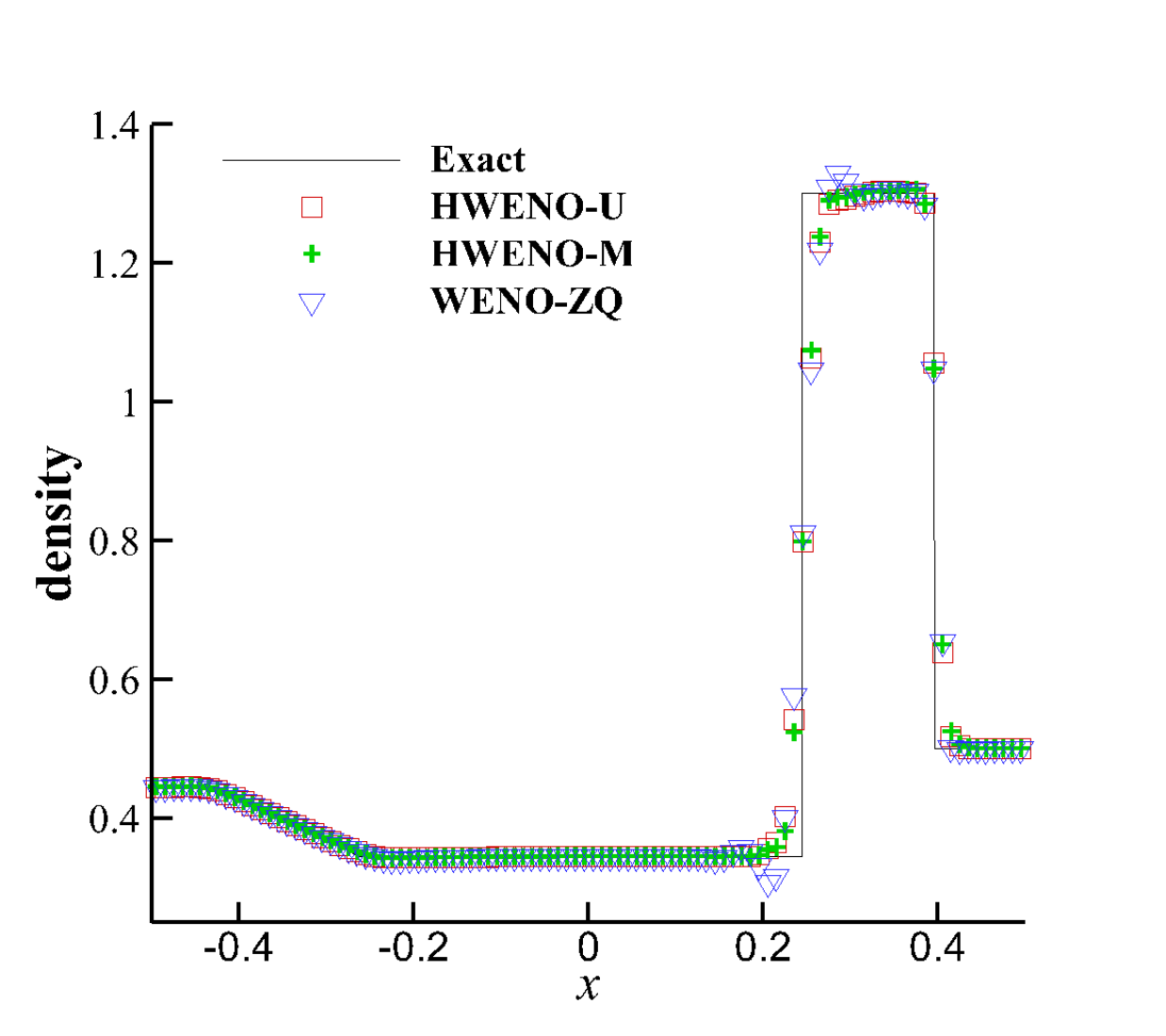}}
	\subfigure[A zoomed-in figure]{\includegraphics[width=7.5cm,angle=0]{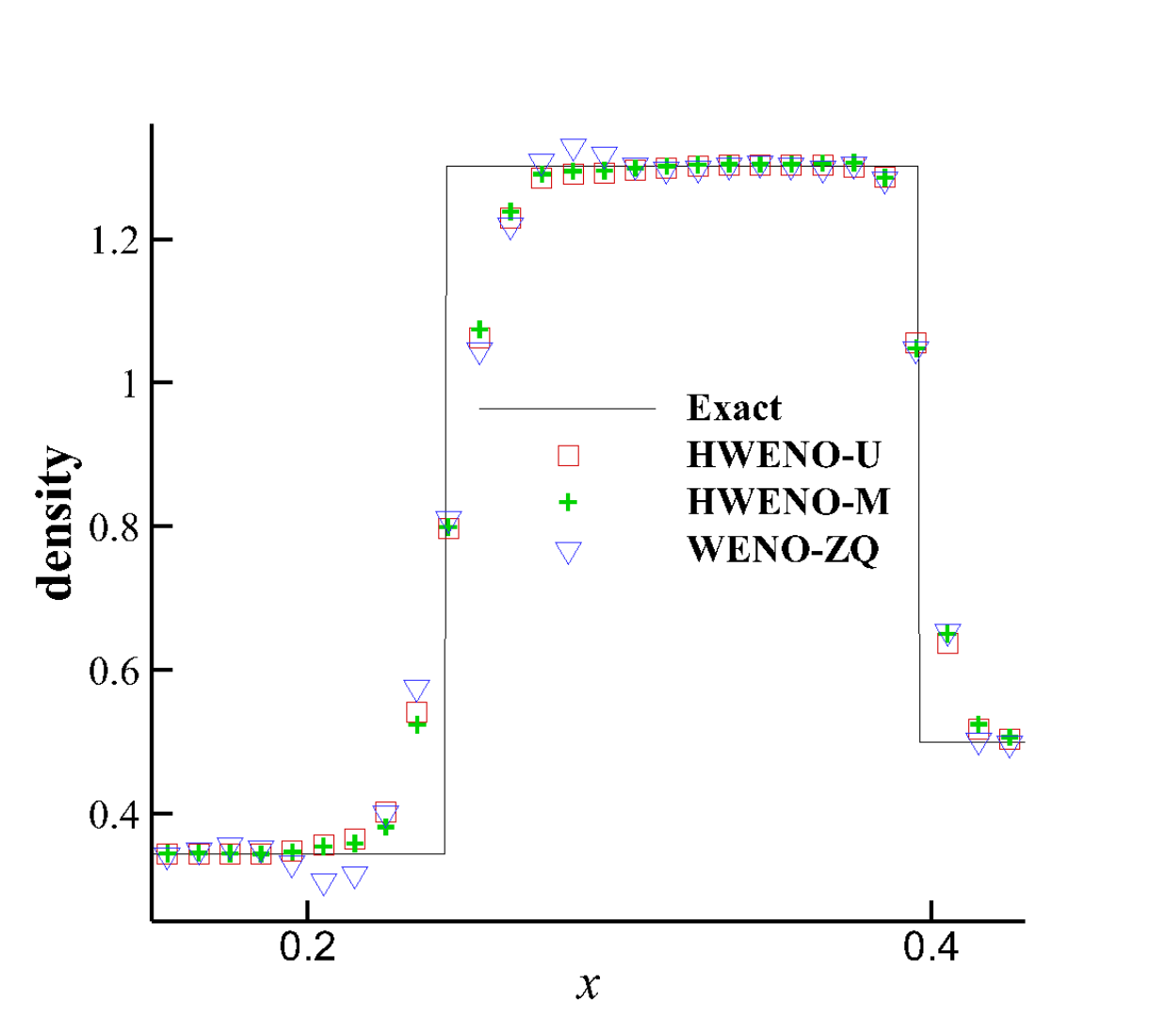}}
	\caption{Example \ref{Example:Lax1D}.
		The results of solution computed by the HWENO-U, HWENO-M and WENO-ZQ schemes.}
	\label{sec3:Lax1D_WENOHWENO}
\end{figure}
\begin{figure}[t]
	\centering
	\subfigure[HWENO-U ($\zeta=10^{-6}$)]
	{\includegraphics[width=7.5cm,angle=0]{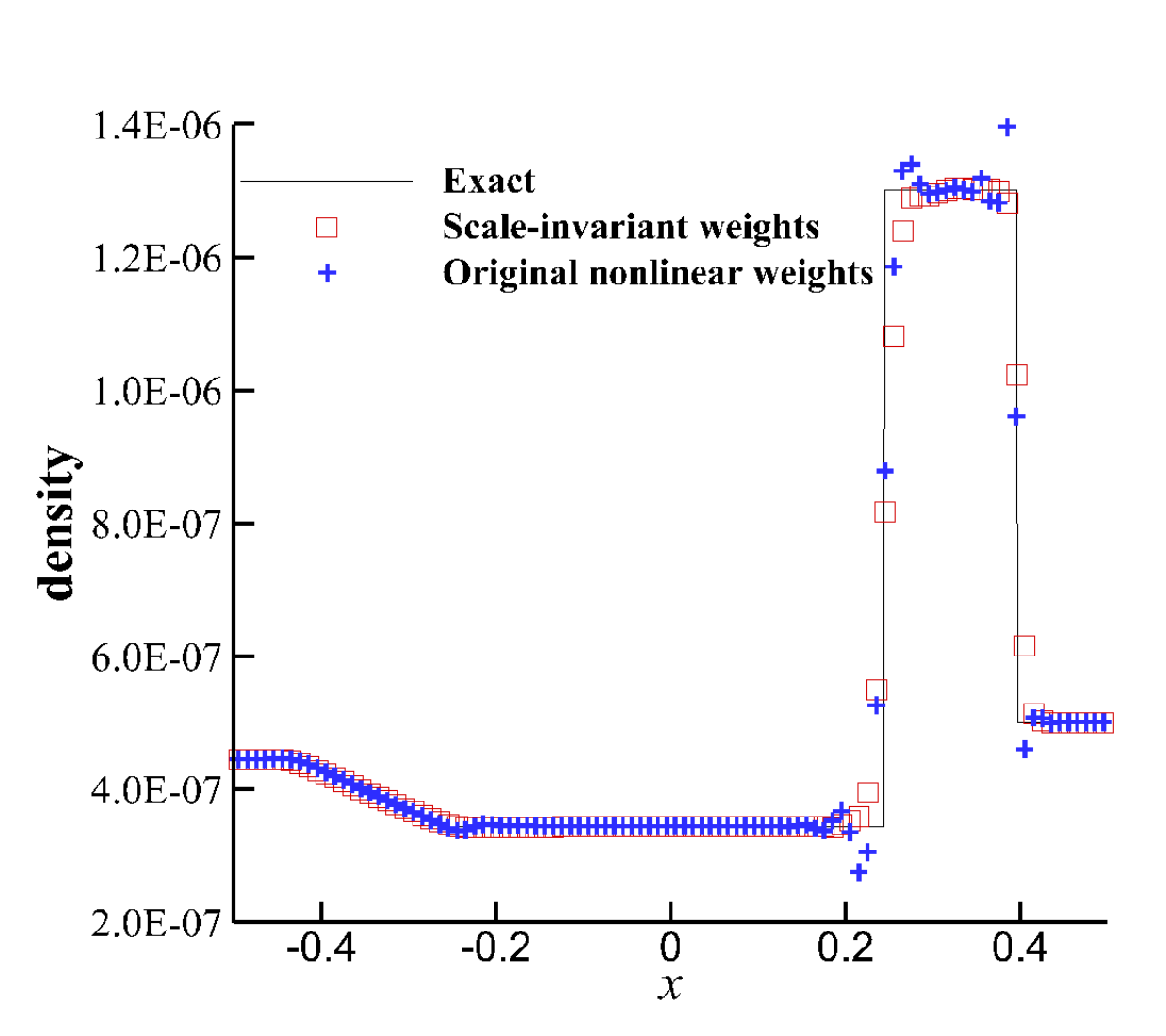}}
	\subfigure[HWENO-U ({$\zeta=10^{6}$})]
	{\includegraphics[width=7.5cm,angle=0]{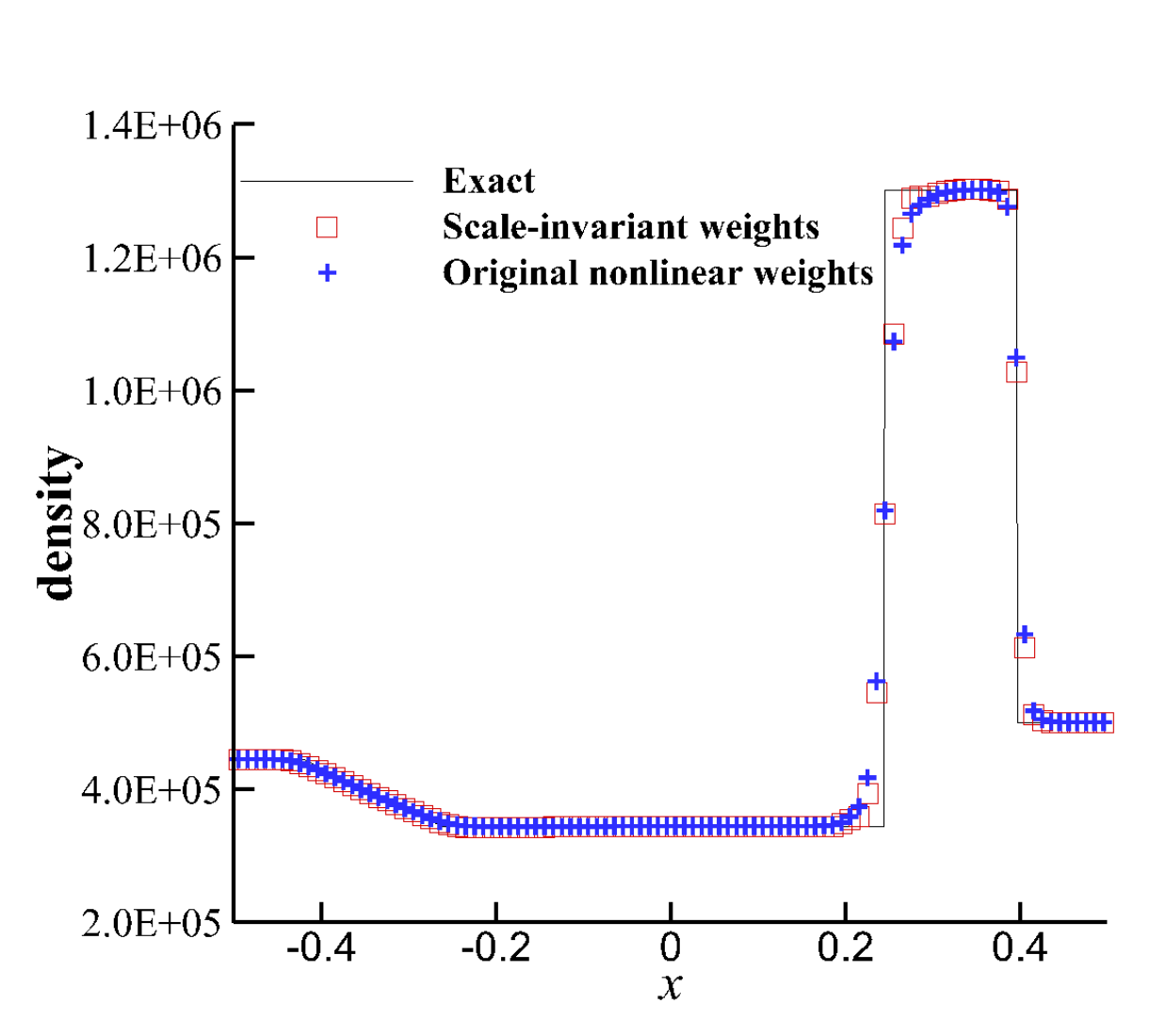}}
	\caption{Example \ref{Example:Lax1D}.
		The results of solution computed by the HWENO-U scheme with the scaled initial conditions using the proposed scale-invariant nonlinear weights \eqref{sec2:non_weight_1d} and the original nonlinear weights \cite{ZhuJunQiuJX2017fvWENOZQ}.}
	\label{sec3:Lax1D_HWENO_U_Scale}
\end{figure}
We solve the Lax problem for one-dimensional Euler equations with the initial conditions:
\begin{equation*}
	(\rho,\mu,p,\gamma)^\mathsf{T}=\begin{cases}
		(0.445, 0.698, 3.528, 1.4)^\mathsf{T}, -0.5\le x<0,
		\\(0.5, 0, 0.571, 1.4)^\mathsf{T},~~~~~~~~~~~~~0\le x\le 0.5.
	\end{cases}
\end{equation*}	
The final time is $T = 0.16$ and outflow boundary conditions are imposed on all boundaries. The computational results of density for the HWENO-U, HWENO-M and WENO-ZQ schemes are displayed in Fig. \ref{sec3:Lax1D_WENOHWENO}, which indicates that the results of the HWENO-U and HWENO-M schemes are more close to the exact solution than the WENO-ZQ scheme. The numerical solution of the WENO-ZQ scheme has obvious overshoots or undershoots in  Fig. \ref{sec3:Lax1D_WENOHWENO}, which is attributed to the nonlinear weights in \cite{ZhuJunQiuJX2017fvWENOZQ} for  violating the scale-invariant property. If the HWENO-U scheme also uses the original nonlinear weights \cite{ZhuJunQiuJX2017fvWENOZQ}, the  overshoots or undershoots also generate.\\ \indent
To prove that the proposed nonlinear weights \eqref{sec2:non_weight_1d} satisfy the scale-invariant property, similar to \cite{ChenWu}, we scale the initial conditions to be $(\zeta\rho,\mu,\zeta p,\gamma)$ with a constant $\zeta>0$. For this Riemann problem, the exact solution at time $T$ is $\zeta\rho(x,T)$. We compute this scaled case by the HWENO-U scheme with the scale-invariant nonlinear weights \eqref{sec2:non_weight_1d}  and original nonlinear weights \cite{ZhuJunQiuJX2017fvWENOZQ}, respectively. The computed  results are shown in Fig. \ref{sec3:Lax1D_HWENO_U_Scale}, which validates the effectiveness of scale-invariant nonlinear weights for the proposed HWENO-U scheme.
\end{example}
\begin{example}\label{Example:Blast1D}
		We solve the interaction of the blast wave problem for one-dimensional Euler equations with the initial conditions:	
		\begin{equation*}
			(\rho,\mu,p,\gamma)^\mathsf{T}=\begin{cases}
				(1,0,1000,1.4)^\mathsf{T},~~~0<x<0.1,
				\\(1,0,0.01,1.4)^\mathsf{T},~0.1<x<0.9,
				\\(1,0,100,1.4)^\mathsf{T},~~0.9<x<1.
			\end{cases}
		\end{equation*}	
		The computing time is $T = 0.038$ and reflective boundary conditions are imposed on all boundaries. The reference solution is generated by the classical WENO scheme \cite{js} using 2001 points. The density computed by HWENO-U, HWENO-M and WENO-ZQ schemes are plotted in Fig. \ref{sec3:Blast1D_WENOHWENO}, which shows the HWENO-U scheme has higher resolutions than the HWENO-M and WENO-ZQ schemes.
		\begin{figure}[!htpb]
			\centering
			\subfigure[Density with 800 cells]{\includegraphics[width=7.5cm,angle=0]{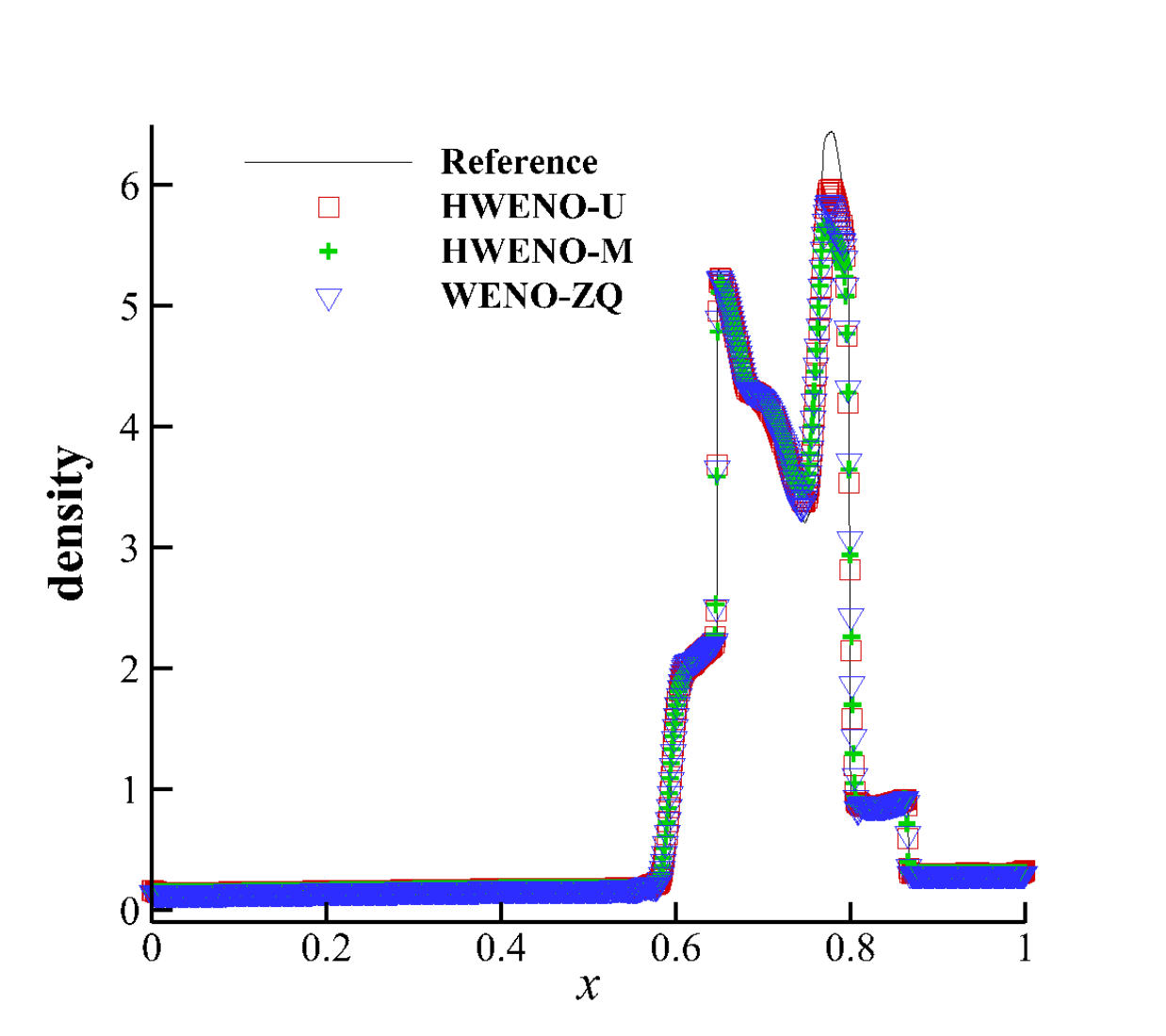}}
			\subfigure[A zoomed-in figure]{\includegraphics[width=7.5cm,angle=0]{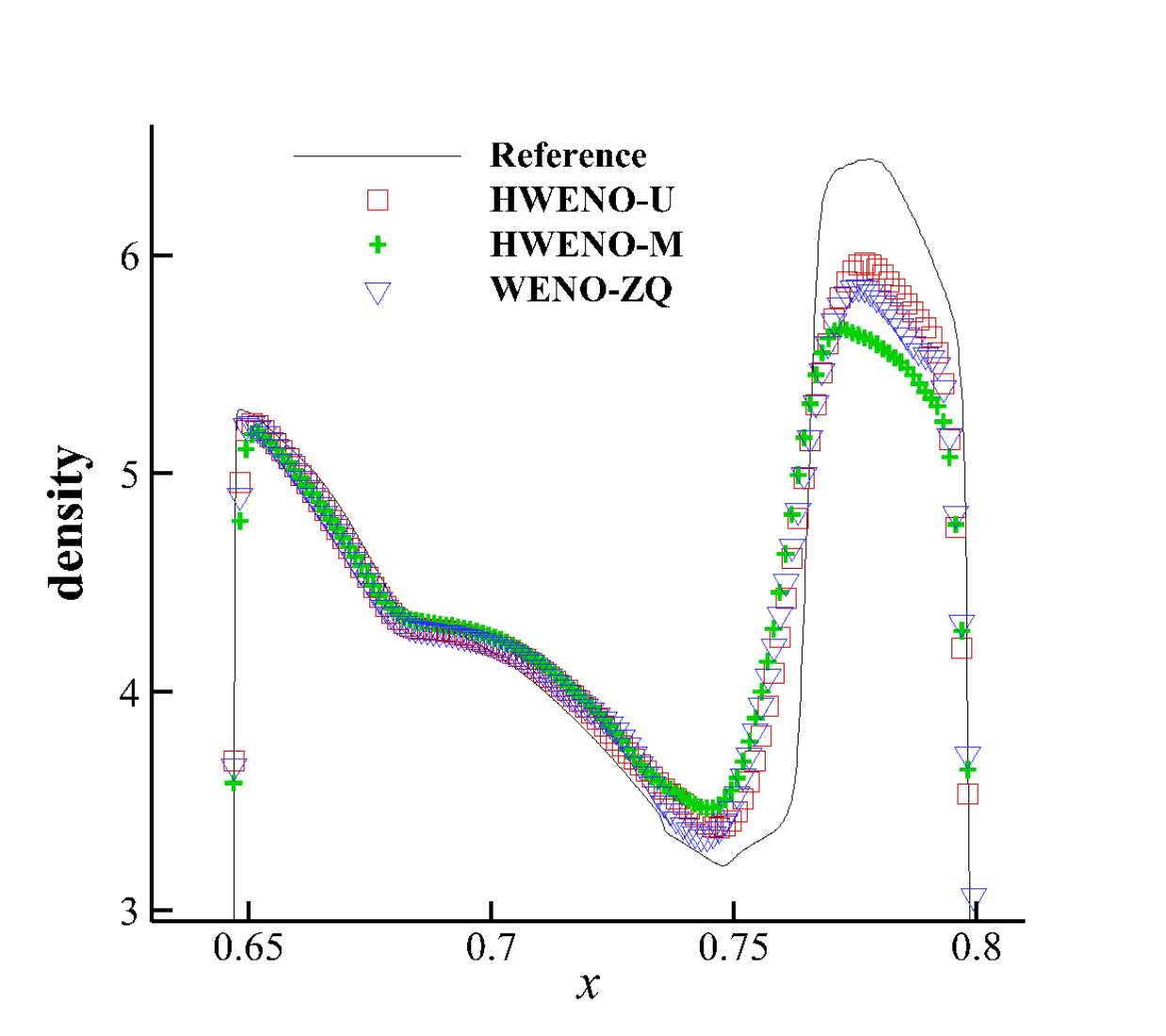}}
			\caption{Example \ref{Example:Blast1D}.
				The results of solution computed by the HWENO-U, HWENO-M and WENO-ZQ schemes.}
			\label{sec3:Blast1D_WENOHWENO}
		\end{figure}
\end{example}
\begin{example}\label{Example:Shu1D}
		We solve the Shu-Osher problem for one-dimensional Euler equations, which describes the interaction between shock and entropy waves. The initial condition is
		\begin{equation*}
			(\rho,\mu,p,\gamma)=\begin{cases}
				(3.857143,2.629369,10.333333,1.4),~-5\le x<-4,
				\\(1+0.2\sin(5x),0,1,1.4),~~~~~~~~~~~~~~~-4\le x\le 5.
			\end{cases}
		\end{equation*}		
		The final time is $T = 1.8$ and outflow boundary conditions are imposed on all boundaries. The density  calculated by the HWENO-U, HWENO-M, and WENO-ZQ schemes are displayed in Fig. \ref{sec3:Shu1D}, indicating that the HWENO-U and WENO-ZQ schemes exhibit similar results but both have better resolutions than the HWENO-M scheme.
		\begin{figure}[t]
			\centering
			\subfigure[Density with 400 cells]{\includegraphics[width=8cm,angle=0]{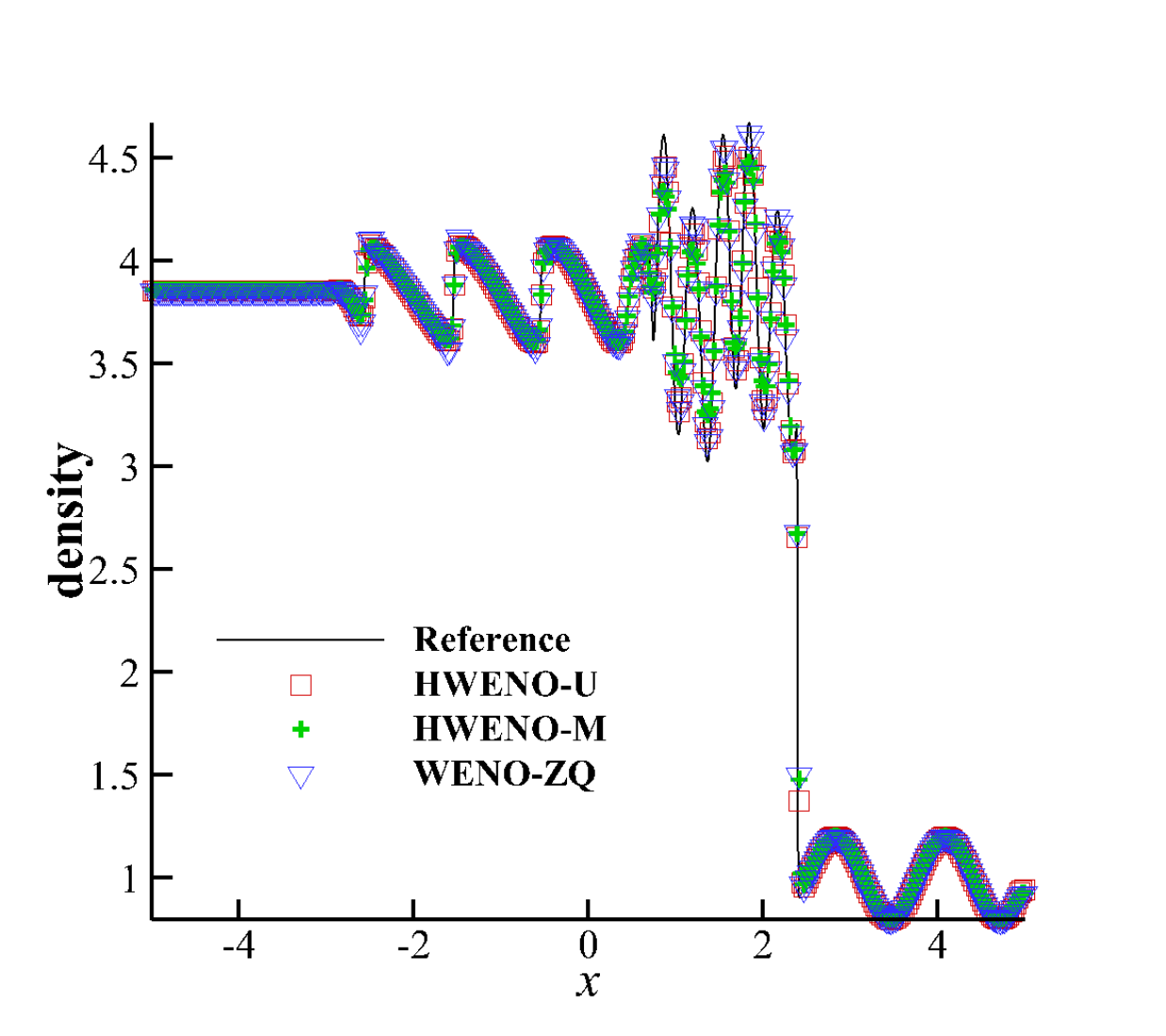}}
			\subfigure[A zoomed-in figure]{\includegraphics[width=8cm,angle=0]{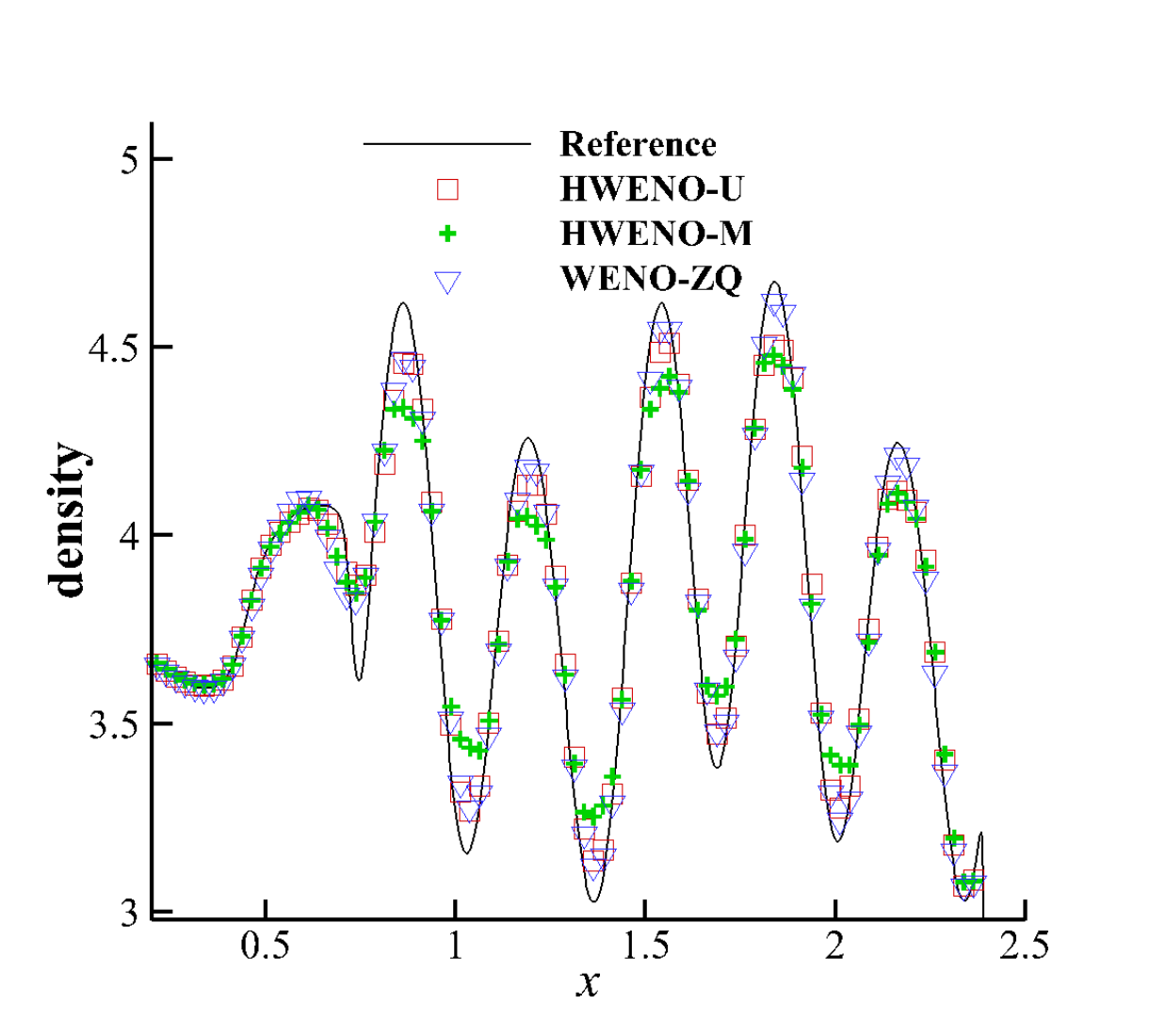}}
			\caption{Example \ref{Example:Shu1D}. The results of solution computed by the HWENO-U, HWENO-M and WENO-ZQ schemes with 400 cells.}
			\label{sec3:Shu1D}
		\end{figure}
\end{example}

\begin{example}\label{Example:DoubleRare}
		We solve the double rarefaction wave problem \cite{linde1997robust} for  one-dimensional Euler equations with the initial condition
		\begin{equation*}
			(\rho,\mu,p,\gamma)=\begin{cases}
				(7,-1,0.2,1.4),~-1<x<0,
				\\(7,1,0.2,1.4),~~~~~~~0<x<1.
			\end{cases}
		\end{equation*}	
		The final time is $T = 0.6$ and outflow boundary conditions are imposed on all boundaries. The results computed by the HWENO-U, HWENO-M and WENO-ZQ schemes are shown in Fig. \ref{sec3:DoubleRare}. Numerically we find that such three schemes work well for this extreme problem without PP limiters, but the two HWENO schemes have more compact reconstructed stencils.
		\begin{figure}[!htpb]
			\centering
			{\includegraphics[width=5.35cm,angle=0]{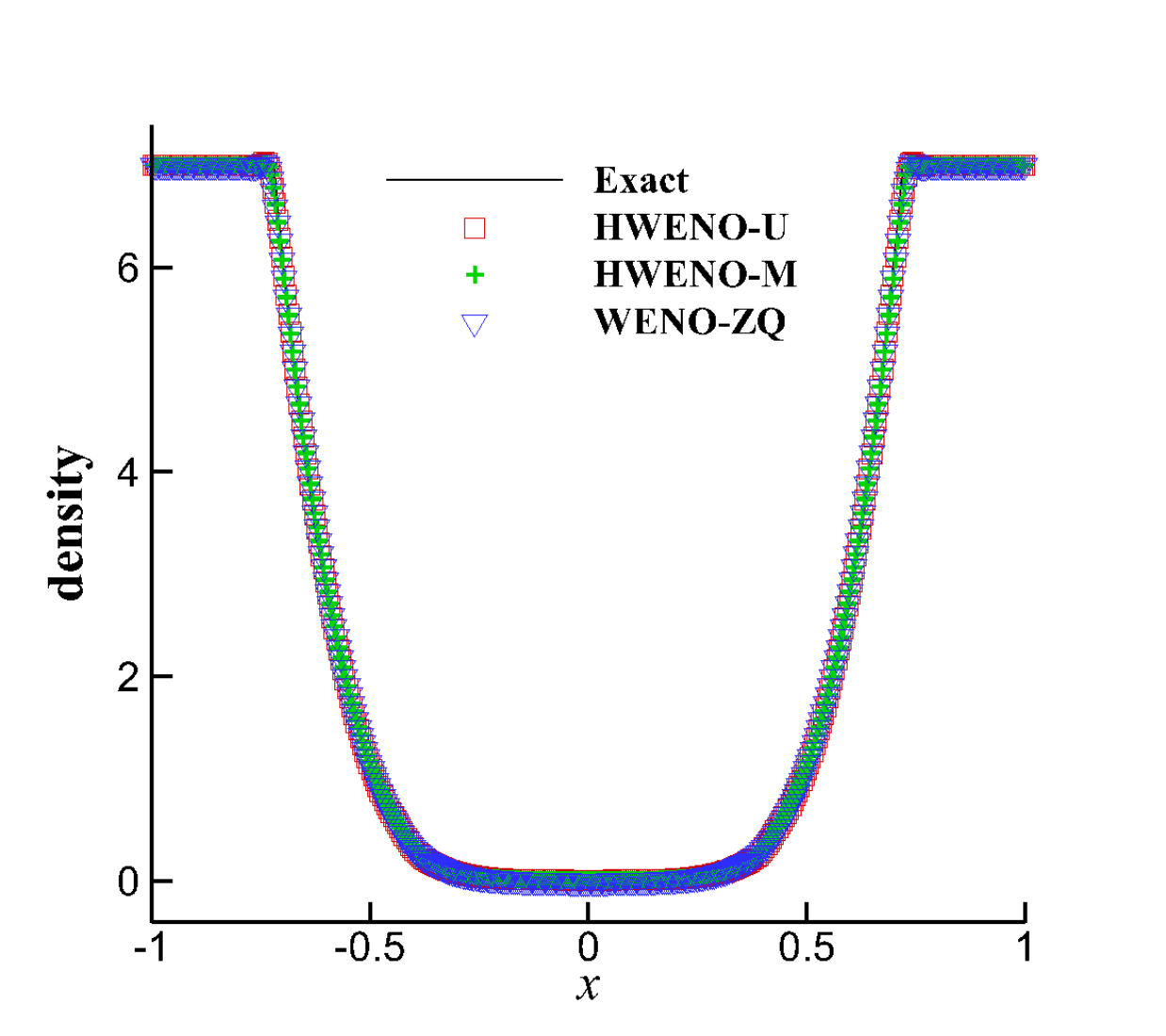}}
			{\includegraphics[width=5.35cm,angle=0]{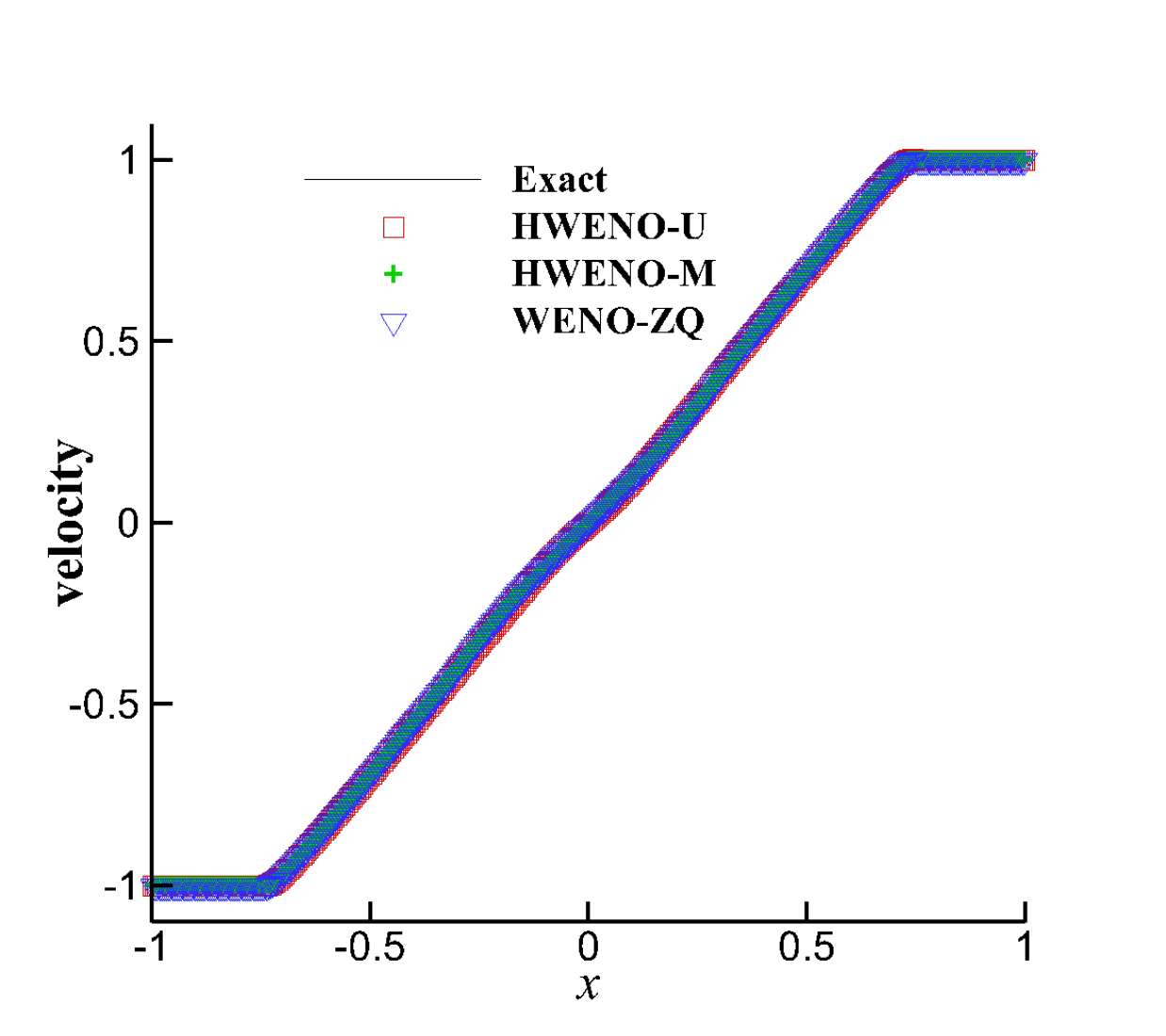}}
			{\includegraphics[width=5.35cm,angle=0]{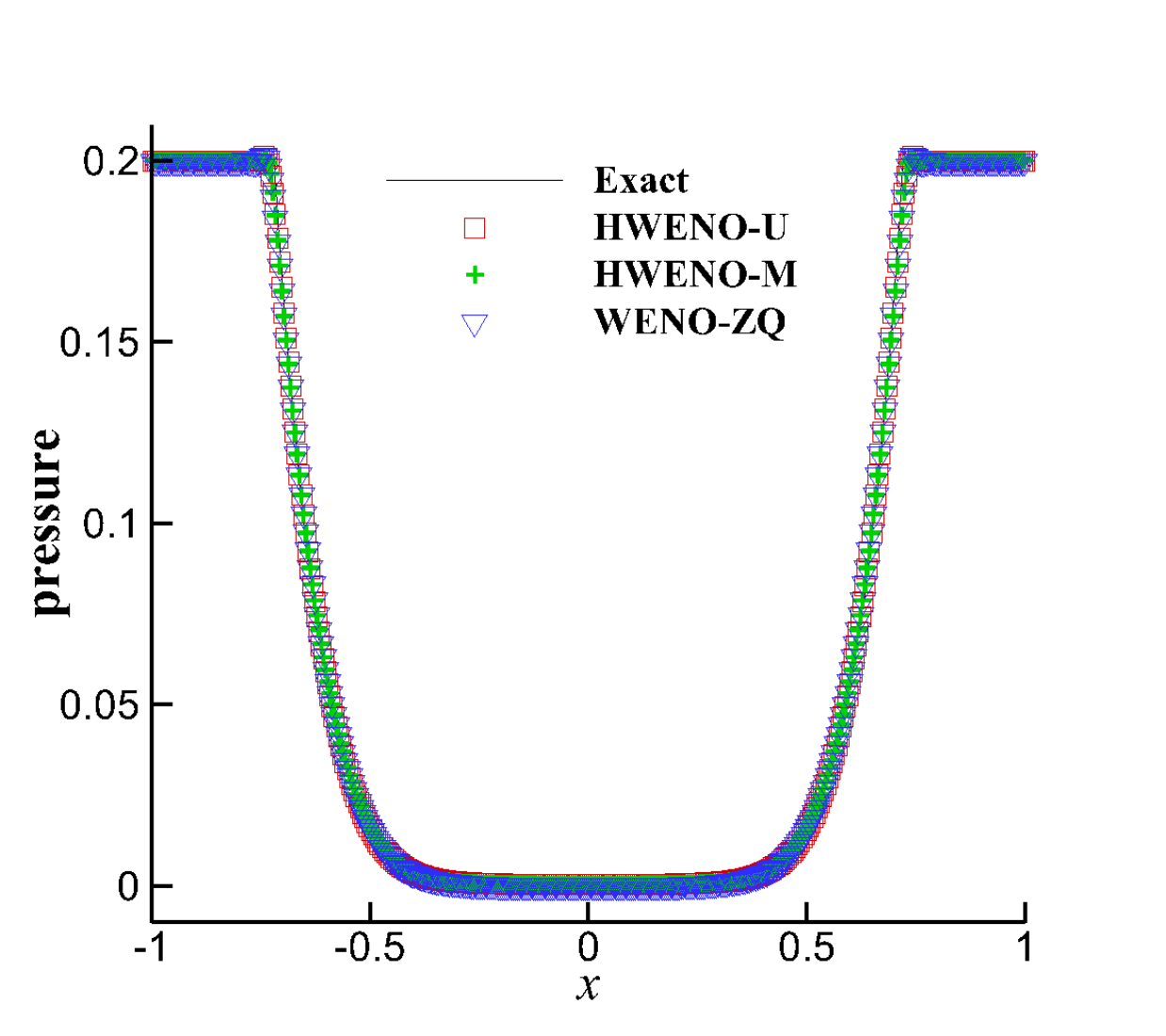}}
			\caption{Example \ref{Example:DoubleRare}. Double rarefaction wave problem with $400$ cells.}
			\label{sec3:DoubleRare}
		\end{figure}
\end{example}
\begin{example}\label{Example:Sedov1D}
		We solve the Sedov blast wave problem for one-dimensional Euler equations with the initial condition
		\begin{equation*}
			(\rho,\mu,E,\gamma)=\begin{cases}
				(1,0,10^{-12},1.4),~~x\in [-2,2] \setminus  \mbox{the center cell},
				\\(1,0,\frac{3200000}{\dx},1.4),~x\in \mbox{the center cell}.
			\end{cases}
		\end{equation*}	
		The final time is $T = 0.001$ and outflow boundary conditions are imposed on all boundaries. The exact solution is provided in \cite{sedov1959similarity,korobeinikov1991problems}. We present the computational density in Fig. \ref{sec3:Sedov1D} for the HWENO-U, HWENO-M and WENO-ZQ schemes. Also, we do not use any PP limiters for the three scheme in this extreme problem, and the results are non-oscillatory with high resolutions.
		\begin{figure}[t]
			\centering
			{\includegraphics[width=5.35cm,angle=0]{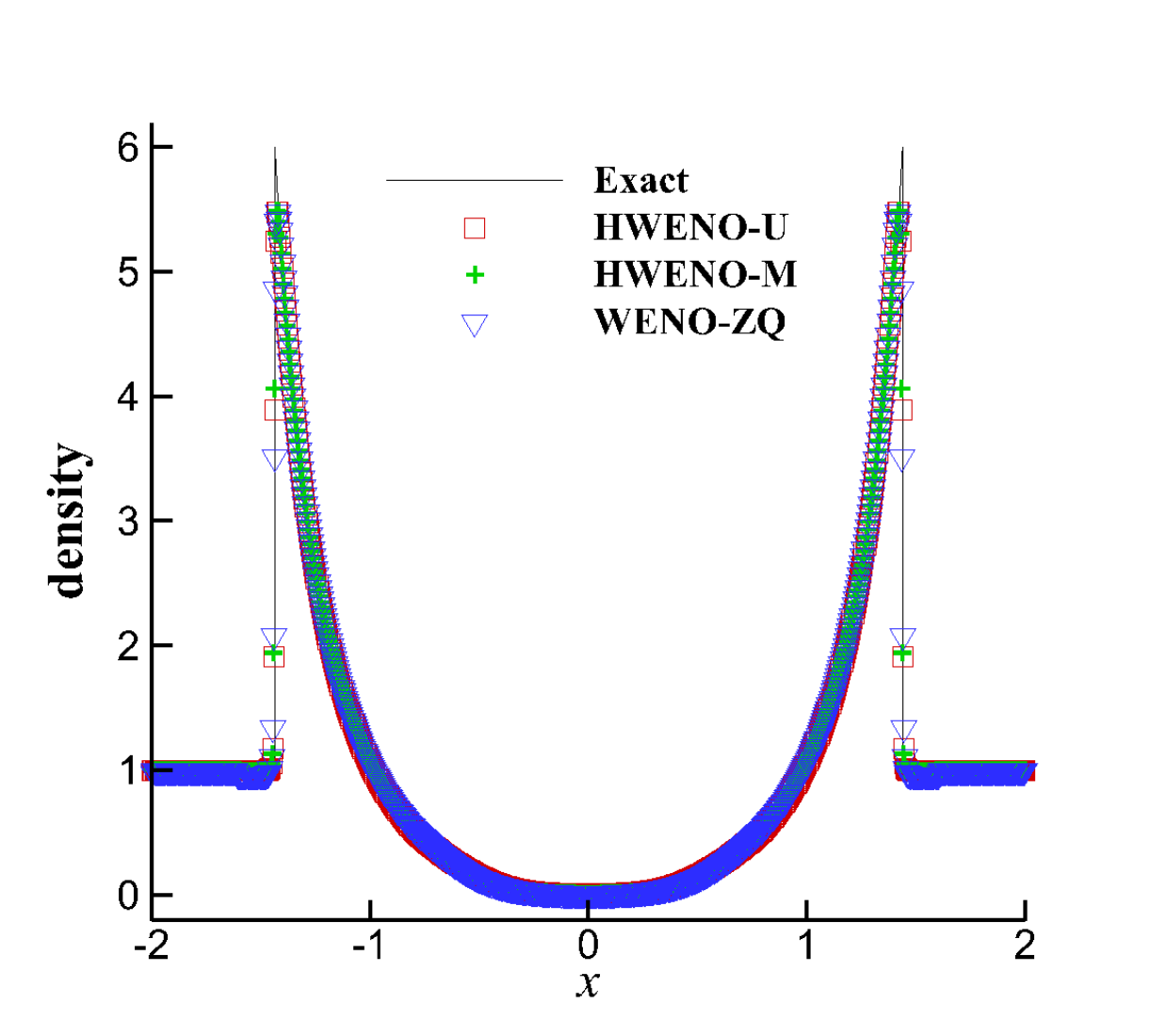}}
			{\includegraphics[width=5.35cm,angle=0]{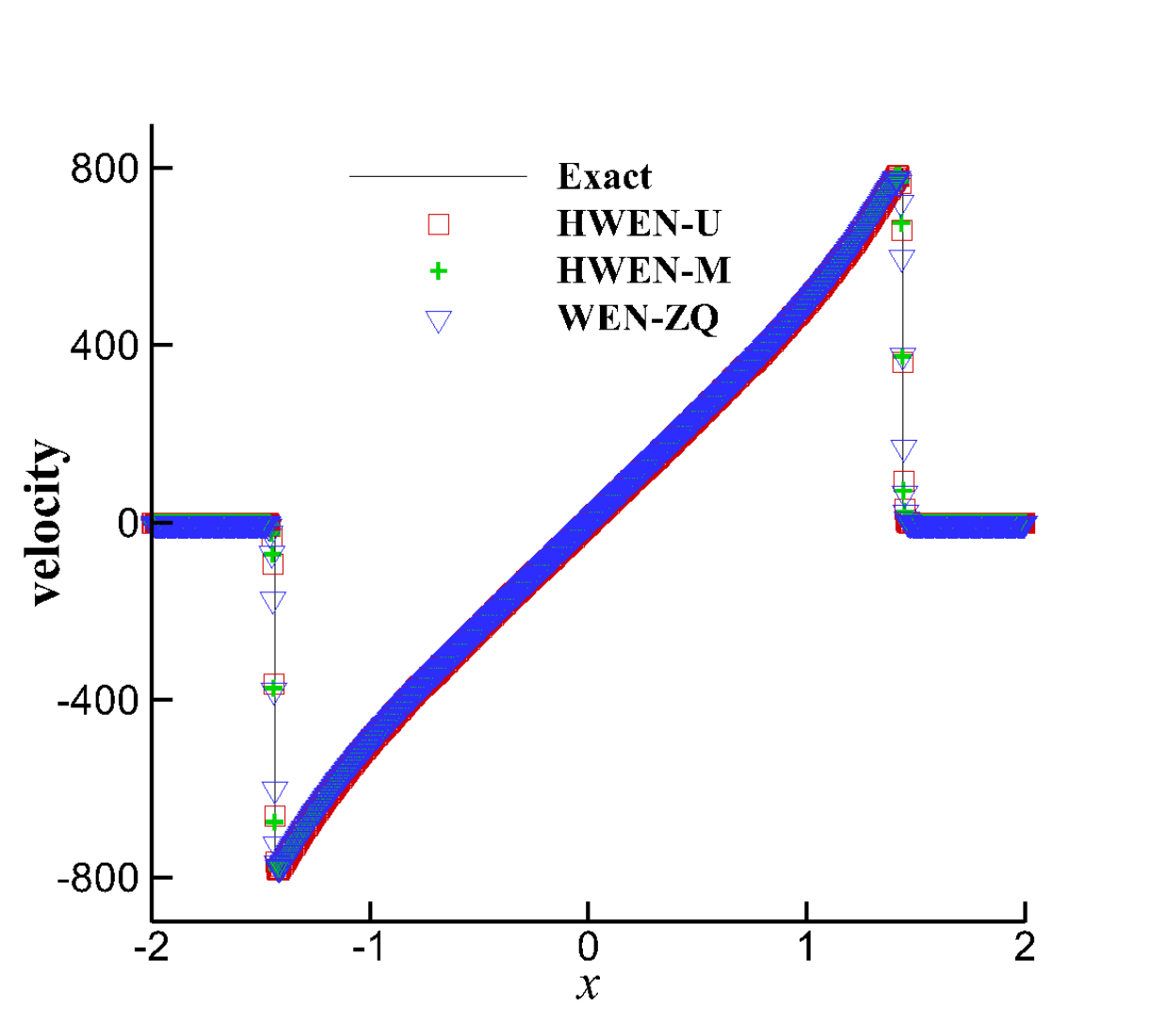}}
			{\includegraphics[width=5.35cm,angle=0]{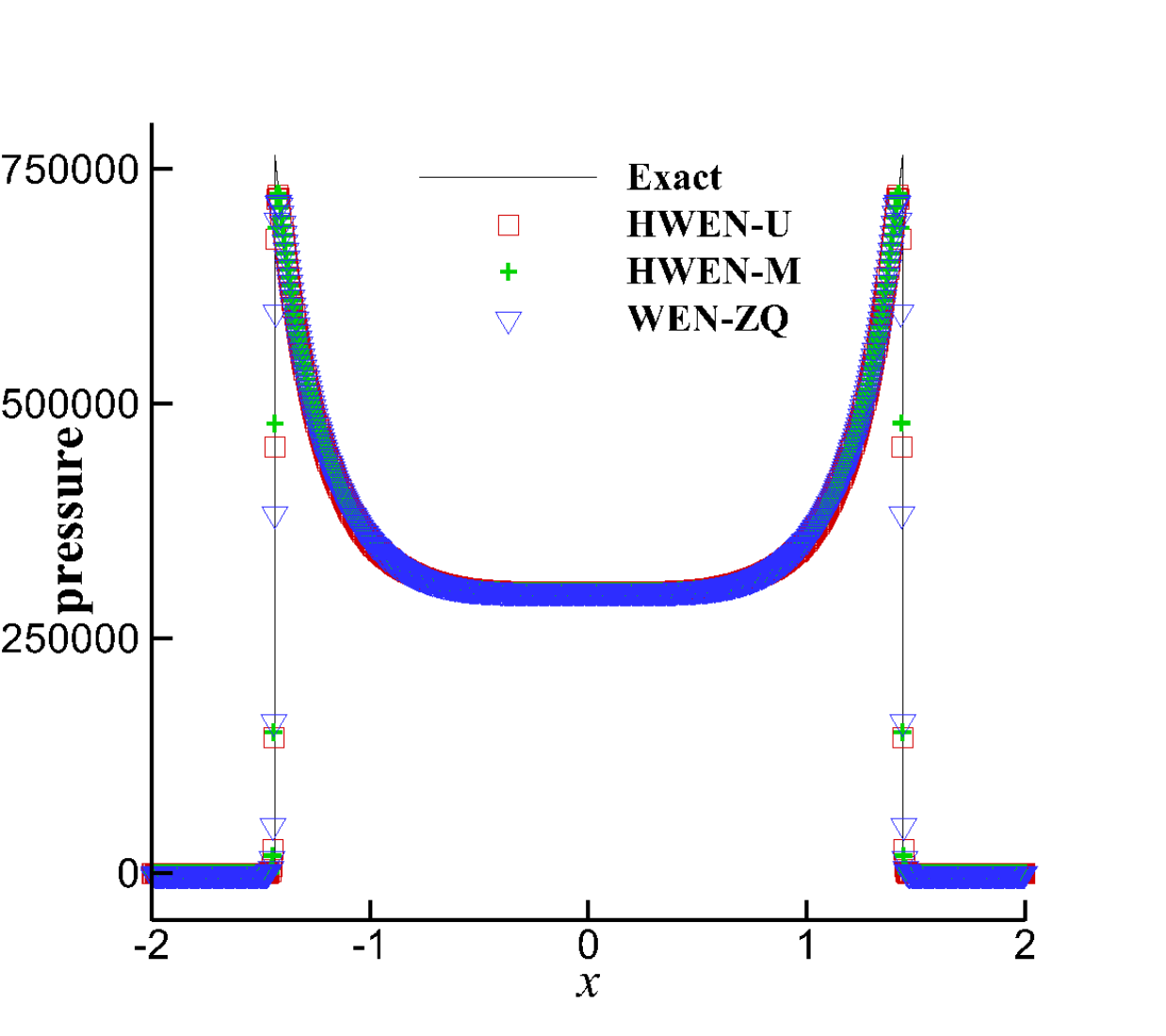}}
			\caption{Example \ref{Example:Sedov1D}. One-dimensional Sedov problem with $800$ cells.}
			\label{sec3:Sedov1D}
		\end{figure}
\end{example}
\begin{example}\label{Example:DoubleMach}
		We solve the double Mach reflection problem \cite{Paul1984} for two-dimensional Euler equations. The computational domain is $[0,4]\times[0,1]$ and the initial condition is
		\begin{equation*}
			(\rho,\mu,\nu,p,\gamma)=
			\begin{cases}
				(8,\frac{33}{4}\sin(\frac{\pi}{3}),-\frac{33}{4}\cos(\frac{\pi}{3}),116.5,1.4),~x<\frac{1}{6}+\frac{y}{\sqrt{3}},
				\\(1.4,0,0,1,1.4),\quad\quad\quad\quad\quad\quad\quad\quad\quad\mbox{otherwise}.
			\end{cases}
		\end{equation*}
		The boundary conditions are set as inflow on the left, outflow on the right and bottom. The reflection boundary condition are applied for the bottom boundary  starting from $x=\frac16$ to $x=4$, while the rest part from $x=0$ to $x=\frac16$ imposes the exact post-shock condition. Besides, the upper boundary is the exact motion of a Mach 10 shock.
The final time is $T = 0.2$. The computational results of density for the HWENO-U and HWENO-M schemes are showed in Fig. \ref{sec3:DoubleMach}. We can see that the two results are similar, but the HWENO-U scheme has simpler procedures with unified stencils.
		\begin{figure}[t]
			\centering
			\subfigure[HWENO-U]{\includegraphics[width=16.5cm,height=5.5cm,angle=0]{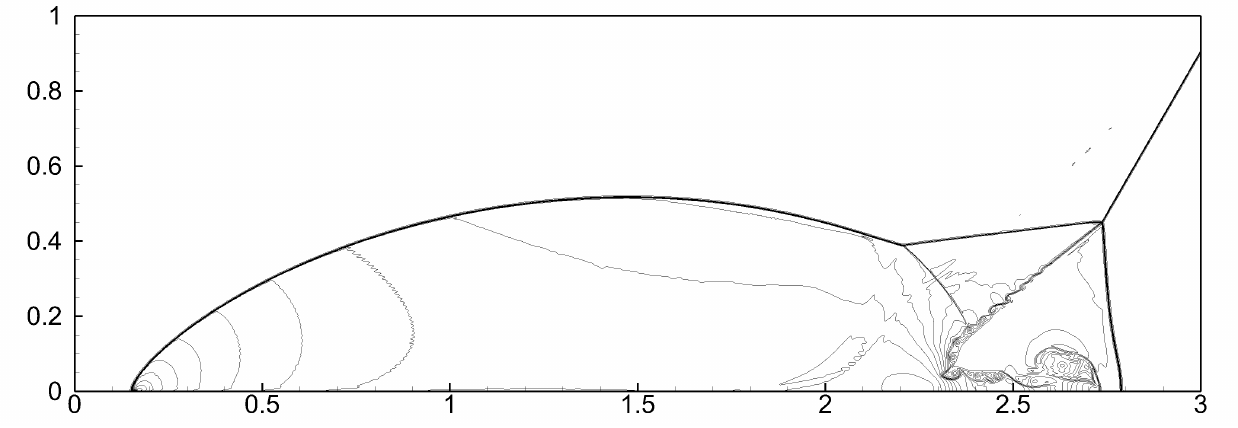}}
			\subfigure[HWENO-M]{\includegraphics[width=16.5cm,height=5.5cm,angle=0]{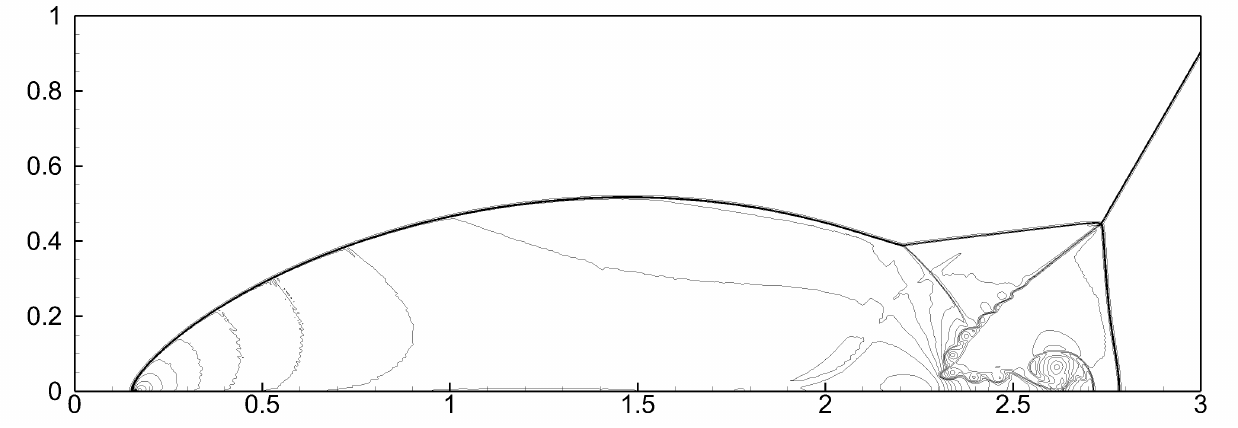}}
			\caption{Example \ref{Example:DoubleMach}.
				Double Mach reflection problem. Contour plots of density with 30 equally spaced lines from  $1.5$ to $22.7$. Uniform meshes: $1920\times480$.}
			\label{sec3:DoubleMach}
		\end{figure}
\end{example}
\begin{example}\label{Example:Step}
		We solve the forward step problem \cite{Paul1984} for two-dimensional Euler equations, which contains a Mach 3 wind tunnel with a step. The computational domain is $[0,0.6]\times[0,1] \cup [0.6,1]\times[0.2,1]$ and the initial condition is a right-going Mach 3 flow. Reflective boundary conditions are applied along the walls of the tunnel, and inflow and outflow boundary conditions are implemented at the entrance and exit respectively. The final time is $T = 4$. The density results computed by the HWENO-U and HWENO-M schemes are shown in Fig. \ref{sec3:Step}. We can observe that both results are comparable for the HWENO-U and HWENO-M schemes.
		\begin{figure}[t]
			\centering
			\subfigure[HWENO-U]{\includegraphics[width=16.5cm,height=5.5cm,angle=0]{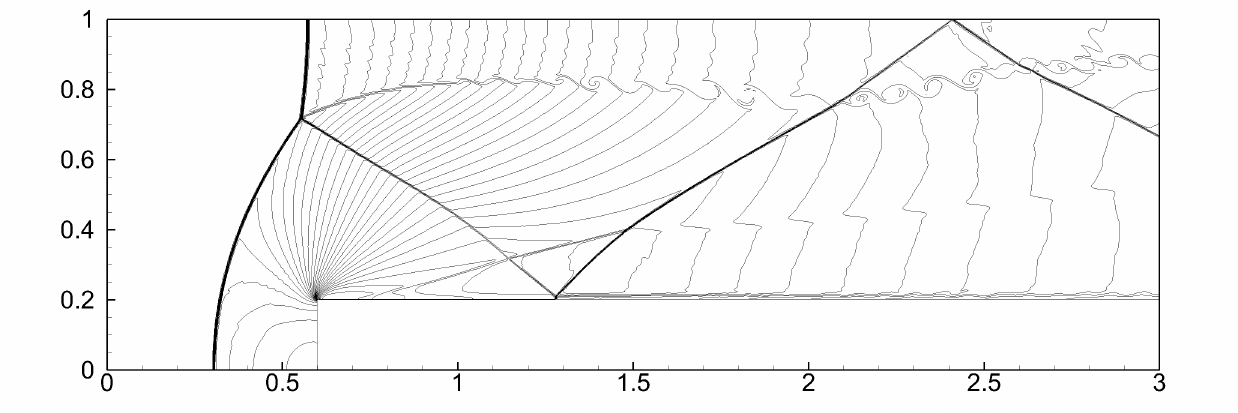}}
			\subfigure[HWENO-M]{\includegraphics[width=16.5cm,height=5.5cm,angle=0]{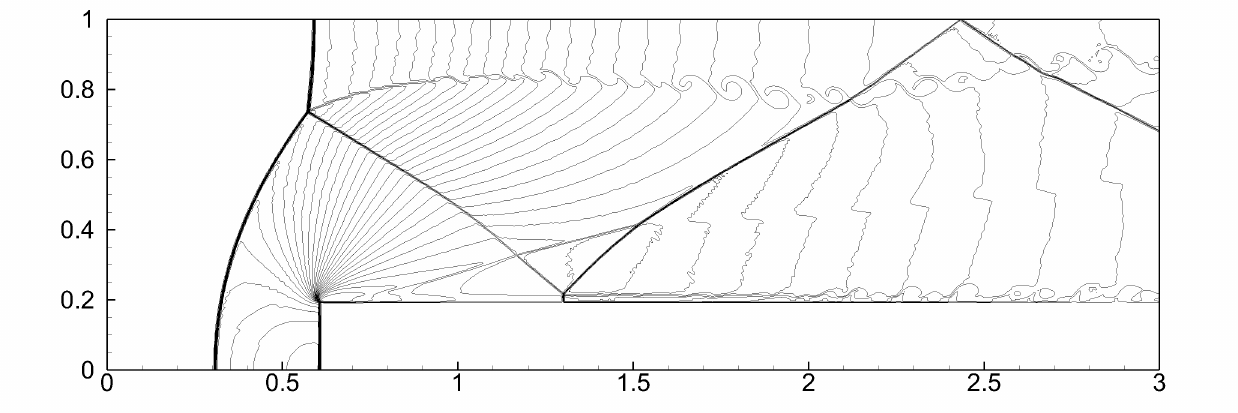}}
			\caption{Example \ref{Example:Step}.
				Step forward problem. Contour plots of density with 30 equally spaced lines from  $0.32$ to $6.15$. Uniform meshes: $960\times320$.}
			\label{sec3:Step}
		\end{figure}
\end{example}
\begin{example}\label{Example:Sedov2D}
		We solve a Sedov blast wave problem \cite{sedov1959similarity,korobeinikov1991problems}  for two-dimensional Euler equations.	 The computational domain is $[0,1.1]\times[0,1.1]$ and the initial condition is
		\begin{equation*}
			(\rho,\mu,\nu,E,\gamma)=\begin{cases}
				(1,0,0,\frac{0.244816}{\dx\dy},1.4),~(x,y)\in[0,\dx]\times[0,\dy],
				\\(1,0,0,10^{-12},1.4),~~~~~\mbox{otherwise}.
			\end{cases}
		\end{equation*}	
		Reflective boundary conditions are employed on the left and bottom, while outflow conditions are applied on the right and upper boundaries. The computational results at the final time $T=1$ are presented in Fig. \ref{sec3:Sedov2D} for both the HWENO-U and HWENO-M schemes with PP limiters. Notably, it is essential to utilize PP limiters in this case, as both schemes would fail to work effectively without them due to negative  densities or pressures. This extreme problem involves very strong shock and the variation of density is pretty large. From Fig. \ref{sec3:Sedov2D}, we can observe obviously that there exist numerical oscillations even using PP limiters for the HWENO-M scheme with original nonlinear weights \cite{ZhuJunQiuJX2017fvWENOZQ}, since the PP limiters can keep the positivity of density and pressure but cannot control numerical oscillations.  On the contrary, the HWENO-M scheme with the scaling-invariant weights in Eq. \eqref{sec2:non_weight_2d} and the proposed HWENO-U scheme behave similar and comparable results as in the reference \cite{ZhangS}. Also, the HWENO-U scheme  has higher resolutions and better performances than the HWENO-M scheme with the nonlinear weights in \cite{ZhuJunQiuJX2017fvWENOZQ} or Eq. \eqref{sec2:non_weight_2d}. The results presented above demonstrate the necessity  of  scale-invariant weights, especially for the problems with sharp scale variations.
		\begin{figure}[t]
			\centering
			\subfigure[HWENO-U with original nonlinear weights \cite{ZhuJunQiuJX2017fvWENOZQ}) ]
			{\includegraphics[width=5.25cm,height=5.25cm,angle=0]{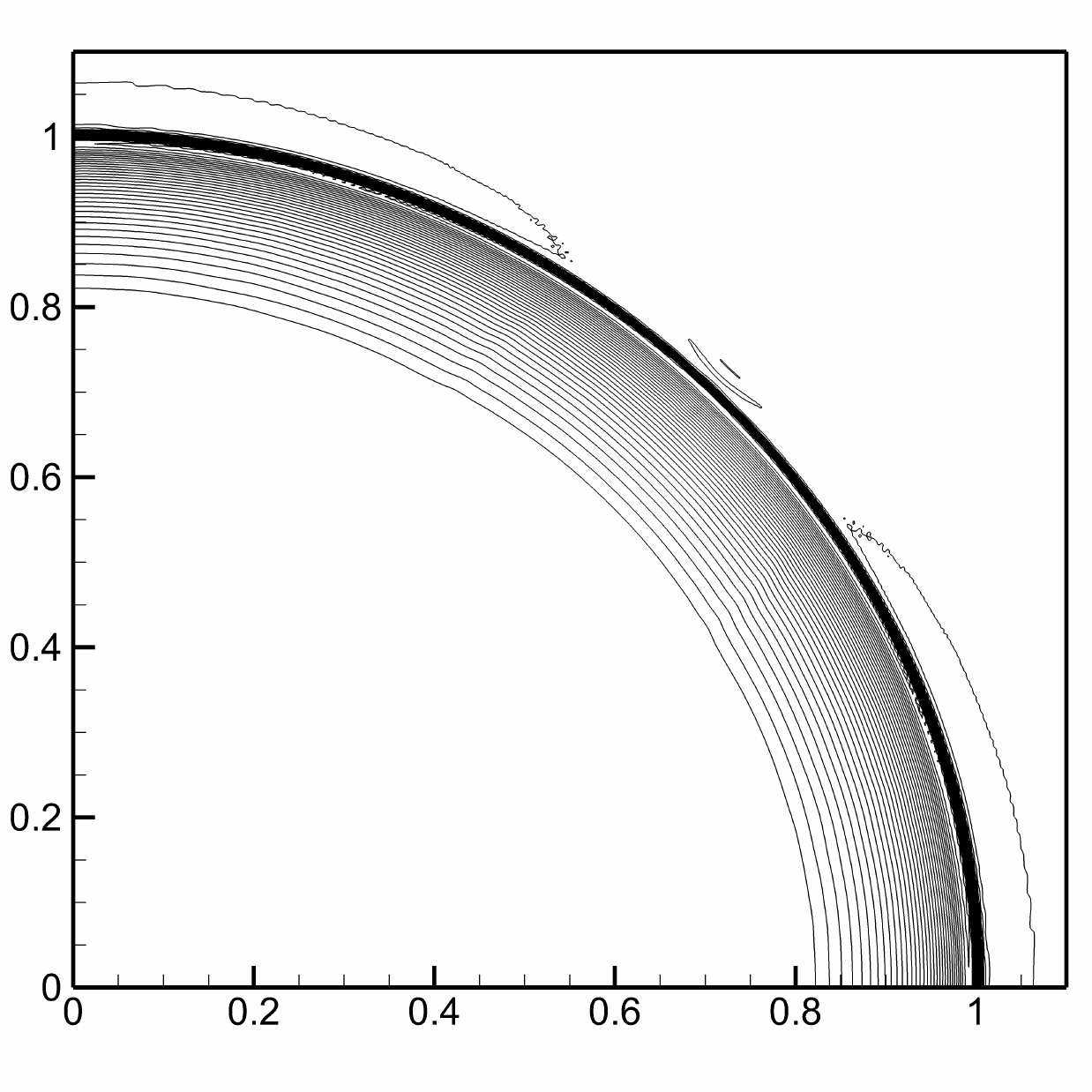}}
			\subfigure[HWENO-M with original nonlinear weights \cite{ZhuJunQiuJX2017fvWENOZQ} ]
			{\includegraphics[width=5.25cm,height=5.25cm,angle=0]{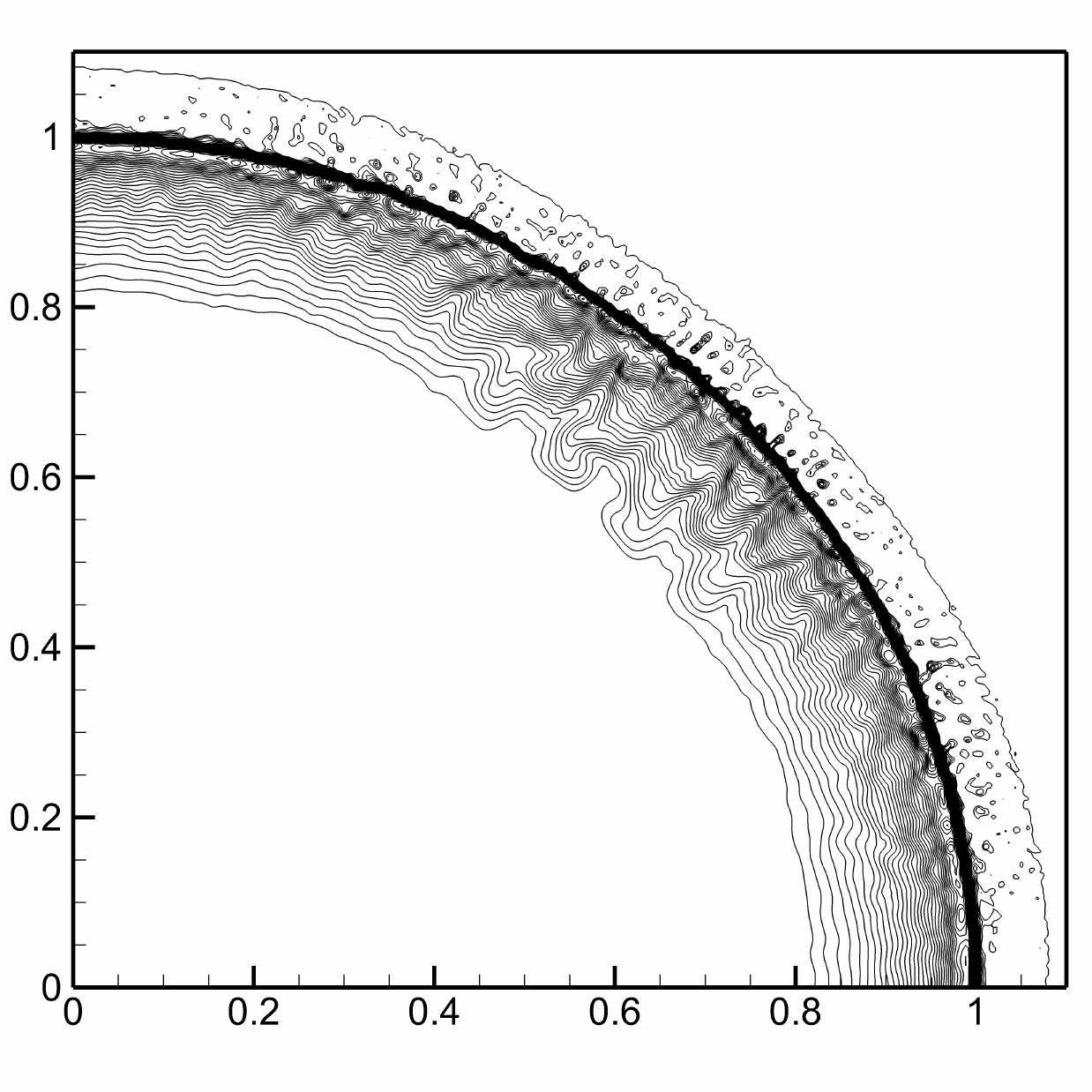}}
			\subfigure[Density at $x=y$]{\includegraphics[width=5.5cm,height=5.5cm,angle=0]{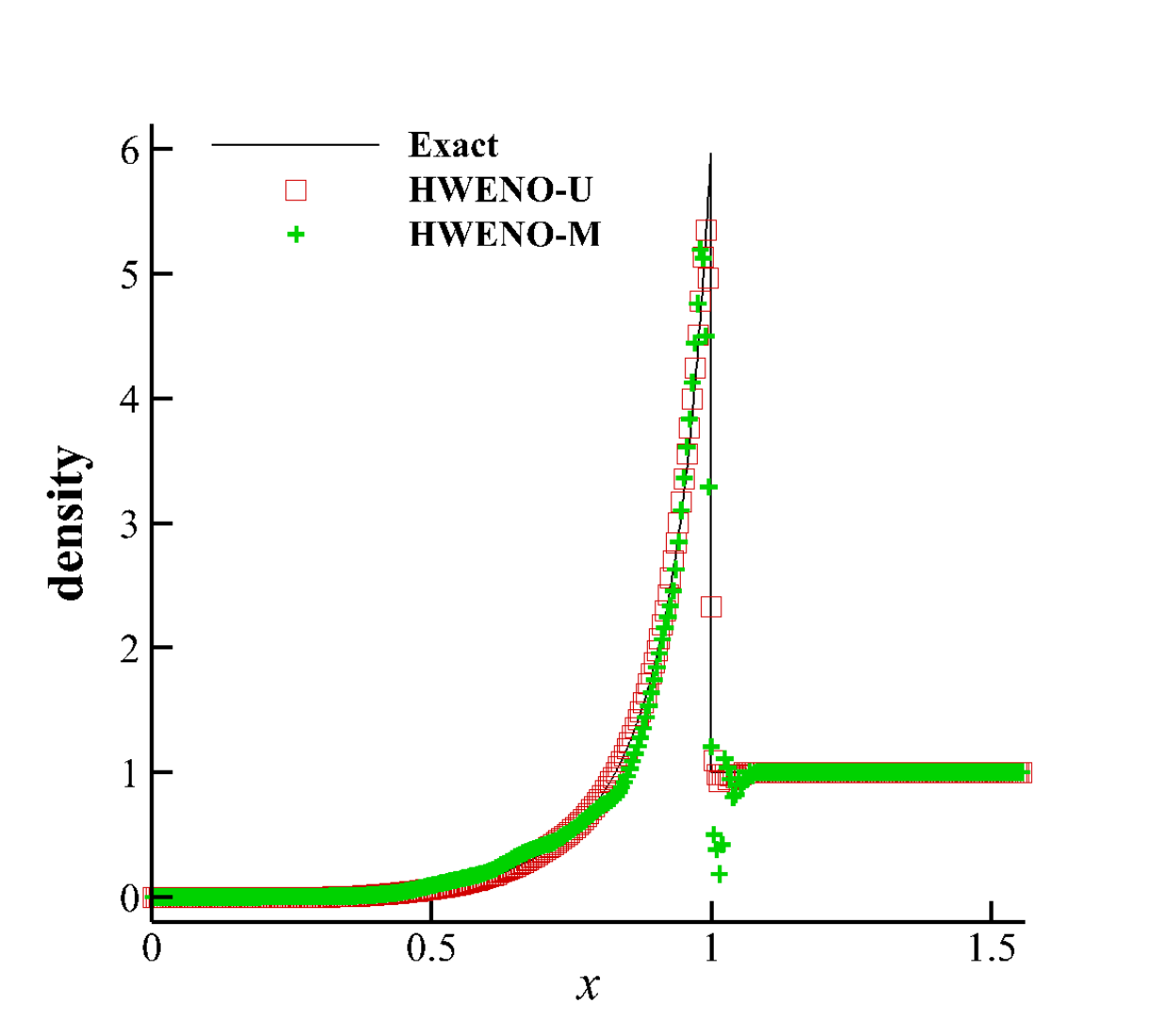}}
			\subfigure[HWENO-U with {the proposed} scale-invariant weights \eqref{sec2:non_weight_2d}] {\includegraphics[width=5.25cm,height=5.25cm,angle=0]{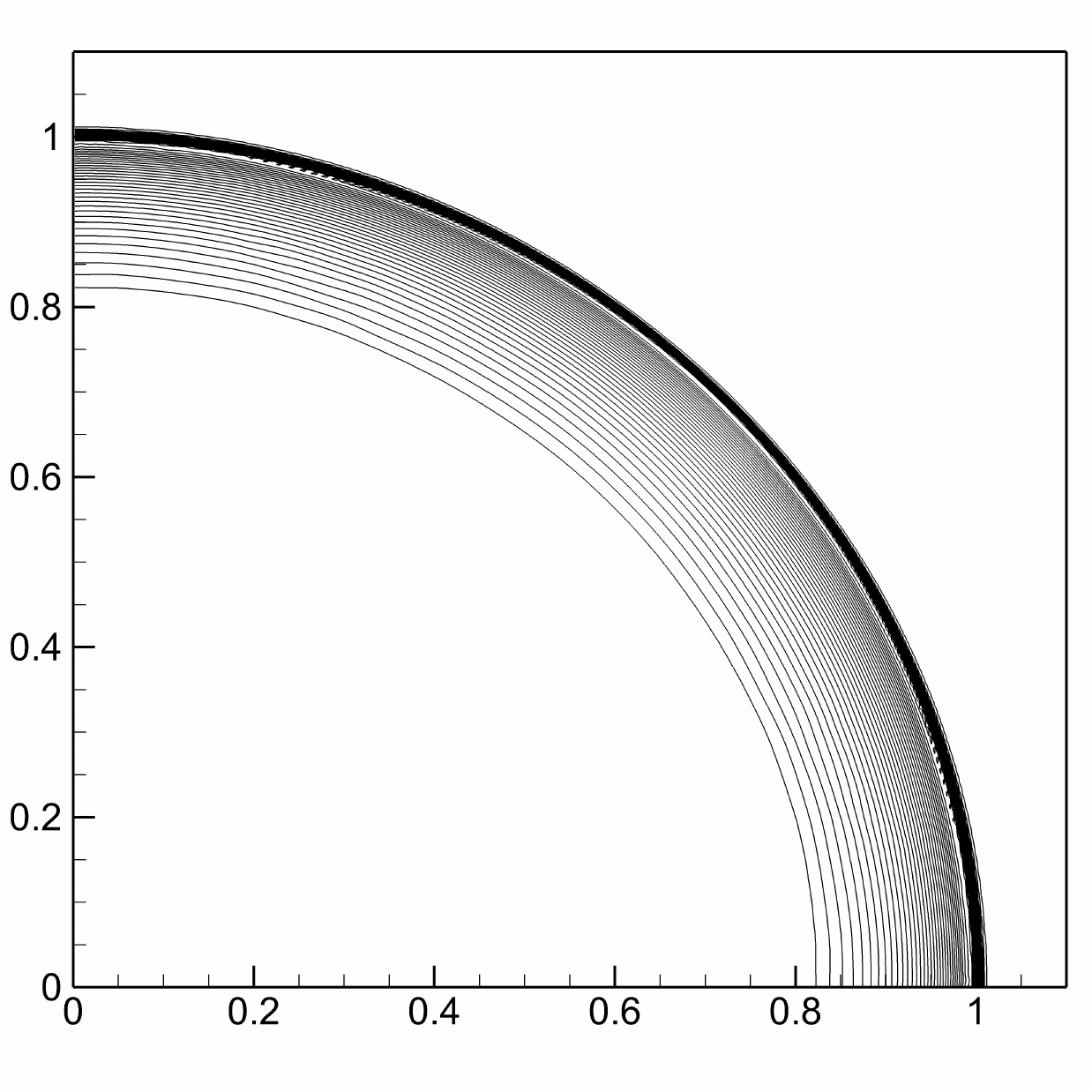}}
			\subfigure[HWENO-M with {the proposed} scale-invariant weights \eqref{sec2:non_weight_2d}] {\includegraphics[width=5.25cm,height=5.25cm,angle=0]{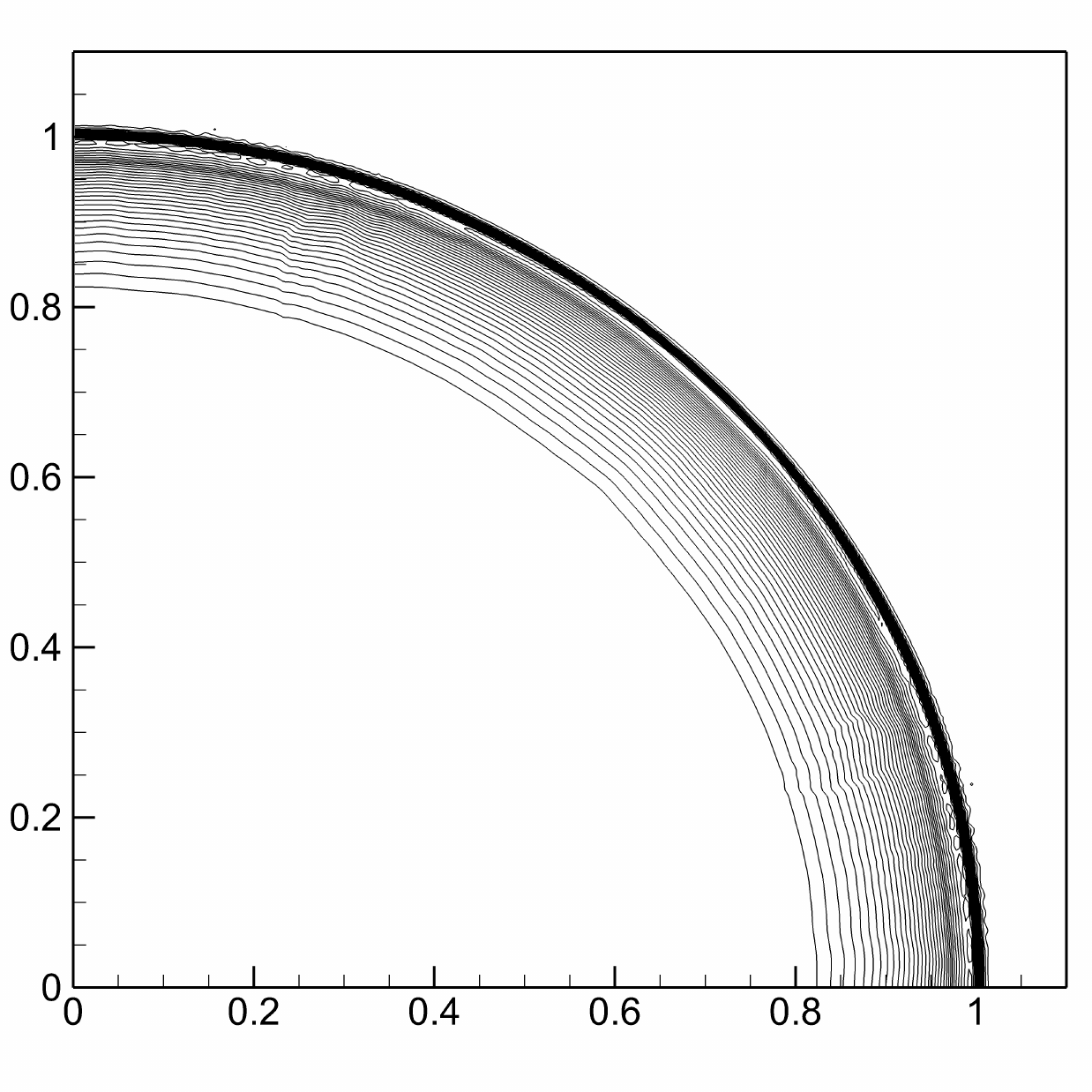}}
			\subfigure[Density at $x=y$]{\includegraphics[width=5.5cm,height=5.5cm,angle=0]{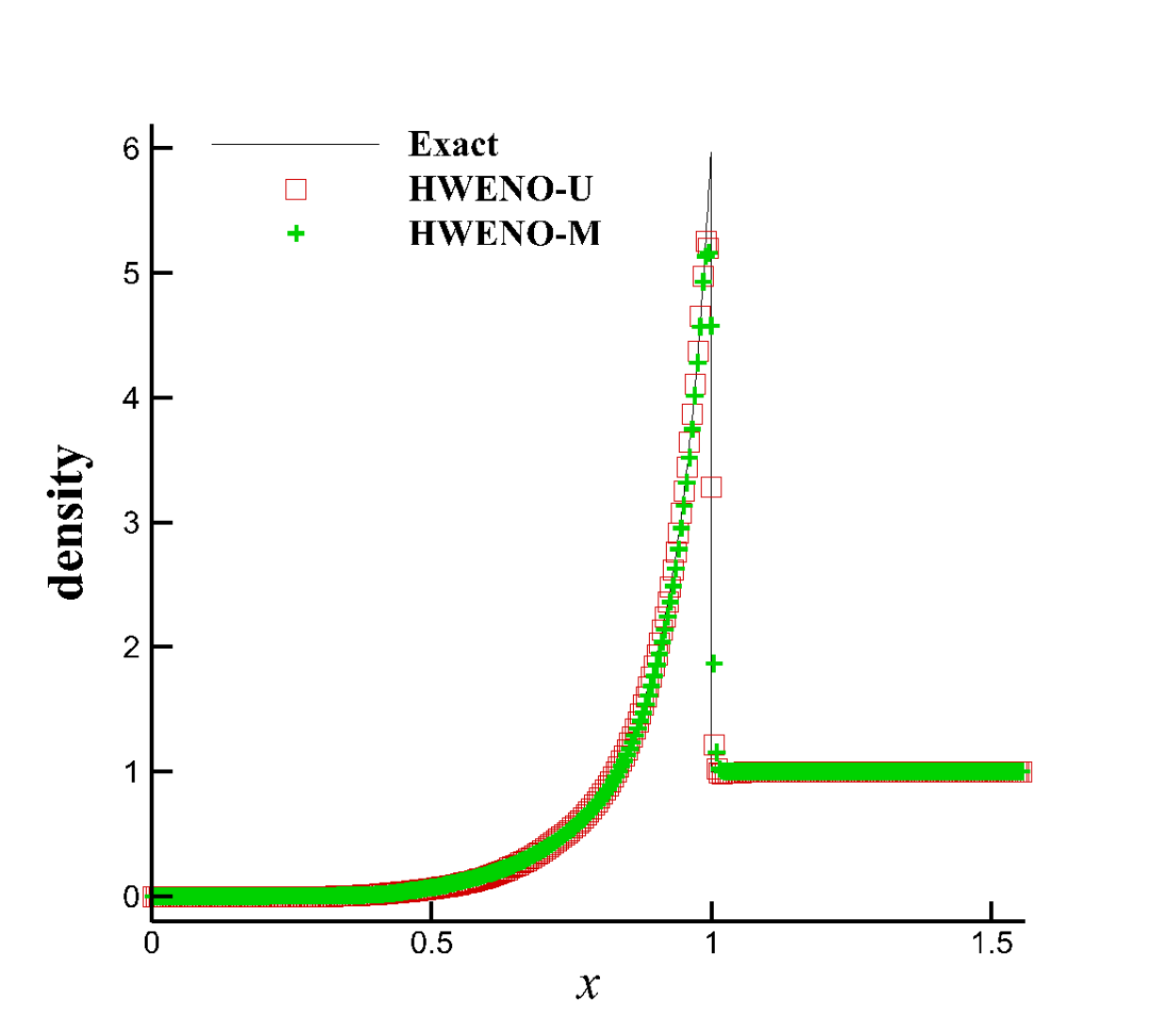}}
			\caption{Example \ref{Example:Sedov2D}.
				Two-dimensional Sedov problem.  Contour plots of density with 40 equally spaced lines from 0.95 to 6.	Uniform meshes: $320\times320$.}
			\label{sec3:Sedov2D}
		\end{figure}
\end{example}
\begin{example}\label{Example:HM2000}
	\begin{figure}[t]
		\centering
		\subfigure[HWENO-U] {\includegraphics[width=8cm,height=4cm,angle=0]{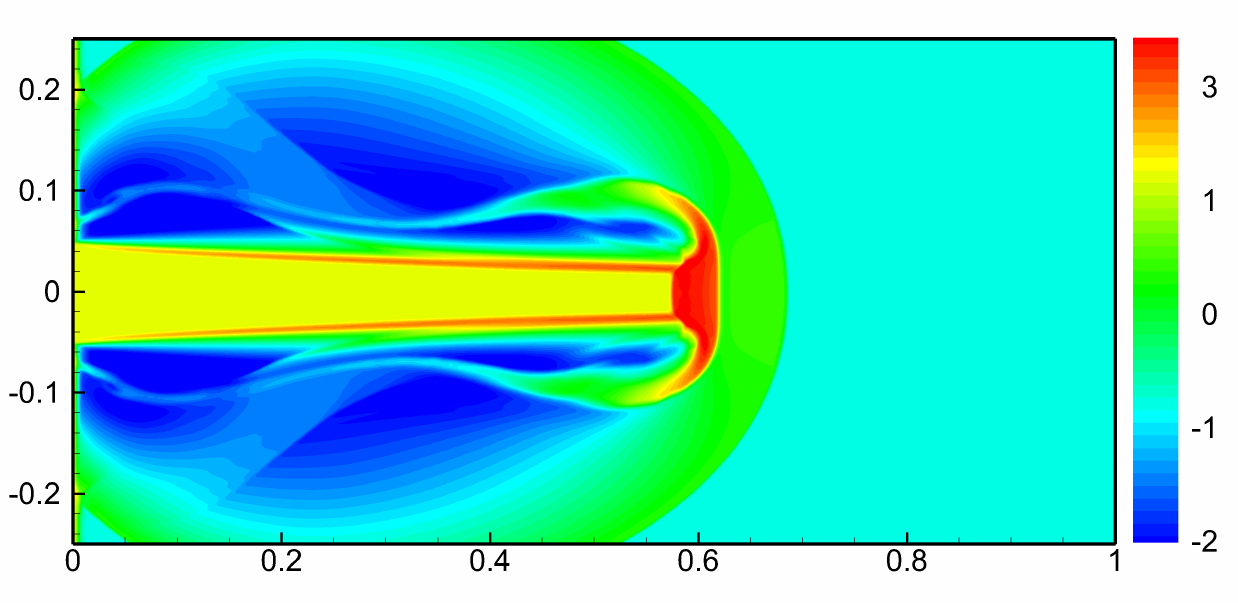}}
		\subfigure[HWENO-M] {\includegraphics[width=8cm,height=4cm,angle=0]{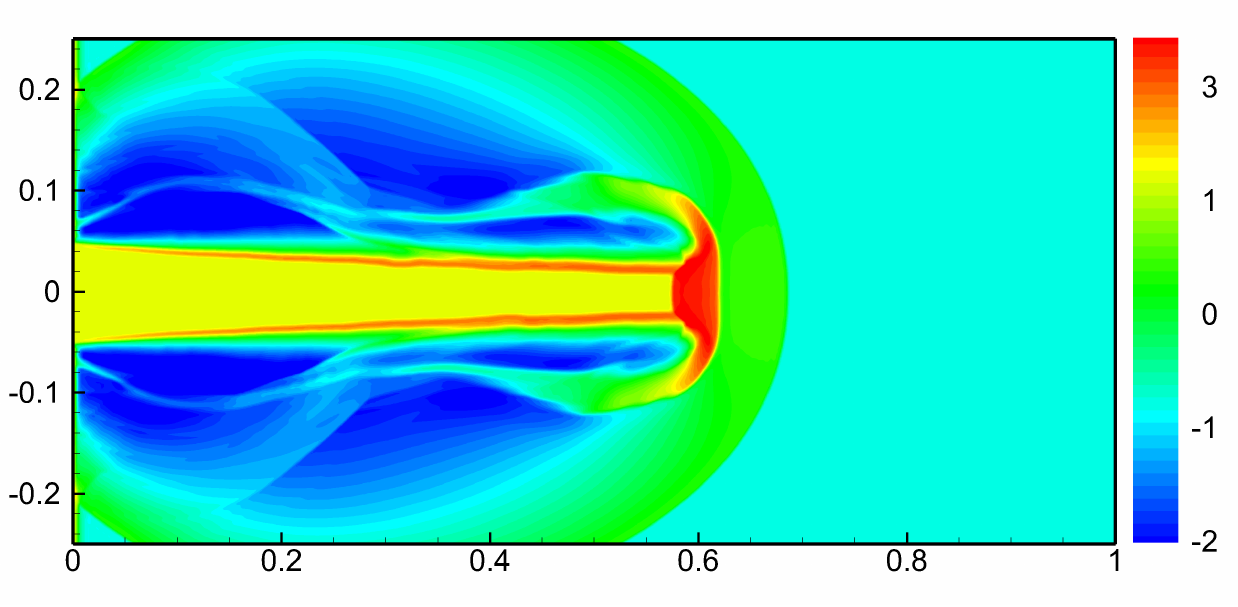}}
		\caption{Example \ref{Example:HM2000}. High Mach 2000 problem. Contour plots of density with 40 equally spaced lines from -2 to 3 and scales are logarithmic. Uniform meshes: $640\times320$.}
		\label{sec3:HM2000}
	\end{figure}
		Finally, we solve the Mach 2000 astrophysical jet problem without a radiative cooling studied in \cite{ha2005numerical, ha2008positive, gardner2009numerical}. The computational domain is $[0,1]\times[-0.25,0.25]$. Initially, it is full of an ambient gas with $(\rho,\mu,\nu,p,\gamma)=(0.5,0,0,0.4127,\frac53)$. Outflow boundary conditions are imposed on the right, top, and bottom. The left boundary conditions are established with the values $(\rho, \mu, \nu, p, \gamma) = (5, 800, 0, 0.4127, \frac{5}{3})$ when $|y|<0.05$. For values outside of this range, the values are $(0.5, 0, 0, 0.4127, \frac{5}{3})$.
		In Figure \ref{sec3:HM2000}, we present the results obtained by the HWENO-U and HWENO-M schemes with PP limiters for a final time $T=0.001$.
The results show that the HWENO-U and HWENO-M schemes have comparable results, which are also similar to that in the reference \cite{ZhangS}.
\end{example}
\section{Concluding remarks}\label{Sec5:conclusion}

In this paper, we introduced a moment-based finite volume HWENO-U scheme with unified stencils on structured meshes. The {novel point} is to incorporate the spatial reconstructions with the modification of the first-order moments into a single step, resulting in a simpler approach {than} the HWENO-M scheme \cite{ZhaoZQiuJX_ArtificalLinearWeight_HWENO2020}, which involves two separate steps. The HWENO modification for the first-order moments in time discretizations serves two significant purposes in the proposed  scheme. Firstly, it ensures the stability of the fully-discrete scheme as that in the Lax-Friedrichs scheme \cite{LaxFridrichs1954}, which is demonstrated through analyses in Subsection \ref{sec2:stabilityHWENO}. Secondly, it helps to overcome spurious  oscillations for using nonlinear HWENO procedures. Furthermore, the proposed scale-invariant nonlinear weight of this paper {not only retains all properties of original one but also is more robust when simulating challenging problems with sharp scale variations}, shown in Examples \ref{Example:Lax1D} and \ref{Example:Sedov2D}.

Overall, the HWENO-U scheme is a simpler and more practical numerical method, which not only inherits the advantages of previous HWENO schemes, including compact stencils, high order accuracy, high resolution, and the use of artificial linear weights, but also employs unified stencils throughout the entire process without any modifications for the governing equations, {resulting in easier and faster implementations} as evidenced in the algorithm descriptions and numerical results. Furthermore, in the two-dimensional case, the framework of the HWENO-U scheme is based on truly two-dimensional reconstructions, making it more straightforward to extend to unstructured meshes, and the relevant works  are ongoing.

\begin{appendix}

\section{Appendix}\label{sec:Appendix1}	

In the one-dimensional case, the coefficients of the reconstructed polynomials $\{p_m(x)\}^2_{m=0}$ in \eqref{sec2:beta1d} are given as follows:
\begin{flalign}
	\hspace{3mm}&
	\begin{cases}
		c_{0,0}= -{\frac {43\,\bar{u}_{i-1}}{384}}+{\frac {235\,\bar{u}_{i}}{192}}-{\frac {43\,\bar{u}_{i+1}}{384}}-{\frac {27\,\bar{v}_{i-1}}{64}}+{\frac {27\,\bar{v}_{i+1}}{64}},\\
		c_{0,1}= -{\frac {63\,\bar{u}_{i-1}}{76}}+{\frac {63\,\bar{u}_{i+1}}{76}}-{\frac {75\,\bar{v}_{i-1}}{19}}-{\frac {75\,\bar{v}_{i+1}}{19}},\\
		c_{0,2}=  {\frac {23\,\bar{u}_{i-1}}{16}}-{\frac {23\,\bar{u}_{i}}{8}}+{\frac {23\,\bar{u}_{i+1}}{16}}+{\frac {45\,\bar{v}_{i-1}}{8}}-{\frac {45\,\bar{v}_{i+1}}{8}},\\
		c_{0,3}= {\frac {5\,\bar{u}_{i-1}}{19}}-{\frac {5\,\bar{u}_{i+1}}{19}}+{\frac {60\,\bar{v}_{i-1}}{19}}+{\frac {60\,\bar{v}_{i+1}}{19}},\\
		c_{0,4}= -\frac{5\,\bar{u}_{i-1}}{8}+\frac{5\,\bar{u}_{i}}{4}-\frac{5\,\bar{u}_{i+1}}{8}-{\frac {15\,\bar{v}_{i-1}}{4}}+{\frac {15\,\bar{v}_{i+1}}{4}};\\
		c_{1,0}= \bar{u}_{i},\ c_{1,1}= \bar{u}_{i}-\bar{u}_{i-1};\\
		c_{2,0}= \bar{u}_{i},\ c_{2,1}= \bar{u}_{i+1}-\bar{u}_{i}.
	\end{cases}	
	&  \nonumber %
\end{flalign}

In the two-dimensional case, the coefficients of the reconstructed polynomials $\{p_m(x,y)\}^4_{m=0}$ in \eqref{sec2:beta2d} are given as follows:
\begin{flalign}
	\hspace{3mm} &	
	\begin{cases}
		c_{0,0} = {\frac {{\bar{{u}}_1}}{576}}-{\frac {133\,{\bar{{u}}_2}}{1152}}+{\frac {{\bar{{u}}_3}
			}{576}}-{\frac {133\,{\bar{{u}}_4}}{1152}}+{\frac {419\,{\bar{{u}}_5}}{288}}-{
			\frac {133\,{\bar{{u}}_6}}{1152}}+{\frac {{\bar{{u}}_7}}{576}}-{\frac {133\,{
					\bar{{u}}_8}}{1152}}+{\frac {{\bar{{u}}_9}}{576}}-{\frac {27\,{\bar{{v}}_4}}{64}}+{
			\frac {27\,{\bar{{v}}_6}}{64}}-{\frac {27\,{\bar{{w}}_2}}{64}}+{\frac {27\,{\bar{{w}}_8}}{64}},
		\\
		c_{0,1} =\frac {\bar{{u}}_1}{48}-\frac{\bar{{u}}_3}{48}-{\frac {397\,{\bar{{u}}_4}}{456}}+{\frac {397\,{
					\bar{{u}}_6}}{456}}+\frac{\bar{{u}}_7}{48}-\frac{\bar{{u}}_9}{48}-{\frac {75\,{\bar{{v}}_4}}{19}}-{
			\frac {75\,{\bar{{v}}_6}}{19}} ,
		\\
		c_{0,2} = \frac{\bar{{u}}_1}{48} -  \frac {397\,{\bar{{u}}_2}}{456} +\frac{\bar{{u}}_3}{48}-\frac{\bar{{u}}_7}{48}+{
			\frac {397\,{\bar{{u}}_8}}{456}}-\frac{\bar{{u}}_9}{48}-{\frac {75\,{\bar{{w}}_2}}{19}}-{
			\frac {75\,{\bar{{w}}_8}}{19}}
		,
		\\
		c_{0,3} = -\frac{\bar{{u}}_1}{48}+\frac{\bar{{u}}_2}{24}-\frac{\bar{{u}}_3}{48}+{\frac {71\,{\bar{{u}}_4}}{48}}-{
			\frac {71\,{\bar{{u}}_5}}{24}}+{\frac {71\,{\bar{{u}}_6}}{48}}-\frac{\bar{{u}}_7}{48}+\frac{\bar{{u}}_8}{24} -\frac{\bar{{u}}_9}{48}+{\frac {45\,{\bar{{v}}_4}}{8}}-{\frac {45\,{\bar{{v}}_6}}{8}}
		,
		\\
		c_{0,4} = {\frac {7\,{\bar{{u}}_1}}{22}}+{\frac {7\,{\bar{{u}}_3}}{22}}+{\frac {7\,{\bar{{u}}_7}}{22}}-{\frac {7\,{\bar{{u}}_9}}{22}}-{\frac {75\,{\bar{{v}}_2}}{11}}+{\frac
			{75\,{\bar{{v}}_8}}{11}}-{\frac {75\,{\bar{{w}}_4}}{11}}+{\frac {75\,{\bar{{w}}_6}}{
				11}}
		,
		\\
		c_{0,5} = -\frac{\bar{{u}}_1}{48}+{\frac {71\,{\bar{{u}}_2}}{48}}-\frac{\bar{{u}}_3}{48}+\frac{\bar{{u}}_4}{24}-{
			\frac {71\,{\bar{{u}}_5}}{24}}+\frac{\bar{{u}}_6}{24}-\frac{\bar{{u}}_7}{48}+{\frac {71\,{\bar{{u}}_8}
			}{48}}-\frac{\bar{{u}}_9}{48}+{\frac {45\,{\bar{{w}}_2}}{8}}-{\frac {45\,{\bar{{w}}_8}}{8}}
		,
		\\
		c_{0,6} = {\frac {5\,{\bar{{u}}_4}}{19}}-{\frac {5\,{\bar{{u}}_6}}{19}}+{\frac {60\,{\bar{{v}}_4}}{19}}+{\frac {60\,{\bar{{v}}_6}}{19}}
		,
		\\
		c_{0,7} =  -\frac{\bar{{u}}_1}{4}+\frac{\bar{{u}}_2}{2}-\frac{\bar{{u}}_3}{4}+\frac{\bar{{u}}_7}{4}-\frac{\bar{{u}}_8}{2}+\frac{\bar{{u}}_9}{4}
		,
		\\
		c_{0,8} = -\frac{\bar{{u}}_1}{4}+\frac{\bar{{u}}_3}{4}+\frac{\bar{{u}}_4}{2}-\frac{\bar{{u}}_6}{2}-\frac{\bar{{u}}_7}{4}+\frac{\bar{{u}}_9}{4}
		,
		\\
		c_{0,9} = {\frac {5\,{\bar{{u}}_2}}{19}}-{\frac {5\,{\bar{{u}}_8}}{19}}+{\frac {60\,{\bar{{w}}_2}}{19}}+{\frac {60\,{\bar{{w}}_8}}{19}}
		,
		\\
		c_{0,10} = -\frac{5\,{\bar{{u}}_4}}{8}+\frac{5\,{\bar{{u}}_5}}{4}-\frac{5\,{\bar{{u}}_6}}{8}-{\frac {15\,{\bar{{v}}_4}}{4}}+{
			\frac {15\,{\bar{{v}}_6}}{4}}
		,
		\\
		c_{0,11} = {\frac {5\,{\bar{{u}}_1}}{22}}-{\frac {5\,{\bar{{u}}_3}}{22}}-{\frac {5\,{\bar{{u}}_7
			}}{22}}+{\frac {5\,{\bar{{u}}_9}}{22}}+{\frac {60\,{\bar{{v}}_2}}{11}}-{\frac {
				60\,{\bar{{v}}_8}}{11}}
		,
		\\
		c_{0,12} = \frac{\bar{{u}}_1}{4}-\frac{\bar{{u}}_2}{2}+\frac{\bar{{u}}_3}{4}-\frac{\bar{{u}}_4}{2}+{\bar{{u}}_5}-\frac{\bar{{u}}_6}{2}+\frac{\bar{{u}}_7}{4}-\frac{\bar{{u}}_8}{2}+\frac{\bar{{u}}_9}{4}
		,
		\\
		c_{0,13} ={\frac {5\,{\bar{{u}}_1}}{22}}-{\frac {5\,{\bar{{u}}_3}}{22}}-{\frac {5\,{\bar{{u}}_7
			}}{22}}+{\frac {5\,{\bar{{u}}_9}}{22}}+{\frac {60\,{\bar{{w}}_4}}{11}}-{\frac {
				60\,{\bar{{w}}_6}}{11}}
		,
		\\
		c_{0,14} = -\frac{5\,{\bar{{u}}_2}}{8}+\frac{5\,{\bar{{u}}_5}}{4}-\frac{5\,{\bar{{u}}_8}}{8}-{\frac {15\,{\bar{{w}}_2}}{4}}+{
			\frac {15\,{\bar{{w}}_8}}{4}};
		\\
c_{1,0}={\bar{{u}}_5}, \ c_{1,1}={\bar{{u}}_5}-{\bar{{u}}_4}, \ c_{1,2}={\bar{{u}}_5}-{\bar{{u}}_2};\\
c_{2,0}={\bar{{u}}_5},\ c_{2,1}={\bar{{u}}_6}-{\bar{{u}}_5},\ c_{2,2}={\bar{{u}}_5}-{\bar{{u}}_2};\\
c_{3,0}={\bar{{u}}_5},\ c_{3,1}={\bar{{u}}_5}-{\bar{{u}}_4},\ c_{3,2}={\bar{{u}}_8}-{\bar{{u}}_5};\\
c_{4,0}={\bar{{u}}_5},\ c_{4,1}={\bar{{u}}_6}-{\bar{{u}}_5},\ c_{4,2}={\bar{{u}}_8}-{\bar{{u}}_5}.\\
	\end{cases}
	\nonumber
\end{flalign}

\end{appendix}


\end{document}